\documentclass[preprint,12pt]{elsarticle}

\usepackage{lineno}
\usepackage[T1]{fontenc}
\usepackage[unicode]{hyperref}
\usepackage[margin=1in]{geometry}
\usepackage{amsmath}
\usepackage{amssymb}
\usepackage{booktabs}
\usepackage{multirow}
\usepackage{float}
\usepackage{subcaption}
\usepackage{bm}
\usepackage{makecell}

\usepackage[final]{changes}
\NewCommandCopy{\addedRI}{\added}
\NewCommandCopy{\replacedRI}{\replaced}
\NewCommandCopy{\deletedRI}{\deleted}
\makeatletter
\renewcommand{\addedRI}[1]{%
\@namedef{Changes@AuthorColor}{red}%
\colorlet{Changes@Color}{red}%
\added{#1}%
}
\renewcommand{\replacedRI}[1]{%
\@namedef{Changes@AuthorColor}{red}%
\colorlet{Changes@Color}{red}%
\replaced{#1}%
}
\renewcommand{\deletedRI}[1]{%
\@namedef{Changes@AuthorColor}{red}%
\colorlet{Changes@Color}{red}%
\deleted{#1}%
}
\makeatother
\NewCommandCopy{\addedRII}{\added}
\NewCommandCopy{\replacedRII}{\replaced}
\NewCommandCopy{\deletedRII}{\deleted}
\makeatletter
\renewcommand{\addedRII}[1]{%
\@namedef{Changes@AuthorColor}{teal}%
\colorlet{Changes@Color}{teal}%
\added{#1}%
}
\renewcommand{\replacedRII}[1]{%
\@namedef{Changes@AuthorColor}{teal}%
\colorlet{Changes@Color}{teal}%
\replaced{#1}%
}
\renewcommand{\deletedRII}[1]{%
\@namedef{Changes@AuthorColor}{teal}%
\colorlet{Changes@Color}{teal}%
\deleted{#1}%
}
\makeatother
\NewCommandCopy{\addedRIII}{\added}
\NewCommandCopy{\replacedRIII}{\replaced}
\NewCommandCopy{\deletedRIII}{\deleted}
\makeatletter
\renewcommand{\addedRIII}[1]{%
\@namedef{Changes@AuthorColor}{blue}%
\colorlet{Changes@Color}{blue}%
\added{#1}%
}
\renewcommand{\replacedRIII}[1]{%
\@namedef{Changes@AuthorColor}{blue}%
\colorlet{Changes@Color}{blue}%
\replaced{#1}%
}
\renewcommand{\deletedRIII}[1]{%
\@namedef{Changes@AuthorColor}{blue}%
\colorlet{Changes@Color}{blue}%
\deleted{#1}%
}
\makeatother
\begin{document}

\begin{frontmatter}
\journal{Journal}
\title{Computationally Efficient Optimisation of Elbow-Type Draft Tube Using Neural Network Surrogates}

\author[inst1,inst2]{Ante Sikirica}
\ead{ante.sikirica@uniri.hr}

\author[inst2]{Ivana Lu\v{c}in}
\ead{ivana.lucin@riteh.hr}

\author[inst2]{Marta Alvir}
\ead{marta.alvir@riteh.hr}

\author[inst1,inst2]{Lado Kranj\v{c}evi\'c\texorpdfstring{\corref{cor1}}{}}
\ead{lado.kranjcevic@riteh.hr}

\author[inst2]{Zoran \v{C}arija}
\ead{zoran.carija@riteh.hr}

\cortext[cor1]{Corresponding author}

\affiliation[inst1]{organization={Center for Advanced Computing and Modelling, University of Rijeka},
            addressline={Radmile Matejčić 2}, 
            city={Rijeka},
            postcode={51000 Rijeka}, 
            country={Croatia}}

\affiliation[inst2]{organization={Faculty of Engineering, University of Rijeka},
            addressline={Vukovarska 58}, 
            city={Rijeka},
            postcode={51000 Rijeka}, 
            country={Croatia}}

\begin{abstract}
This study aims to provide a comprehensive assessment of single-objective and multi-objective optimisation algorithms for the design of an elbow-type draft tube, as well as \addedRI{to} introduce a computationally efficient optimisation workflow. The proposed workflow leverages deep neural network surrogates trained on \replacedRI{data obtained from numerical simulations.}{computational fluid dynamics data.} \deletedRI{defined by Latin hypercube sampling. Surrogate models are used to predict the performance of draft tube designs without the need for computationally demanding numerical simulations.} \addedRI{The use of surrogates allows for a more flexible and faster evaluation of novel designs.}
\addedRI{The success history-based adaptive differential evolution with linear reduction and the multi-objective evolutionary algorithm based on decomposition were identified as the best-performing algorithms and used to determine the influence of different objectives in the single-objective optimisation and their combined impact on the draft tube design in the multi-objective optimisation.} \addedRI{The results for the single-objective algorithm are consistent with those of the multi-objective algorithm when the objectives are considered separately. Multi-objective approach, however, should typically be chosen, especially for computationally inexpensive surrogates.} A multi-criteria decision analysis method was used to obtain optimal multi-objective results, showing an improvement of 1.5\% and 17\% for the pressure recovery factor and drag coefficient, respectively. The difference between the predictions and the numerical results is less than 0.5\% for the pressure recovery factor and 3\% for the drag coefficient.
\deletedRI{Multi-objective evolutionary algorithm based on decomposition and success history-based adaptive differential evolution with linear reduction were identified as the best-performing algorithms and were used to determine the influence of the objectives in the single-objective optimisation and the combined effects on the draft tube design in the multi-objective optimisation.}
\addedRI{As the demand for renewable energy continues to increase, the relevance of data-driven optimisation workflows, as discussed in this study, will become increasingly important, especially in the context of global sustainability efforts.}
\end{abstract}



\end{frontmatter}

\section{Introduction}

The growing interest in renewable energy sources is a major driving force behind hydropower research and development aimed at improving efficiency. An essential instrument in this effort is the design optimisation technique. Typically, design optimisation is employed to improve runner designs, often utilising novel mathematical and computational methods \cite{wang2011, yang2014}. However, this strategy can incur significant deployment costs for an existing hydropower turbine. While this approach may be justifiable in the context of a long-term revitalisation of a turbine, in the short term, optimising specific segments where simpler redesigns are needed may be more cost-effective.

Draft tubes play a crucial role in the overall performance of hydropower turbines \addedRI{\cite{warnick1984}}. These components act as conduits for the water leaving the turbine and influence the flow characteristics and pressure recovery \addedRI{\cite{warnick1984}}. \deletedRI{The notion of optimising draft tube design during maintenance operations offers a path to improve the overall efficiency of the turbine within a short time frame and at a minimal cost.} \addedRI{As the water from the runner enters the draft tube, its swirling motion can form a structure known as a vortex rope. These vortex structures induce flow instabilities and can cause pressure fluctuations, vibrations and critical damage. In a review paper by \citet{kumar2021}, this was investigated from an experimental perspective. Recently, \citet{zhou2023} proposed a draft tube design methodology to mitigate vortex rope creation, thus avoiding pressure fluctuations.}

Model testing and optimisation are typically used to ensure the effectiveness of the draft tube. Traditionally, this process was based on engineering experience and could be time-consuming and susceptible to oversights \cite{puente2003}. \addedRI{Furthermore, turbines and specific components were designed to work optimally at a single point, but the trend has shifted towards turbines that can work at various operating points and conditions \cite{marjavaara2007, lyutov2015}. New trends for flexible operation require new strategies to find the best design. Consequently, optimisation methods and guidelines that are reliable and computationally efficient are needed. \citet{abbas2017} and \citet{tiwari2020} conducted reviews of different computational fluid dynamics (CFD)-based approaches for design optimisation.}

The conventional design cycle involves initial design selection, evaluation, and iterative redesign to attain desired characteristics. \deletedRI{Cross-sections, characterised by a set of parameters, are essential in the design process \mbox{\cite{marjavaara2003, sosa2015, arispe2018}}. The trend has recently shifted towards encompassing a more comprehensive range of working points and operating conditions \mbox{\cite{marjavaara2007, lyutov2015}}. Introduced complexity has incentivised the use of different parameterisation techniques to reduce the number of variables and increase the efficiency of the design process \cite{marjavaara2007}.} \addedRI{Over the years, various parameterisation techniques have been used to reduce the complexity of the optimisation problem. The fundamental goal is to reduce the number of optimisation variables and increase the efficiency of the design process \cite{marjavaara2007}. The main problem is that parameterisation dictates the shapes that can be obtained, so innovative designs are difficult to obtain. The use of cross-sections can provide some flexibility and has been explored in several papers \cite{marjavaara2003, sosa2015, arispe2018}. However, a large number of cross-sections are required for a suitable geometry representation, so the number of optimisation variables increases rapidly and, with it, the complexity of the problem. In addition, such an approach can lead to shapes that are difficult to implement \cite{lucin2022}. Therefore, an approach based on parametric curves may be more advantageous as it allows geometric flexibility but does not significantly increase the number of optimisation variables.} 

\deletedRI{Moreover,} The use of new and often complicated optimisation frameworks is necessary \cite{wang2011, mcnabb2014, daniels2021}\deletedRI{\cite{lucin2022}}, as the dynamic interaction between different parameters and operating conditions requires a comprehensive approach. \deletedRI{Consequently, optimisation methods and guidelines that meet robustness requirements and are computationally efficient are needed. Overviews of different computational fluid dynamics (CFD) based approaches for design optimisation were conducted by \mbox{\citet{abbas2017}} and \mbox{\citet{tiwari2020}}. In addition, it is important to consider the tendency for vortex rope formation during the draft tube design process. A review paper by \mbox{\citet{kumar2021}} examined this aspect from an experimental point of view. Recently, \mbox{\citet{zhou2023}} proposed a draft tube design methodology to address vortex rope and pressure fluctuations.} In a study by \citet{marjavaara2006}, draft tube redesign was attained by combining the design of experiments (DOE) technique, the response surface method (RSM)\deletedRII{,} and \deletedRIII{the} single-objective optimisation (SO). \addedRIII{The} surrogate models were based on 36 CFD simulations and \addedRIII{were} individually optimised to maximise \addedRIII{the} pressure recovery factor\addedRI{,} $C_p$\addedRI{,} and minimise \addedRIII{the} energy loss\addedRI{,} $\zeta$. \deletedRI{Similarly,}\citet{demirel2017} combined the RSM\addedRI{, based on 50 design points,} with multi-objective (MO) optimisation. An improvement of 4.3\% and 20\% was reported for the pressure recovery factor \deletedRIII{$C_p$} and head loss, $\Delta H$, respectively, but with a model error of 8\%. \replacedRI{Study by}{The paper by} \citet{marjavaara2007} employed the multi-objective non-dominated sorting genetic algorithm II (NSGA-II) in conjunction with local search to generate a design that meets pressure recovery requirements under different operating conditions. \addedRI{A total of 34 design points were used and a simplified draft tube was used for testing.} The authors noted that radial basis neural networks (RBNN) performed better as surrogates than the RSM approach. Similarly, \citet{mcnabb2014} used RBNN but relied on a hierarchical evolutionary algorithm (EA) that combines different construction and production costs as objectives using a weighted sum approach. \citet{shojaeefard2014} evaluated the applicability of artificial neural networks (ANN) coupled with NSGA-II using 64 design points to train the ANN model with two input features and two outputs, namely pressure recovery and energy loss. The surrogate achieved excellent accuracy, in part due to the simplified definition of the optimisation problem. \deletedRI{ \mbox{\citet{lyutov2015}} proposed a simultaneous MO optimisation of the turbine runner and the draft tube, with segments described by 28 and 9 geometric parameters, respectively. Conducted optimisations provided an average improvement of 0.3\% when optimising coupled problems compared to just runner optimisation. Although the improvement is quantifiable, it is questionable whether this approach has merit if computational effort and simplifications are taken into account.} \addedRI{In all noted studies, a small number of design points were used. Additionally, the geometries used for the optimisation process were mostly simplified or considerably constrained due to parametrisation. To address these limitations, \citet{hammond2022} have recently investigated the applicability of various machine learning algorithms for turbine design. Methods and guidelines that can improve current design approaches and increase optimisation performance are described.}

\addedRI{The multi-objective optimisation of draft tubes has been assessed in multiple studies \cite{demirel2017,marjavaara2007,shojaeefard2014}. However, clear guidelines were never established, and different objective functions were often used.} Recently, \citet{orso2020} performed an MO optimisation of a draft tube, with the pressure recovery factor \deletedRIII{$C_p$} and the drag coefficient\addedRI{,} $C_d$\addedRI{,} as objectives. \addedRI{The drag coefficient is effectively equal to the energy loss in this context.} The RSM model was trained on 150 \deletedRIII{DOE} samples defined by nine characteristic variables\replacedRIII{. This resulted}{, resulting} in improvements of 2.54\% and 5.7\% for $C_p$ and $C_d$, respectively, with no reported surrogate errors. \citet{daniels2021} employed a MO Bayesian methodology and reported improvements of 3.7\% and 22.3\% for \replacedRI{pressure recovery}{$C_p$} and \replacedRI{energy loss}{$\zeta$}, respectively. \citet{fleischli2021} compared the efficiency of EA and a general discrete adjoint method, finding that the latter, although more complex to implement and occasionally limited in flexibility, reduced computational time by 30\% while providing performance improvements. \replacedRIII{The authors}{They} concluded that local and global optimization approach\addedRIII{es} can be beneficial for multi-point\deletedRIII{operating} and multi-objective optimizations. \citet{sheikh2022} introduced an optimisation framework that combines Bayesian optimisation and a design-by-morphing concept, where the design is defined as a linear combination of transformations that represent cross-sectional shapes.

\deletedRI{This paper expands on existing methods for draft tube design by providing a comprehensive assessment of various single-objective and multi-objective algorithms.} Previous research on draft tube optimisation has been limited in terms of adaptability, cost-effectiveness assessment, and method diversity. Additionally, many studies \addedRIII{have} used their \deletedRIII{own} techniques, making it difficult to determine the most appropriate approach \addedRI{that can be used in practice. When the geometry is flexible, only a single surrogate and a single optimisation strategy were explored. On the other hand, when different numerical strategies were explored, considerable geometric simplifications were used.} This gap in the literature is now being filled. \deletedRI{The selection of stochastic optimisation algorithms includes the multi-objective evolutionary algorithm based on decomposition (MOEA/D), the strength Pareto evolutionary algorithm 2 (SPEA2) and NSGA-II as well as single-objective fireworks algorithm (FWA), particle swarm optimisation (PSO) and success history-based adaptive differential evolution with linear reduction (L-SHADE). The optimal results for the multi-objective algorithms are determined with the technique for order of preference by similarity to the ideal solution (TOPSIS). The main objective is to investigate the behaviour, characteristics, and practical utility of these algorithms.} \deletedRI{Additionally,}A simple, easy-to-understand, and validate \addedRIII{draft tube optimisation} workflow is \deletedRIII{also} proposed. This involves the use of deep neural network (DNN) surrogates trained on CFD data obtained through the use of Latin hypercube sampling (LHS). Machine learning in this context replaces the traditionally used RSM\deletedRI{and ensures improved accuracy}. \addedRI{Although they might require more computational effort for training, tuned DNNs can handle high-dimensional data effectively and are highly adaptable and able to learn from new data. An important aspect of this study is thus the assessment of the applicability of DNNs as an alternative to conventional methods. The optimisation algorithms are assessed in different contexts. The primary optimisation objectives that are evaluated are the pressure recovery factor and the drag coefficient. The selection of stochastic optimisation algorithms includes the multi-objective evolutionary algorithm based on decomposition (MOEA/D), the strength Pareto evolutionary algorithm 2 (SPEA2) and NSGA-II as well as single-objective fireworks algorithm (FWA), particle swarm optimisation (PSO) and success history-based adaptive differential evolution with linear reduction (L-SHADE). The optimal results for the multi-objective algorithms are determined with the technique for order of preference by similarity to the ideal solution (TOPSIS).} Based on the \deletedRI{obtained} results and insights, suggestions and guidelines for future research are given. The applicability for a wide range of operating points is also discussed. \addedRI{Clear guidance is provided on selecting the appropriate optimisation algorithm and objectives.}

\section{Materials and methods}

\subsection{Workflow}

A systematic workflow is envisioned \replacedRIII{to determine the most favourable draft tube design efficiently}{to efficiently determine the most favourable draft tube design}. The proposed workflow hinges on machine learning and data obtained using computational fluid dynamics \addedRIII{simulations}. The surrogate model is subsequently used as a rapid evaluator, predicting the performance of new draft tube designs without additional CFD \replacedRIII{evaluations}{simulations}. This significantly accelerates design exploration.

The first step is to determine near-random design candidates through Latin hypercube sampling. This approach ensures that a wide range of potential configurations is evaluated and a non-biased dataset is generated. Subsequently, a rigorous CFD assessment of these candidates is performed. Simulations are carried out in an automated and standardised manner, with several safeguards to ensure the validity of the results. Particular attention is given to the pressure data, which are used to calculate two coefficients, the pressure recovery factor $C_p$ and the drag coefficient $C_d$. \addedRIII{The pressure recovery factor is commonly used as a measure of the efficiency of the draft tube. Since the fluid decelerates with the cross-sectional area of the draft tube, the static pressure, and consequently the pressure recovery factor, increases. In this sense, the draft tube effectively recovers the kinetic energy of the runner and thus increases the net head of the turbine. The drag coefficient quantifies the energy that cannot be recuperated, i.e., the energy that is converted into an unusable form. This is mainly due to frictional losses, secondary losses, flow separation, etc. \cite{orso2020, daniels2020}.} The pressure recovery factor is calculated as follows \addedRI{\cite{white1990}}:
\begin{linenomath}
\begin{equation}
\label{cp_eq}
    C_p = \frac{p_{s,2} - p_{s,1}}{0.5 \cdot \rho u^2}~.
\end{equation}
\end{linenomath}
The drag coefficient is defined as \addedRI{\cite{white1990}}:
\begin{linenomath}
\begin{equation}
\label{cd_eq}
    C_d = \frac{p_{t,1}-p_{t,2}}{0.5 \cdot \rho u^2}~.
\end{equation}
\end{linenomath}
\replacedRI{In noted equations (\ref{cp_eq},~\ref{cd_eq}),}{Indices 1 and 2 in the equations refer to the inlet and outlet, respectively,} $p_{s}$ is the static pressure, $p_{t}$ is the total pressure, $\rho$ is the fluid density, and $u$ is the absolute velocity. \addedRI{The indices 1 and 2 used for pressure in the equations refer to the values at the inlet and outlet, respectively.}

\begin{figure}[H]
    \centering
    \includegraphics[width=\textwidth]{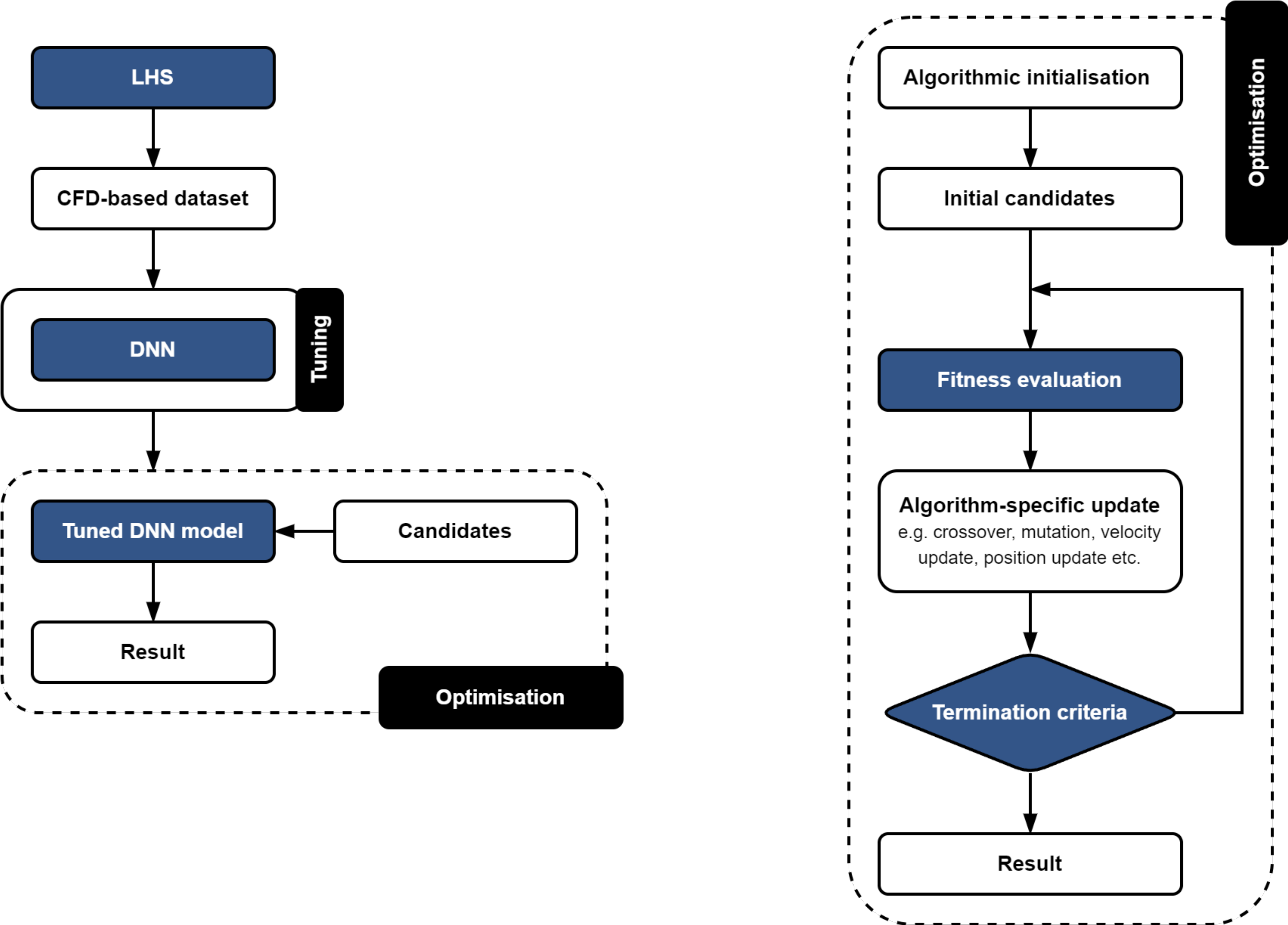}
    \caption{ \protect\addedRI{A proposed workflow combining LHS for geometry selecetion, CFD, DNN surrogate modelling and optimisation (left). CFD is used to generate data. The tuned DNNs are coupled with the SO and MO algorithms. The optimal MO solution is determined using the TOPSIS method. Typical optimisation process (right). First, the algorithm is initialised, and an initial set of candidates is generated. The fitness is calculated, and each algorithm performs algorithm-specific steps (crossover, mutation, location update, velocity and position updates, etc.). The process is repeated until the termination criterion is satisfied (number of iterations, generations).}}
    \label{fig:workflow}
\end{figure}

The numerical data collected serves as the basis for the datasets on which the DNN surrogate models are trained. DNNs are fine-tuned to improve their ability to make accurate predictions. The final step involves utilising the DNNs with different optimisation algorithms to determine the optimal draft tube design, considering the objectives and constraints set. Different SO and MO optimisation algorithms are compared. The entire workflow is implemented in Python \cite{python}. Figure \ref{fig:workflow} shows a simplified scheme of the proposed workflow.

\subsection{Geometry and parametric model}

\begin{figure}[H]
    \centering
    \includegraphics[trim={6cm 1.5cm 6cm 1.5cm}, clip, height=6cm]{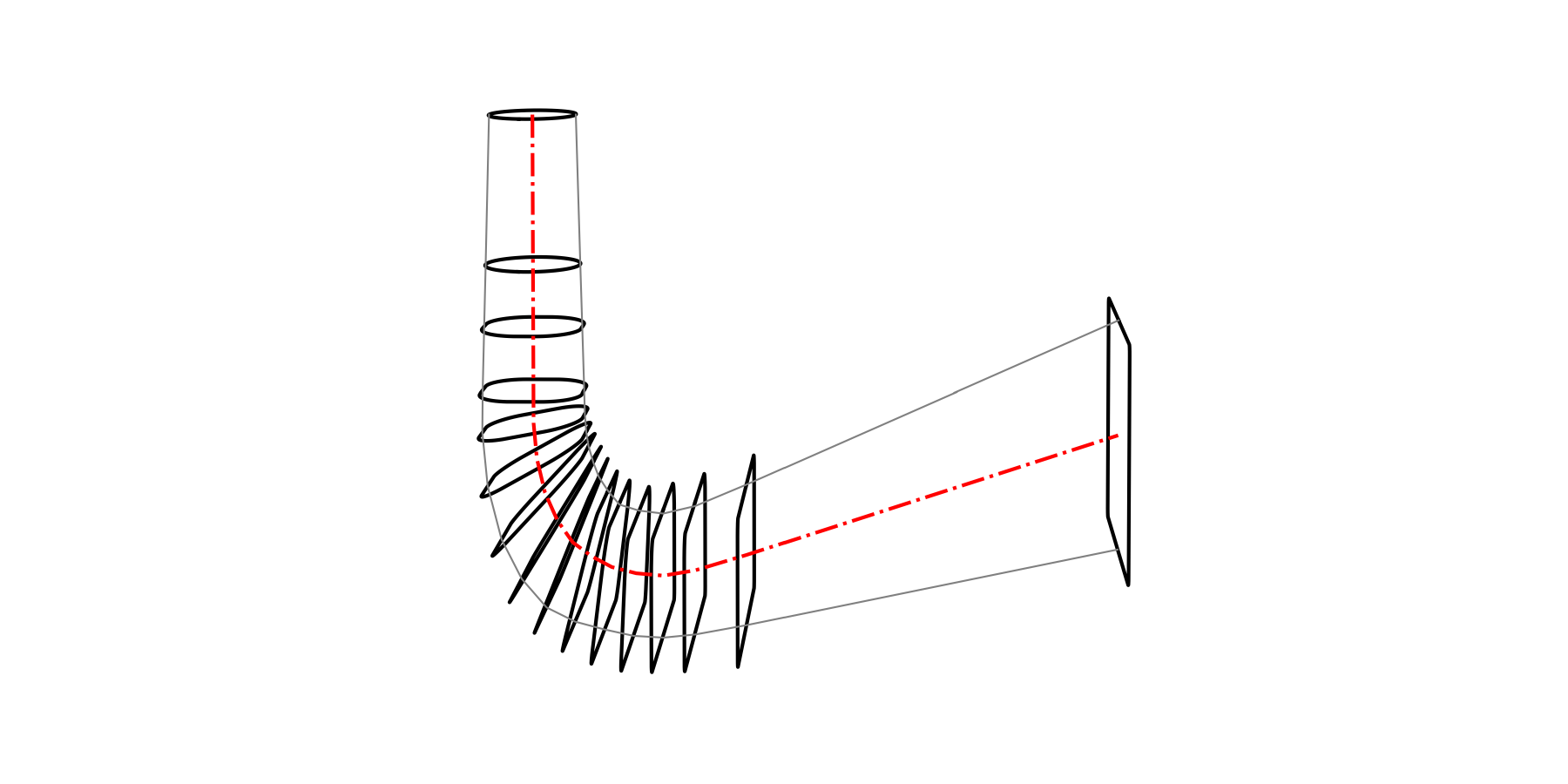}
    \caption{HPP Rijeka draft tube reference cross-sections.}
    \label{fig:model_3D}
\end{figure}

The Rijeka hydroelectric power plant (HPP), Rijeka, Croatia, is situated on the Rječina River and has two Francis turbines installed 1.6 m above sea level. At a flow rate \deletedRIII{of} $Q = 2 \cdot 10.5$ m$^3$/s, the HPP Rijeka generates 36.8 MW of electricity \addedRI{\cite{HERijeka}}. The baseline elbow-type draft tube design, which is part of the installed Francis turbines, is considered in this study. The design can be described by 16 reference cross-sections, as illustrated in Figure \ref{fig:model_3D}.

\begin{figure}[H]
\centering
\begin{subfigure}[b]{0.33\textwidth}
    \centering
    \includegraphics[trim={3cm 3cm 3cm 3cm}, clip, height=4cm]{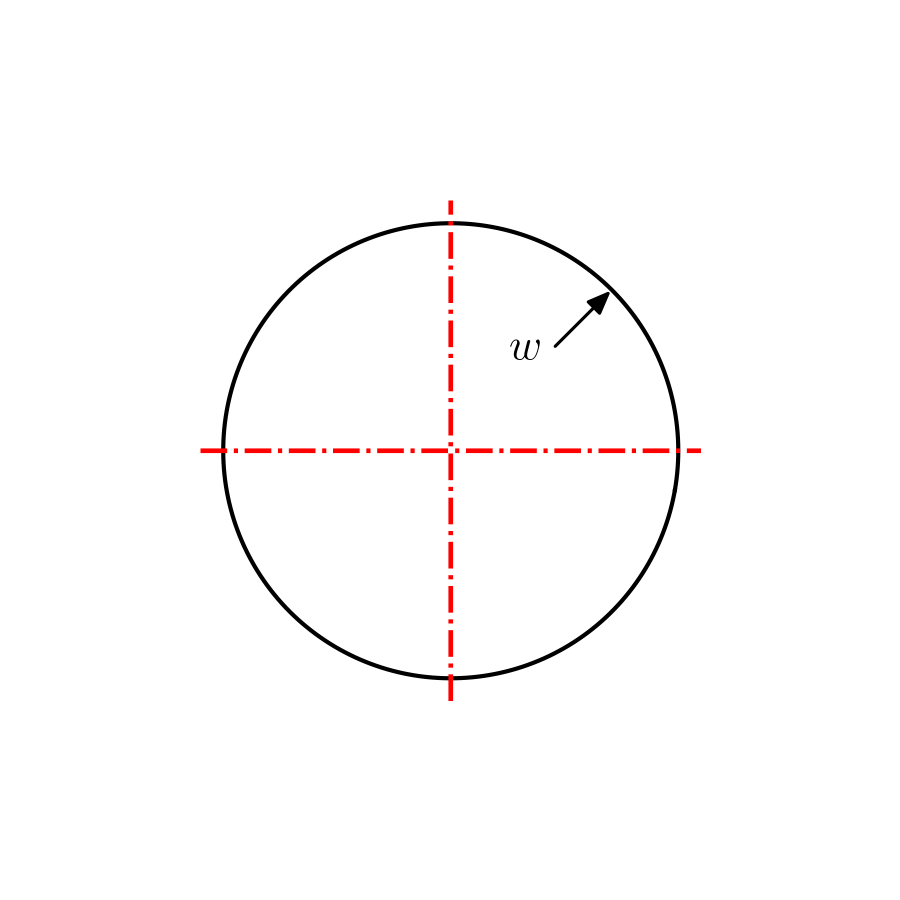}
    \caption{}
    \label{fig:model_circular_section}
\end{subfigure}
\begin{subfigure}[b]{0.66\textwidth}
    \centering
    \includegraphics[trim={3cm 3cm 3cm 3cm}, clip, height=4cm]{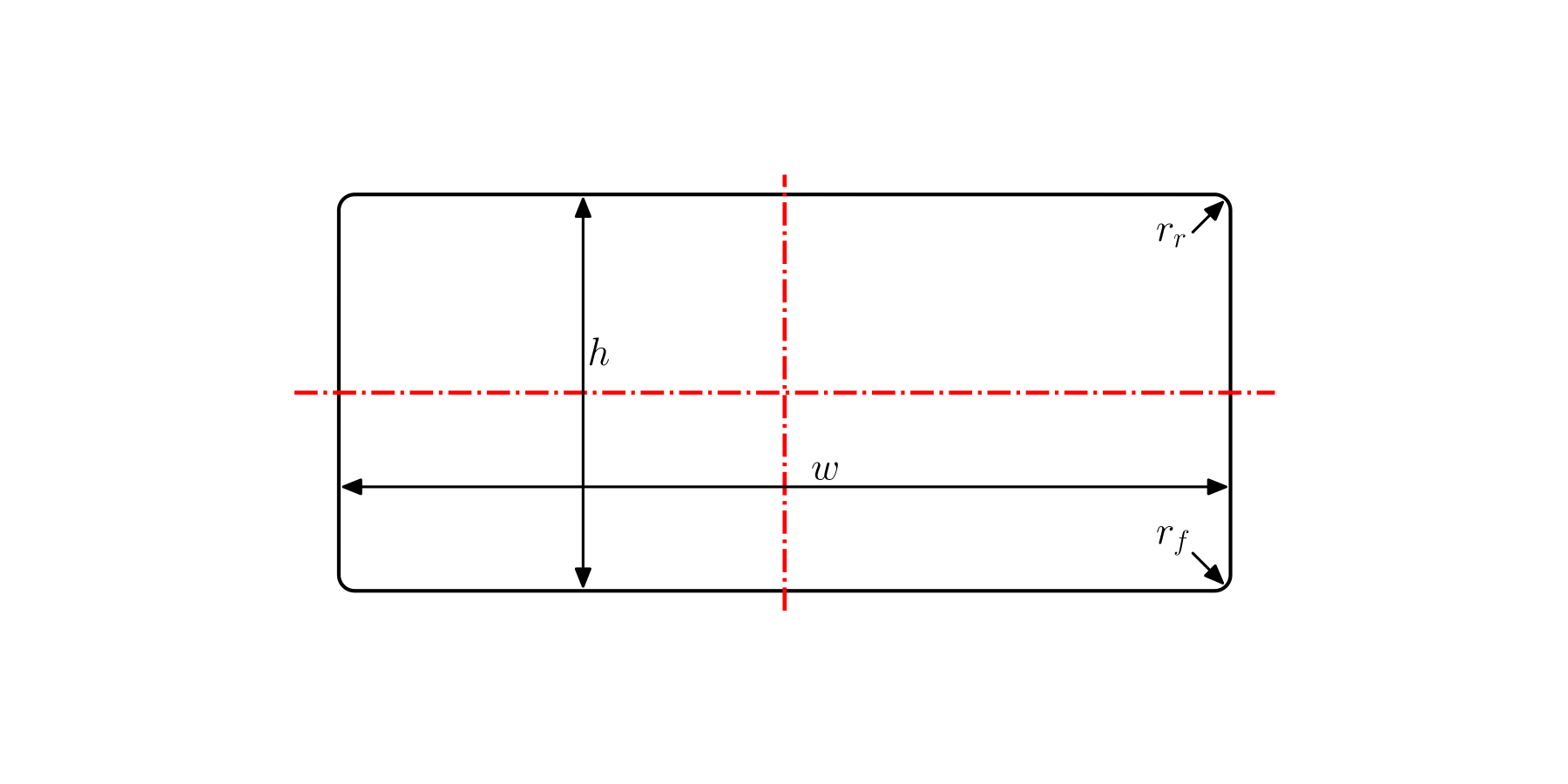}
    \caption{}
    \label{fig:model_rectangular_section}
\end{subfigure}
\caption{Circular (a) and rectangular (b) cross-sections of the elbow-type draft tube with characteristic dimensions.}
\label{fig:model_sections}
\end{figure}

The first three cross-sections are circular and gradually change into a rounded rectangle. They are defined by a single parameter, radius $w$. The subsequent ellipsoidal and rectangular cross-sections are defined by their width $w$ (horizontal semi-axis) and height $h$ (vertical semi-axis). The corners of the rectangular cross-sections are curved according to the radii of curvature $r_r$ and $r_f$ for the roof and floor of the draft tube, respectively. Each cross-section has a specific spatial position and angle (Figure \ref{fig:model_3D}). Typical cross-sections (first and last) are shown in Figure \ref{fig:model_sections}.

In order to simplify the design optimisation process and to ensure shape flexibility, a smoothed design described by 84 extrapolated cross-sections was created. The resulting shape was subsequently parameterised by fitting B-spline curves along the centreline through the draft tube's roof, floor and width. This means that for any point $x$ along the centreline, a corresponding cross-section can be determined from three fitted functions. The remaining parameters (i.e. roof and floor radii) are calculated by interpolating between adjacent cross-sections for a given $x$.

A B-spline curve $\bm{C}(t)$ is a sequence of polynomial curve segments. The analytical formulation of a B-spline can be written as:
\begin{linenomath}
\begin{equation}
    \bm{C}(t) = \sum_{i=0}^{n} \bm{P}_{i} N_{i,k} (t)~\addedRIII{,}
\end{equation}
\end{linenomath}
where $\{\bm{P}_i : i=0,1,...,n\}$ are the control points, $N_{i,k} (t)$ are normalized blending (basis) functions, $k$ is the order of the polynomial curve segments and \addedRI{$t_i$} a knot \addedRI{\cite{piegl1996}}. Knots form a non-decreasing sequence of real numbers, i.e. a knot vector $\{t_i : i=0,...,n+k\}$. Basis functions $N_{i,k} (t)$ \deletedRI{for $k > 1$} are \addedRI{calculated} as follows \addedRI{\cite{piegl1996}}:
\addedRI{
\begin{linenomath}
\begin{equation}
\label{eq_bspline_one}
N_{i,1} (t) = 
    \begin{cases}
    1,&  \text{if}~t_i \leq t < t_{i+1}~,\\
    0,&  \text{else}~,
    \end{cases}
\end{equation}
\begin{equation}
\label{eq_bspline_generic}
N_{i,k} (t) = \frac{t - t_i}{t_{i+k-1} - t_i} \cdot N_{i,k-1} (t) + \frac{t_{i+k} - t}{t_{i+k} - t_{i+1}} \cdot N_{i+1,k-1} (t)~.
\end{equation}
\end{linenomath}
The equation (\ref{eq_bspline_generic}) can be rewritten as:
\begin{linenomath}
\begin{equation}
N_{i,k} (t) = w_{i,k-1} (t) \cdot N_{i,k-1} (t) + [1 - w_{i+1,k-1} (t)] \cdot N_{i+1,k-1} (t)~,
\end{equation}
\end{linenomath}
where
\begin{linenomath}
\begin{equation}
w_{i,k} (t) =
    \begin{cases}
    \dfrac{t - t_i}{t_{i+k} - t_i},&  \text{if}~t_i < t_{i+k}~,\\
    0,&  \text{else}~.
    \end{cases}
\end{equation}
\end{linenomath}
}

\begin{figure}[H]
\centering
\begin{subfigure}[b]{0.66\textwidth}
    \centering
    \includegraphics[trim={1cm 1cm 1cm 1cm}, clip, height=5cm]{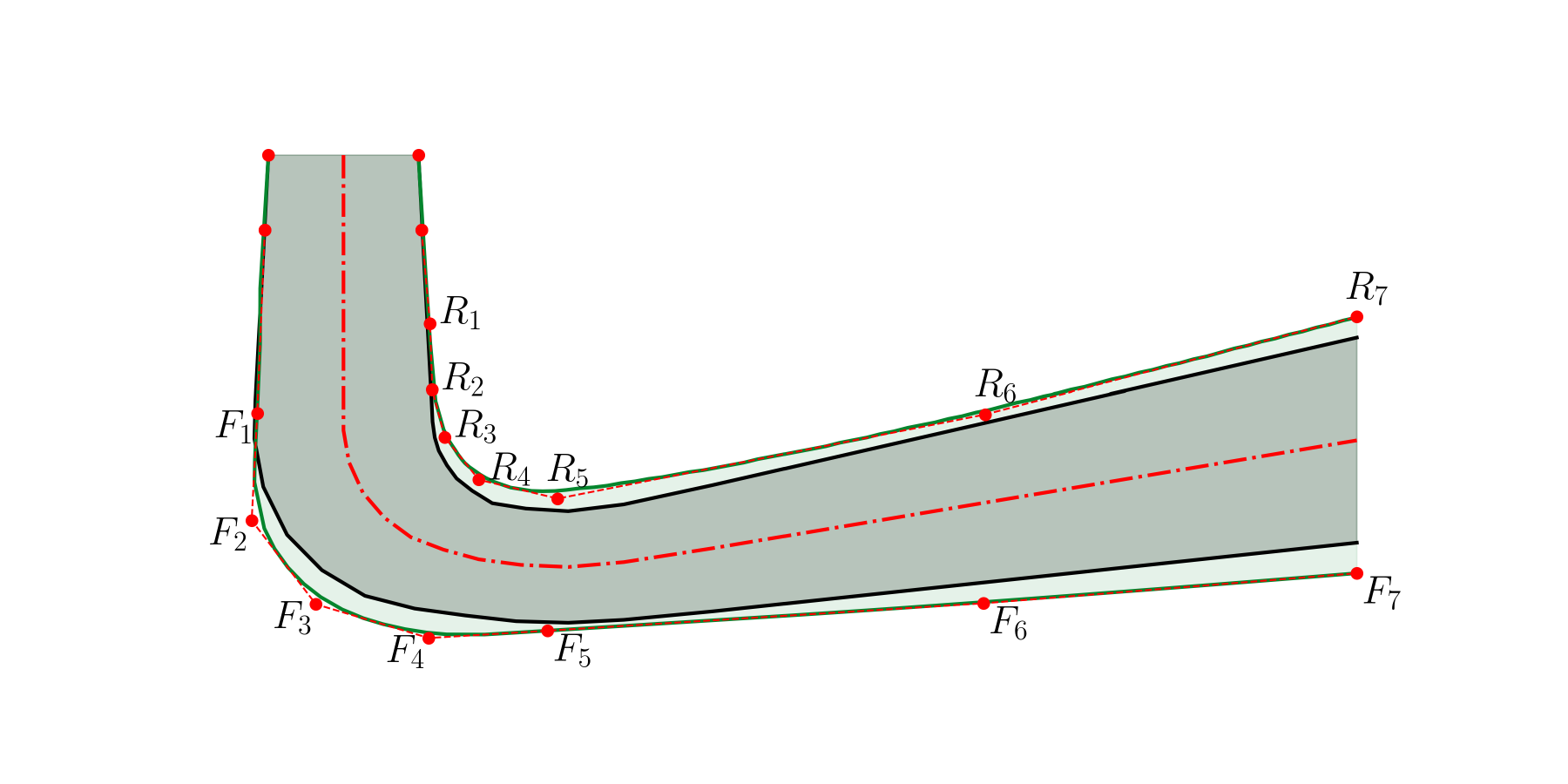}
    \caption{}
    \label{fig:model_control_points_tb}
\end{subfigure}
\begin{subfigure}[b]{0.33\textwidth}
    \centering
    \includegraphics[trim={1cm 1cm 1cm 1cm}, clip, height=5cm]{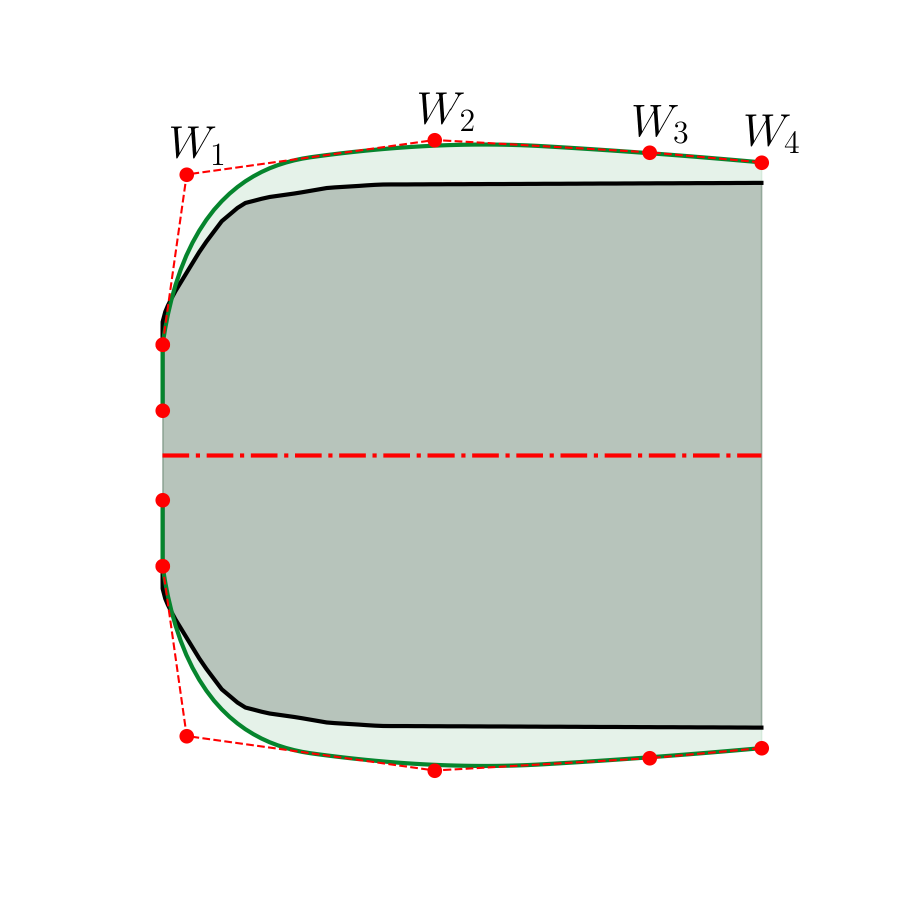}
    \caption{}
    \label{fig:model_control_points_w}
\end{subfigure}
\caption{Draft tube parameterisation. The reference design of the draft tube is outlined by black lines. The red dots are control points. Dashed lines represent control polygons. The roof and floor of the draft tube are described by B-spline curves with nine control points (a), while the width of the draft tube is parameterised by a B-spline curve with six control points (b).}
\label{fig:model_control_points}
\end{figure}

In this paper, the internal knots are evenly distributed (uniform B-spline). The B-spline is quadratic. The B-spline curves of the roof and floor are governed by nine control points each. The width-defining B-spline curve has six control points. The first two control points for all splines are kept constant. The parameterisation of the draft tube and the corresponding control points are given in Figure \ref{fig:model_control_points}.

\subsection{Numerical setup}

\begin{figure}[H]
    \centering
    \includegraphics[trim={0cm 0cm 0cm 0cm}, clip, height=5cm]{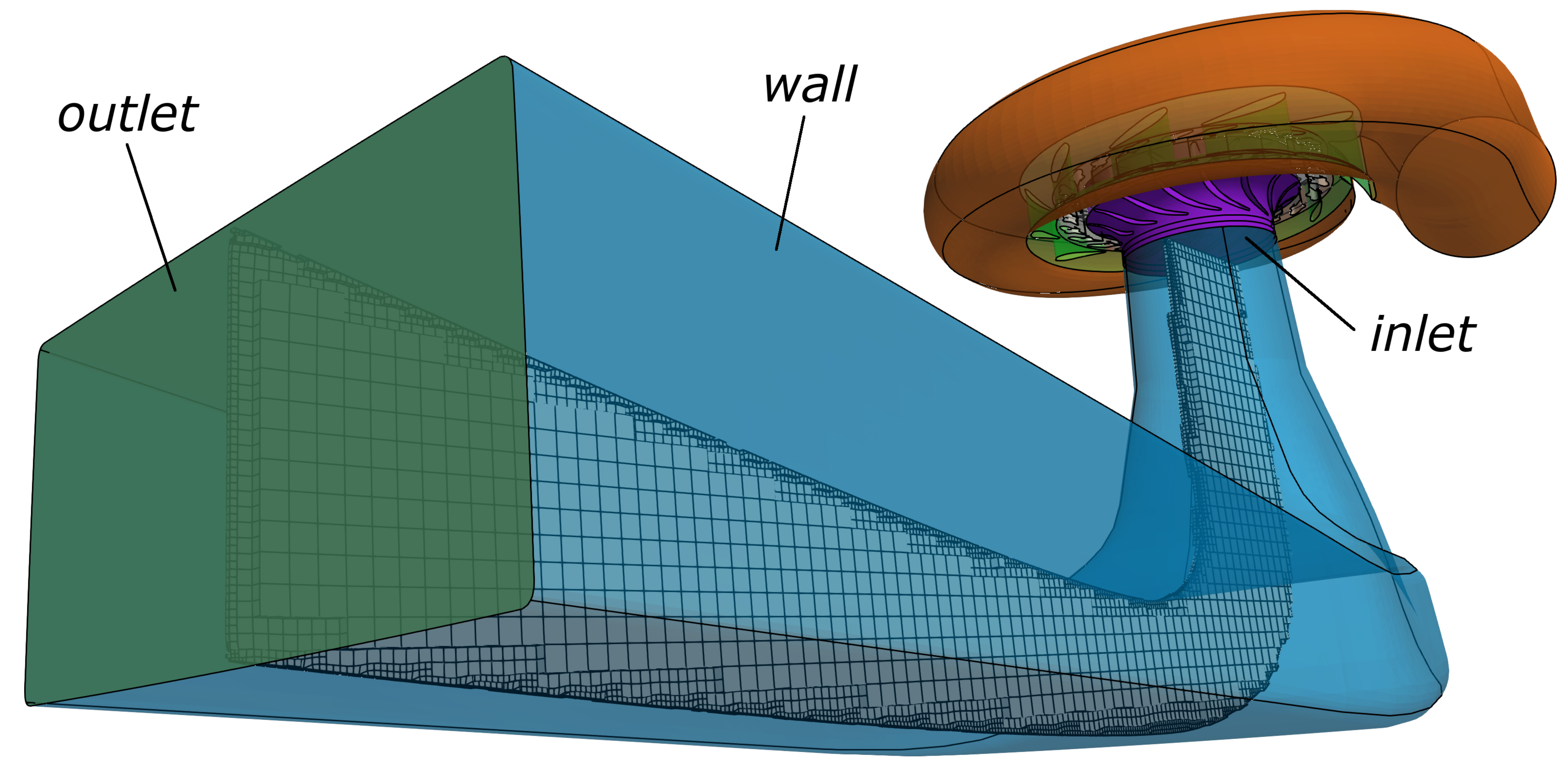}
    \caption{Computational domain and a sectional cut of the numerical grid. The runner and upstream segments are displayed for context.}
    \label{fig:domain}
\end{figure}

The computational domain is created through FreeCAD's Python API \cite{freecad}. \textit{snappyHexMesh} \cite{ofv2012} is used to generate the numerical grid. The maximum and minimum cell sizes are $l_{max} = f_s \cdot 0.064$ m and $l_{min} = f_s \cdot 0.001$ m, respectively, where $f_s$ is a scaling factor. The dimensionless wall distance, $y^+$, is on average greater than 70, which implies that the first grid point is in the log-law region, and hence wall functions must be used. Due to the variety of geometries generated in the evaluation process, $y^+$ is determined at the end of each simulation. Depending on the value, the results are discarded, and the simulation is rerun with updated grid parameters. A sectional cut of the numerical grid and computational domain in the context of the entire turbine is shown in Figure \ref{fig:domain}.

Relevant numerical data are calculated with the open source toolbox OpenFOAM \cite{ofv2012}. A steady-state incompressible solver, \textit{simpleFoam}, is used. The realisable k-epsilon turbulence model is employed. The fluid is water with a kinematic viscosity $\nu = 1 \cdot 10^{-6}$ m$^2$/s. The domain is bounded by three patches: inlet, outlet, and draft tube wall. A fixed flow rate $Q = 10$ m$^3$/s is set at the inlet instead of a velocity profile, which is a deliberate simplification that should not affect the validity of the proposed workflow. The turbulent viscosity $\nu_t$, the turbulent kinetic energy $k$ and the dissipation of the turbulent kinetic energy $\epsilon$ were estimated based on the domain size and the inflow velocity \cite{ferziger2019}. Table \ref{tab:boundary_conditions} provides a summary of the boundary conditions.

\begin{table}[H]
\centering
\footnotesize
\caption{A summary of employed boundary conditions for CFD simulations.}
\label{tab:boundary_conditions}
\begin{tabular}{@{}cccc@{}}
\toprule
Field & Inlet & Outlet & Wall \\ \midrule
\textbf{u} & \textit{flowRateInletVelocity} & \textit{inletOutlet} & \textit{noSlip} \\
$p$ & \textit{zeroGradient} & \textit{fixedValue} & \textit{zeroGradient} \\
$\nu_t$ & \textit{calculated} & \textit{calculated} & \textit{nutkWallFunction} \\
$k$ & \textit{fixedValue} & \textit{zeroGradient} & \textit{kqRWallFunction} \\
$\epsilon$ & \textit{fixedValue} & \textit{zeroGradient} & \textit{epsilonWallFunction} \\ \bottomrule
\end{tabular}
\end{table}

Second-order numerical schemes are used for convective terms; \textit{linearUpwind} and \textit{limitedLinear} schemes are used for \textbf{u} and the remaining fields, respectively. All remaining terms are second-order accurate. The simulations are initiated with first-order accurate schemes. The simulations run for $10^4$ iterations or until the convergence criteria of $1 \cdot 10^{-4}$ are met. To ensure stability, the standard deviation of the physical value (i.e. pressure) for the last 5\% iterations is calculated and considered acceptable if it is below 10\% of the mean.

\subsection{Design of experiments}

The design of experiments is a systematic statistical approach used to explore the relationship between independent and dependent variables by sampling a decision space. It is a powerful data collection and analysis paradigm that can be used in various experimental scenarios. Latin hypercube sampling is a statistical method used to explore the decision space efficiently, i.e., to generate near-random samples from the decision space $\mathbb{R}^m$ \cite{mckay2000}. For each variable $\{x_i : i = 1,...,m\}$, the decision space is divided into $n$ intervals of equal probability. This ensures that the samples are distributed without clustering or overlap. The intervals are randomly sampled, and values are assigned to the respective variables, thus defining a candidate \cite{menvcik2016}. The number of intervals is defined by the desired sample size and the complexity of the problem.

\begin{figure}[H]
\centering
\begin{subfigure}[b]{0.66\textwidth}
    \centering
    \includegraphics[trim={1cm 1cm 1cm 1cm}, clip, height=5cm]{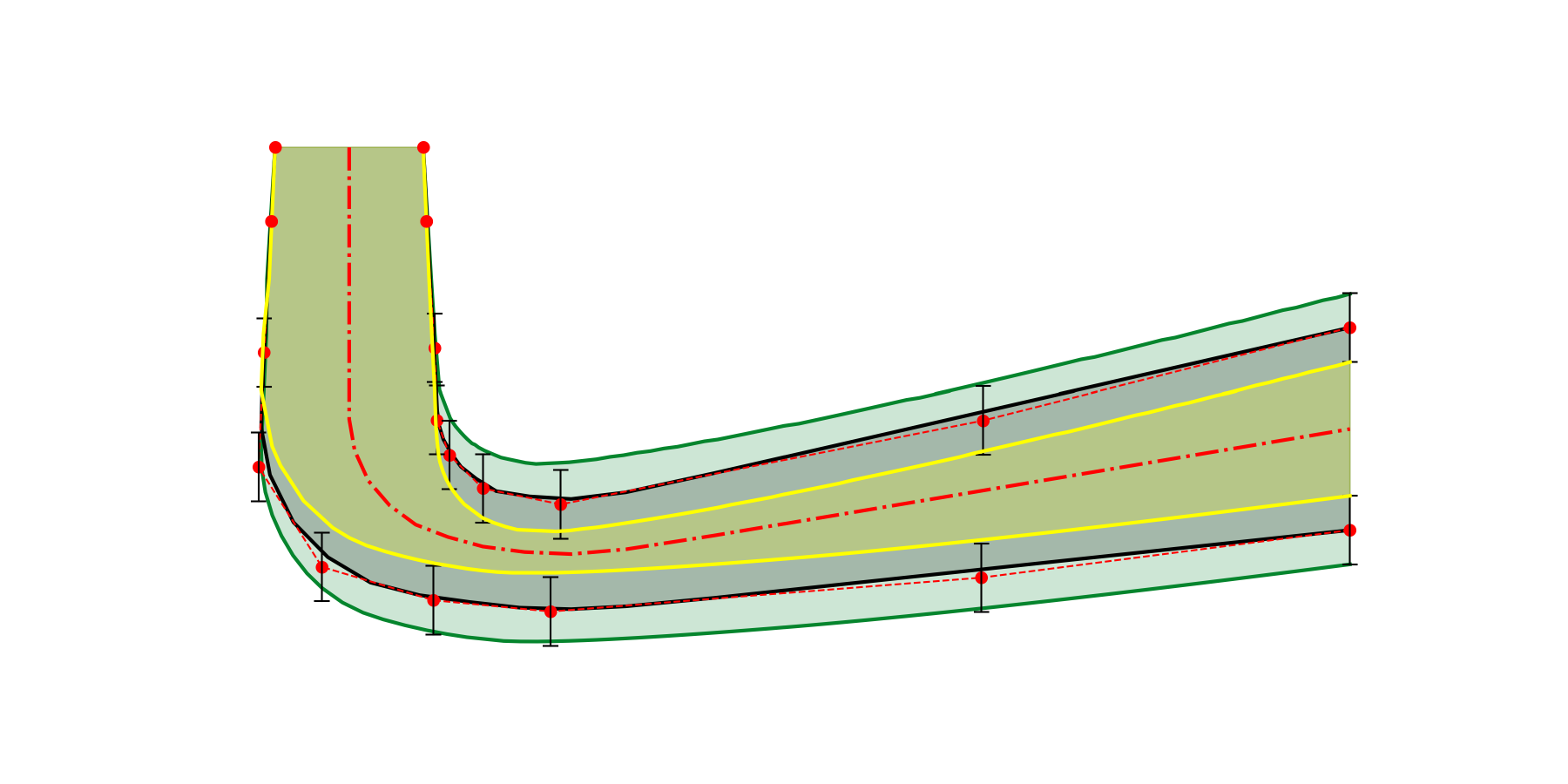}
    \caption{}
    \label{fig:model_constraints_tb}
\end{subfigure}
\begin{subfigure}[b]{0.33\textwidth}
    \centering
    \includegraphics[trim={1cm 1cm 1cm 1cm}, clip, height=5cm]{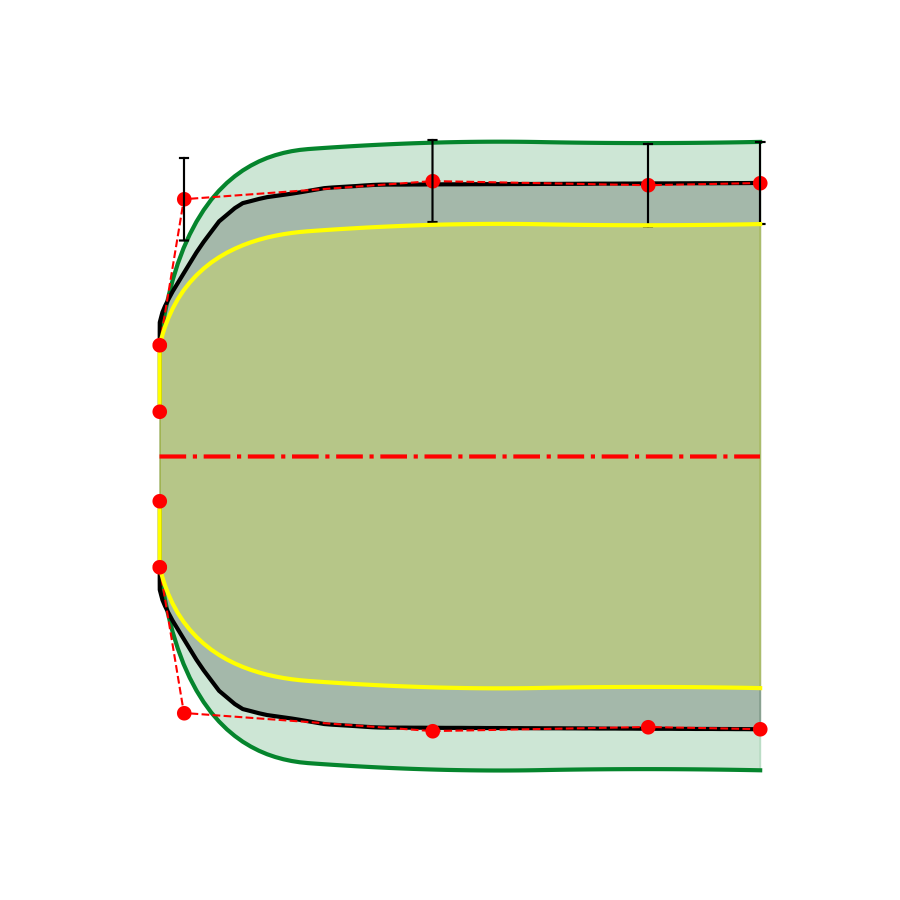}
    \caption{}
    \label{fig:model_constraints_w}
\end{subfigure}
\caption{DOE parameters. The red dots represent roof and floor B-spline control points (a) and width control points (b). Dashed lines represent control polygons. Offset limits for each control point are indicated in black. Designs shown in green and yellow are control point governed maxima and minima, respectively.}
\label{fig:model_constraints}
\end{figure}

In this study, the point offsets for the B-spline control points are DOE parameters and optimisation variables. The default offset limit for control points is $-0.25$ m $\leq x_i \leq$ $0.25$ m. The initial positions of the control points and their respective offset limits are shown in Figure \ref{fig:model_constraints}.

Two design scenarios are assessed. The control points of the width curve ($W_i$) are kept constant in the first scenario and are considered as parameters in the second scenario. In both scenarios, the control points of the roof ($R_i$) and floor ($F_i$) B-spline curves are treated as parameters\replacedRIII{. Thus}{, thus}, the number of parameters is 14 and 18, respectively. A total of 5000 samples were generated and evaluated for each scenario. The design parameters for the scenarios considered are given in Table \ref{tab:model_constraints}.

\begin{table}[H]
\centering
\footnotesize
\caption{DOE parameters and their respective limits. \protect\addedRIII{Parameters are control point offsets.} All units are in meters.}
\label{tab:model_constraints}
\begin{tabular}{@{}cccccc@{}}
\toprule
& \multicolumn{2}{c}{Scenario I} & \multicolumn{3}{c}{Scenario II} \\ \cmidrule[0.06em](r){1-1} \cmidrule[0.06em](lr){2-3} \cmidrule[0.06em](l){4-6}
$x_i$ & $\{R_i : i = 1,...,7\}$ & $\{F_i : i = 1,...,7\}$ & $\{R_i : i = 1,...,7\}$ & $\{F_i : i = 1,...,7\}$ & $\{W_i : i = 1,...,4\}$ \\
\replacedRI{$lb$}{LB} & -0.25 & -0.25 & -0.25 & -0.25 & -0.25 \\
\replacedRI{$ub$}{UB} &  0.25 &  0.25 &  0.25 &  0.25 &  0.25 \\ \bottomrule
\end{tabular}
\end{table}

\subsection{Surrogate modelling}

The present study deals with a multi-output regression problem, i.e. surrogate models are trained to predict the pressure recovery factor and the drag coefficient. Rudimentary deep neural networks (DNN) were implemented using TensorFlow v2.10.0 \cite{tensorflow}. The design and parameters of the DNNs were fine-tuned with the hyperparameter optimisation framework Optuna v3.1.0 \cite{optuna}. The optimal DNNs were networks with the highest mean coefficient of determination (R$^2$) for $C_p$ and $C_d$ after 500 trials. The hyperparameters to be optimised were the number of hidden and dropout layers, the number of neurons, the dropout rate, the learning rate, and the activation functions and initialisers. The last two hyperparameters were kept constant in all hidden layers. Table \ref{tab:hyperparameters} provides a summary of these hyperparameters.

\begin{table}[H]
\centering
\footnotesize
\caption{Assessed hyperparameters of deep neural networks and their respective value ranges.}
\label{tab:hyperparameters}
\begin{tabular}{@{}lc@{}}
\toprule
No. hidden layers, $L$              & $x,~1 \leq x \in \mathbb{N} \leq 5$ \\
No. neurons in layer $i$, $N_i$     & $2 \cdot x,~1 \leq x \in \mathbb{N} \leq 16$ \\
Dropout rate in layer $i$, $p$      & $0.1 \cdot x,~0 \leq x \in \mathbb{N}^0 \leq 8$ \\
Learning rate, $\alpha$             & $x \cdot 10^{-y},~x=\{ 1,2,4,6,8 \},~2 \leq y \in \mathbb{N} \leq 4$ \\  
Activation functions                & ELU, LeakyReLU, ReLU, SoftPlus, Swish, Tanh \\
Initialisers                        & Glorot, He, Lecun (Normal \& Uniform) \\ \bottomrule
\end{tabular}
\end{table}

The datasets were preprocessed, and outliers were removed using the unsupervised local outlier factor (LOF) algorithm, which identifies outliers in relation to their neighbours. Less than 3\% of the data in both datasets was excluded. \addedRI{In \ref{sec:appendix:a}, Figures \ref{fig:distribution_fixed} and \ref{fig:distribution_free} show the distribution of values for each feature, including targets.} The data was then split into a training set and a test set in an 80/20 ratio. All data were scaled with MinMaxScaler to a $[0,1]$ range. The test set was kept as the validation set, and the training set was further subdivided using 5-fold cross-validation. EarlyStopping was used to prevent overfitting, and the number of epochs and batch size were kept constant at 512 and 32, respectively. A stochastic gradient descent-based optimiser, Adam, was used.

The accuracy of the tuned models was assessed using the mean absolute percentage error (MAPE), relative root mean squared error (rRMSE) and coefficient of determination. MAPE is commonly employed to assess the accuracy of predictions. It is calculated as follows:
\begin{linenomath}
\begin{equation}
    \text{MAPE} = \frac{1}{n} \sum\limits_{i=1}^n  \left|\frac{y_i - \hat{y}_i}{y_i}\right| \cdot 100~\addedRIII{,}
\end{equation}
\end{linenomath}
where $n$ is the sample size, $\hat{y}_i$ is the predicted value and ${y_i}$ is the observed value. The predictive error of the regression models was also assessed using the relative root mean squared error (rRMSE). The rRMSE provides the appropriate context for the error (scale) and can be derived as:
\begin{linenomath}
\begin{equation}
    \text{rRMSE} = \frac{{\rm RMSE}}{\max( y ) - \min( y )} = \frac{\sqrt{ \frac{1}{n} \sum\limits_{i=1}^n \bigl( y_i - \hat{y}_i \bigl) ^2}}{\max( y ) - \min( y )}~.
\end{equation}
\end{linenomath}
The coefficient of determination, R$^2$, is a measure of how much of the variance of the dependent variable is explained by the model. It is calculated as follows:
\begin{linenomath}
\begin{equation}
    \text{R}^2 = 1 - \frac{ \sum\limits_{i=1}^n \bigl( y_i - \hat{y}_i \bigl) ^2 }{ \sum\limits_{i=1}^n \bigl( y_i - \bar{y} \bigl) ^2 }~\addedRIII{,}
\end{equation}
\end{linenomath}
where $\bar{y}$ is the mean of all values.

\subsection{Optimisation algorithms}

Fine-tuned surrogate models are used in conjunction with single- and multi-objective optimisation algorithms. SO optimisation allows us to understand how different design configurations affect the objectives separately. In contrast, the MO optimisation provides a more comprehensive assessment of draft tube performance since $C_p$ and $C_d$ are evaluated together. The optimisation vectors consist of 14 or 18 variables depending on the surrogate and the problem. The number of iterations in all instances was 500.

Single-objective optimisation algorithms FWA, PSO and L-SHADE are evaluated. Algorithms are tasked to maximise $C_p$:
\begin{linenomath}
\begin{equation}
    \begin{aligned}
    & \underset{\bm{x}}{\text{minimise}}
    & & f(\bm{x}) = -C_p (\bm{x}) \\
    & \text{subject to}
    & & lb \leq \bm{x} \leq ub
    \end{aligned}
\end{equation}
\end{linenomath}
and minimise $C_d$:
\begin{linenomath}
\begin{equation}
    \begin{aligned}
    & \underset{\bm{x}}{\text{minimise}}
    & & f(\bm{x}) = C_d (\bm{x}) \\
    & \text{subject to}
    & & lb \leq \bm{x} \leq ub~\addedRIII{,}
    \end{aligned}
\end{equation}
\end{linenomath}
where $\bm{x}$ is an $n$-dimensional vector of design variables in decision space $\mathbb{R}^n$\addedRI{, while $lb$ and $ub$ are lower and upper bounds, respectively}.

FWA is known for its ability to handle complex, multi-modal problems and \deletedRIII{for} its strong global search capability with fast convergence. However, it can be computationally demanding and sensitive to the parameter settings \cite{tan2015}. PSO is a simple and efficient choice for well-structured design spaces. It can struggle with complex search spaces and local optima \cite{gad2022}. On the other hand, L-SHADE is excellent for large-scale optimisations due to its adaptive population size strategy. However, it requires adequate computational resources and parameter tuning for optimal performance \cite{tanabe2014}. The implementation of these algorithms in the Python module Indago v0.2.7 \cite{indago} is used. The parameters for each algorithm are given in Table \ref{tab:SO_parameters}. The parameters were chosen in accordance with the \replacedRIII{values}{suggestions found} in the \replacedRIII{implementation papers}{literature}.

\begin{table}[H]
\centering
\footnotesize
\caption{Parameters used by the single-objective optimisation algorithms.}
\label{tab:SO_parameters}
\begin{tabular}{@{}cccccccccc@{}}
\toprule
\multicolumn{3}{c}{FWA} & \multicolumn{3}{c}{PSO} & \multicolumn{4}{c}{L-SHADE} \\ \cmidrule[0.06em](r){1-3} \cmidrule[0.06em](lr){4-6} \cmidrule[0.06em](l){7-10}
$n$ & $m_1$ & $m_2$ & $n$ & $w$   & $AKB_{model}$ & $n^{init}$ & $H$ & $r^{arc}$ & $p$ \\
20  & 10    & 10    & 20  & $AKB$ & $languid$     & 200        & 6   & 2.6       & 0.11 \\ \bottomrule
\end{tabular}
\end{table}

The multi-objective algorithms MOEA/D, SPEA2 and NSGA-II are also assessed, i.e. an optimisation problem is formulated:
\begin{linenomath}
\begin{equation}
    \begin{aligned}
    & \underset{\bm{x}}{\text{minimise}}
    & & F(\bm{x}) = \bigl( -C_p (\bm{x}),~C_d (\bm{x}) \bigl)^T \\
    & \text{subject to}
    & & lb \leq \bm{x} \leq ub~.
    \end{aligned}
\end{equation}
\end{linenomath}

MOEA/D is an efficient MO evolutionary algorithm that works by decomposing an MO problem into several SO sub-problems that are solved simultaneously \cite{zhang2007}. Its performance is highly dependent on the decomposition strategy. SPEA2 is a robust algorithm that uses a fitness assignment strategy based on the concept of dominance strength. It can handle any objective function, and its elitist selection mechanism helps convergence \cite{zitzler2001}. Nevertheless, it may be inadequate for problems with many objectives or a large population. NSGA-II is a simple and efficient algorithm, especially for problems with multiple arbitrary objectives. However, for optimal performance, it may be necessary to optimise its parameters \cite{deb2002}. It can also be comparatively space-inefficient and struggle in scenarios with many decision variables and objectives. The MO algorithms implemented in jMetalPy v1.5.5 \cite{jmetalpy} are used. Parameters are given in Table \ref{tab:MO_parameters}.

\begin{table}[H]
\centering
\footnotesize
\caption{Parameters used by the multi-objective optimisation algorithms.}
\label{tab:MO_parameters}
\begin{tabular}{@{}cccccccccccccccc@{}}
\toprule
\multicolumn{8}{c}{MOEA/D} & \multicolumn{4}{c}{SPEA2} & \multicolumn{4}{c}{NSGA-II} \\ \cmidrule[0.06em](r){1-8} \cmidrule[0.06em](lr){9-12} \cmidrule[0.06em](l){13-16}
$N$ & $p_m$ & $\eta$ & $CR$ & $F$ & $n_r$ & $T$ & $\delta$ & $N$ & $p_m$ & $p_c$ & $\eta$ & $N$ & $p_m$ & $p_c$ & $\eta$ \\
200 & 0.1   & 20     & 1.0  & 0.5 & 2     & 20  & 0.9      & 200 & 0.1   & 0.9   & 20     & 200 & 0.1   & 0.9   & 20 \\ \bottomrule
\end{tabular}
\end{table}

Multi-objective results are provided as Pareto fronts. The Pareto front is a collection of non-dominated solutions such that improvement in pressure recovery factor can only be achieved at the expense of drag coefficient and vice versa. \addedRI{TOPSIS is used to identify the best MO candidate, i.e. a configuration that provides the most balanced combination of $C_p$ and $C_d$ contributions. TOPSIS is a multi-criteria decision-making technique that selects the best solution from a set of alternatives based on their geometric distance from a positive ideal solution and their distance from a negative ideal solution \cite{hwang1981}.}

\section{Results and discussion}

\subsection{Grid convergence}

A grid convergence study was conducted to evaluate the spatial convergence of steady-state simulations \cite{roache1998}. Three grids were created for the reference draft tube design, namely coarse, medium, and fine with $2.1 \cdot 10^{5}$, $4.2 \cdot 10^{5}$ and $7.4 \cdot 10^{5}$ cells, respectively. The scaling factor, $f_s$, was 0.8 for the fine grid, 1.2 for the medium grid, and 1.8 for the coarse grid. The first grid point was in the log-law region in all cases. The grid convergence index (GCI) for the pressure recovery factor and the drag coefficient is reported in Table \ref{tab:gci}. The results indicate that the solutions are in the asymptotic range of convergence, i.e. $GCI_{m,f} / (r^p \cdot GCI_{c,m}) \approx 1$ for both $C_p$ and $C_d$. Therefore, the medium grid where $f_s = 1.2$ was chosen for CFD simulations.

\begin{table}[H]
\centering
\footnotesize
\caption{Results of the grid convergence study for $C_p$ and $C_d$.}
\label{tab:gci}
\begin{tabular}{@{}ccccccc@{}}
\toprule
$f(x)$ & $\varepsilon_{c,m}$ (\%) & $\varepsilon_{m,f}$ (\%) & $GCI_{c,m}$ (\%) & $GCI_{m,f}$ (\%) & $p_{GCI}$ & $\frac{GCI_{m,f}}{r^p \cdot GCI_{c,m}}$ \\ \midrule
$C_p$ &  1.575 & 0.563 & 1.084 & 0.387 & 2.553 & 0.994 \\
$C_d$ & 10.055 & 4.252 & 9.944 & 4.206 & 2.015 & 1.044 \\ \bottomrule
\end{tabular}
\end{table}

\subsection{Surrogate model tuning and validation}

Two deep neural network models were trained. The first model was trained on a dataset where the width of the draft tube was kept constant while the control point offsets for the roof and floor B-spline changed. The second model was trained on a dataset where the width change was considered. A unified model trained on data from both DOE experiments was not pursued because a combined dataset would be biased and skewed due to objective differences in the sampling approach. The hyperparameters of the surrogate models are summarised in Table \ref{tab:DNN_parameters}. The learning rate $\alpha = 0.002$ and the kernel initialiser LecunNormal are used by both models. Interestingly, neither network uses dropout regularisation. This can be attributed to several factors. The use of EarlyStopping prevents overfitting and, therefore, competes with dropout in this aspect. Since EarlyStopping is always applied, a combination of both approaches has a negative impact on performance. Additionally, the adaptive nature of the design, where the number of neurons per layer and the total number of layers vary, reduces the need for dropout regularisation. Figure \ref{fig:DNN} shows the tuned network architectures.

\begin{table}[H]
\centering
\footnotesize
\caption{Tuned hyperparameters for each network obtained using Optuna framework.}
\label{tab:DNN_parameters}
\begin{tabular}{@{}lcc@{}}
\toprule
 & Scenario I & Scenario II \\ \midrule
No. hidden layers, $L$              & 5 & 5 \\
No. neurons in layer $i$, $N_i$     & 22, 22, 20, 24, 4 & 26, 24, 32, 12, 6 \\
Dropout rate in layer $i$, $p$      & - & - \\
Learning rate, $\alpha$             & 0.002 & 0.002 \\  
Activation functions                & Swish & ELU \\
Initializers                        & LecunNormal & LecunNormal \\
Batch size                          & 32 & 32 \\
Epochs                              & 512 & 512 \\ \bottomrule
\end{tabular}
\end{table}

The performance of the tuned DNNs was assessed using 5-fold cross-validation. Subsequently, the models were trained on all available data, excluding the validation set, and evaluated. The performance metrics for the 5-fold cross-validation and the final models are summarised in Table \ref{tab:DNN_metrics}. The MAPE and rRMSE values are consistently below 5\% (below 4.5\% and 3\% for the final models), indicating a reasonably low error rate in predicting the target values. Moreover, the R$^2$ score is higher than 0.9, which implies that the models capture considerable variance in the target variables. Overall performance is improved when the models are trained on all available data.

\begin{table}[H]
\centering
\footnotesize
\caption{Metrics and standard deviations for DNNs trained on different datasets. The values are the worst outcomes for $C_p$ and $C_d$. Values in bold are metrics when the model is trained on all available data (excluding validation).}
\label{tab:DNN_metrics}
\begin{tabular}{@{}ccccc@{}}
\toprule
 & Scenario I (5-fold) & Scenario I & Scenario II (5-fold) & Scenario II \\ \midrule
MAPE & 2.829 $\pm$ 0.133 & \textbf{2.548} & 3.169 $\pm$ 0.217 & \textbf{2.827} \\
rRMSE & 4.166 $\pm$ 0.165 & \textbf{3.857} & 4.703 $\pm$ 0.313 & \textbf{4.252} \\
R$^2$ & 0.935 $\pm$ 0.006 & \textbf{0.946} & 0.913 $\pm$ 0.009 & \textbf{0.929} \\ \bottomrule
\end{tabular}
\end{table}

\subsection{Optimisation results}

The results of optimisation using DNN surrogates and different optimisation algorithms are presented below. Two DNN models, Scenario I and Scenario II, are used in cases where the objective is to generate a generic draft tube design (a) with no inherent constraints (apart from set lower and upper limits). In addition, cases are defined where the objective is to generate an optimal design within the limits of the reference design (b). This implies that any changes to the reference design are merely additive, while there are no such requirements in the first case. These requirements are set as limits for the optimisation variables in the optimisation process and are summarised in Table \ref{tab:scenarios_tests}.

\begin{table}[H]
\centering
\footnotesize
\caption{The upper and lower limits for the optimisation variables in all test cases. All units are in meters.}
\label{tab:scenarios_tests}
\begin{tabular}{@{}ccccc@{}}
\toprule
& Scenario I.a & Scenario II.a & Scenario I.b & Scenario II.b \\ \cmidrule[0.06em](r){1-1} \cmidrule[0.06em](lr){2-2} \cmidrule[0.06em](lr){3-3} \cmidrule[0.06em](lr){4-4} \cmidrule[0.06em](l){5-5}
$x_i$ & $\{X_i : i = 1,...,14\}$ & $\{X_i : i = 1,...,18\}$ & $\{X_i : i = 1,...,14\}$ & $\{X_i : i = 1,...,18\}$   \\
\replacedRI{$lb$}{LB} & -0.25,~for $1 \leq i \leq 14$ & -0.25,~for $1 \leq i \leq 18$ & \footnotesize\makecell[c]{-0.25,~for $1 \leq i \leq 7$ \\[-2mm] 0,~for $8 \leq i \leq 14$} & \footnotesize\makecell[c]{-0.25,~for $1 \leq i \leq 7$ \\[-2mm] 0,~for $8 \leq i \leq 14$ \\[-2mm] -0.25,~for $15 \leq i \leq 18$} \\
\replacedRI{$ub$}{UB} & 0.25,~for $1 \leq i \leq 14$ & 0.25,~for $1 \leq i \leq 18$ & \footnotesize\makecell[c]{0,~for $1 \leq i \leq 7$ \\[-2mm] 0.25,~for $8 \leq i \leq 14$} & \footnotesize\makecell[c]{0,~for $1 \leq i \leq 7$ \\[-2mm] 0.25,~for $8 \leq i \leq 14$ \\[-2mm] 0,~for $15 \leq i \leq 18$} \\ \bottomrule
\end{tabular}
\end{table}

\subsubsection{Single-objective optimisation results}

An assessment of convergence for FWA, PSO and L-SHADE was performed. The pressure recovery factor and drag coefficient were evaluated separately. In Scenario I.a, for both $C_p$ and $C_d$, L-SHADE converged first after about 100 iterations. Although the difference was minimal, FWA outperformed PSO when optimising for $C_p$. However, for $C_d$, FWA remained at a local minimum. The results for the pressure recovery factor were comparable, while the difference for the drag coefficient was $\approx 5\%$ (Figure \ref{fig:cvg_scenarioIa}). In terms of computational efficiency, L-SHADE requires significant resources as it evaluates a large population. PSO is the most cost-effective as it can achieve satisfactory convergence with a much smaller evaluation budget. However, considering the simplicity and short computation times of the objective functions, L-SHADE is the optimal choice. For Scenario I.b, a similar pattern can be observed (Figure \ref{fig:cvg_scenarioIb}). In this case, the FWA still lags behind but can achieve satisfactory convergence with results comparable to those of the competition. The results indicate that the limits of the optimisation variables have a significant impact on FWA's convergence. Adjustments to algorithm parameters and setting stricter limits should improve performance.

L-SHADE is the fastest in terms of overall convergence in Scenario II.a. Differences between results are negligible (Figure \ref{fig:cvg_scenarioIIa}). This is also true for Scenario II.b. All algorithms are adequate when optimising for $C_d$. However, L-SHADE is most effective for $C_p$ and can avoid local minima (Figure \ref{fig:cvg_scenarioIIb}).

\subsubsection{Multi-objective optimisation results}

Three algorithms, MOEA/D, SPEA2 and NSGA-II, are considered in this assessment. The fronts for Scenario I.a show a remarkable similarity and overlap for different algorithms (Figure \ref{fig:pf_scenarioIa}). The same conclusion can be drawn for Scenario II.a (Figure \ref{fig:pf_scenarioIIa}). In both cases, the Pareto-optimal solutions outperform the reference draft tube design. It is worth noting that the Pareto front shifts downwards and to the right (a decrease in $C_d$ and an increase in $C_p$) when the width of the draft tube is allowed to change. The sloping shape of the Pareto front is a direct consequence of the prescribed limits for the optimisation variables.

\begin{figure}[H]
    \centering
    \includegraphics[width=0.9\textwidth]{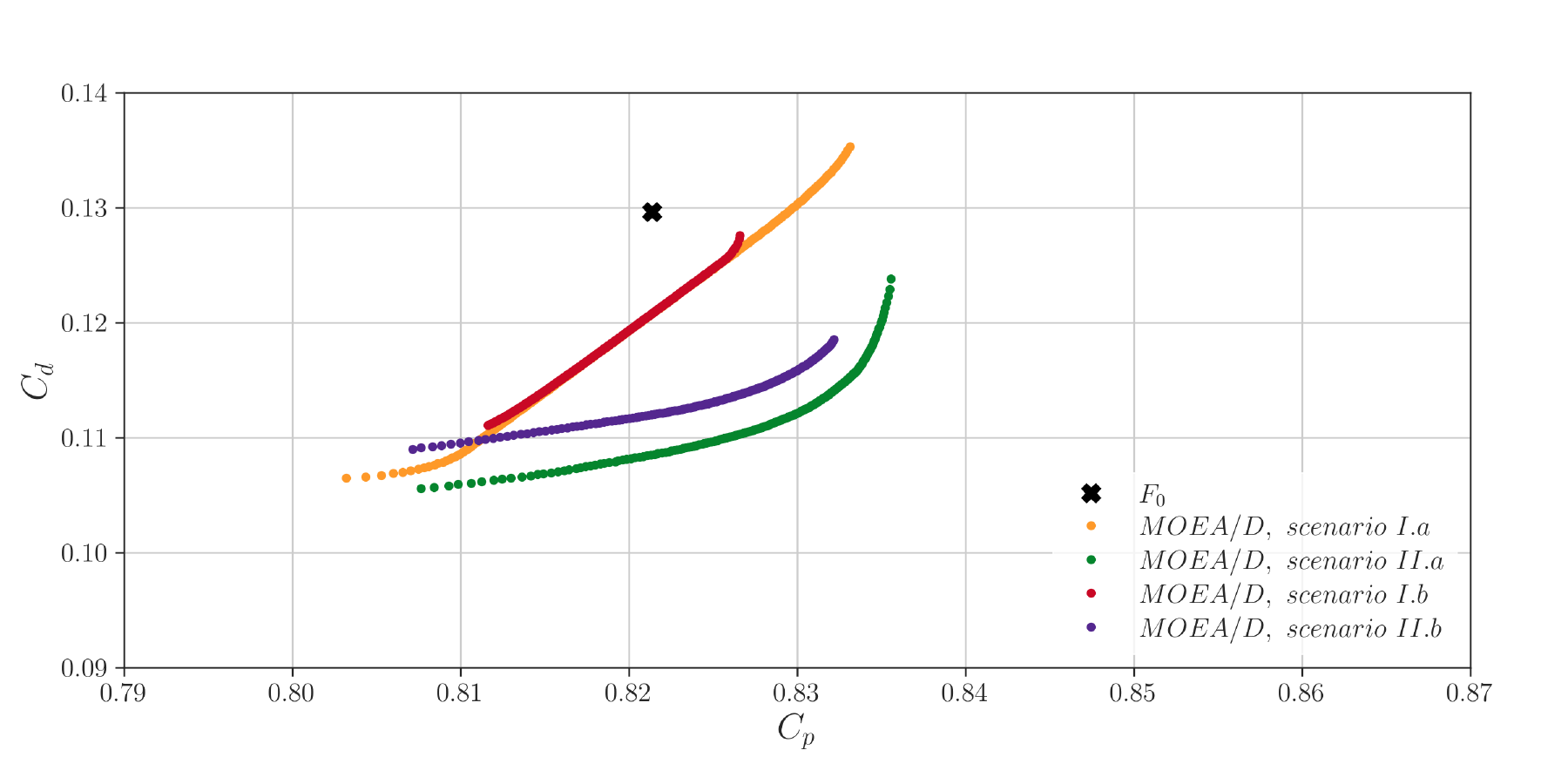}\vspace{-0.5cm}
    \caption{\protect\addedRI{The extent of the Pareto fronts for different test cases when using MOEA/D. The central front segments for Scenario I.a and Scenario I.b coincide. In both scenarios, the width of the draft tube is constant. However, the range of potential non-dominated solutions is smaller in Scenario I.b due to the imposed bound limits on the design of the draft tube. For Scenario II.a and Scenario II.b, the width of the draft tube is considered in the optimisation process. As the shape of the design is limited to the original draft tube extents in Scenario II.b, the obtained non-dominated solutions perform worse than in Scenario II.a; flexible design limits allow for more unconventional designs and improved Pareto-optimal solutions. The black marker indicates the results for the original draft tube design. As shown, it is possible to find an optimal solution on the Pareto front where both $C_p$ and $C_d$ are improved compared to the original design.}}
    \label{fig:pf_MOEAD}
\end{figure}

The results for Scenario I.b (Figure \ref{fig:pf_scenarioIb}) and Scenario II.b (Figure \ref{fig:pf_scenarioIIb}) demonstrate a similar trend. The shapes of the fronts are similar to the cases where there are no design constraints. However, due to the limits imposed on the optimisation variables, the fronts are narrower and have shifted upward and to the left (a decrease in $C_p$ and an increase in $C_d$), indicating worse results. Nevertheless, all Pareto-optimal solutions in these scenarios outperform the reference draft tube design. Figure \ref{fig:pf_MOEAD} compares the Pareto fronts for four cases using MOEA/D. The results reveal a dynamic relationship between variable limits and the resulting Pareto front in the context of draft tube design optimisation. MOEA/D has henceforth been adopted as the default MO algorithm.

\subsubsection{Parameter influence}

The influence of the key parameters of the algorithms L-SHADE and MOEA/D on the convergence and quality of the solutions for draft tube optimisation was investigated. Scenario II.a was considered. Table \ref{tab:opt_algo_parameters} summarises the tests performed for each algorithm. For L-SHADE, the initial population size $n^{init}$, the history size $H$ and the archive ratio $r^{arc}$ were evaluated. For MOEA/D, the population size $N$, the crossover rate $CR$, the mutation factor $F$, the number of replaced solutions $n_r$, neighbourhood size $T$ and the probability of updating the neighbour $\delta$ were evaluated.

\begin{table}[H]
\centering
\footnotesize
\caption{L-SHADE and MOEA/D parameter evaluation matrix.}
\label{tab:opt_algo_parameters}
\begin{tabular}{@{}ccccccccc@{}}
\toprule
\multicolumn{3}{c}{L-SHADE} & \multicolumn{6}{c}{MOEA/D} \\ \cmidrule[0.06em](r){1-3} \cmidrule[0.06em](l){4-9}
$n^{init}$ & $H$ & $r^{arc}$ & $N$ & $CR$ & $F$ & $n_r$ & $T$ & $\delta$ \\
50  & 6  & 1.0 & 50  & 0.7 & 0.5 & 1 & 10 & 0.6 \\
100 & 10 & 1.5 & 100 & 0.8 & 0.6 & 2 & 20 & 0.7 \\
200 & 15 & 2.0 & 200 & 0.9 & 0.7 & 4 & 40 & 0.8 \\
400 & 20 & 2.6 & 400 & 1.0 & 0.8 & 8 & 80 & 0.9 \\ \bottomrule
\end{tabular}
\end{table}

As expected, increasing $n^{init}$ leads to faster convergence for L-SHADE. The default value $n^{init}=200$ achieves an acceptable balance between convergence and computational efficiency. Values smaller than $n^{init}<100$ are not recommended. The history size $H$ has no or negligible influence on convergence. The convergence improves with the archive ratio; the suggested value $r^{arc}=2.6$ is suitable, although higher values could improve convergence (Figure \ref{fig:comparisson_LSHADE}).

The change in population size $N$ for MOEA/D has no discernible impact on the results. However, the fronts are captured more accurately as $N$ and the total number of evaluations increases. For $CR=0.7$ and $CR=0.9$, the fronts are slightly extended downward. The behaviour near the anchor for $CR=0.7$ indicates that the standard value may be more suitable. The mutation factor and the neighbourhood size have no visible influence on the front. The downward shift of the Pareto fronts for $n_r$ and $\delta$ can be observed (Figure \ref{fig:comparisson_MOEAD}). This behaviour can be attributed to algorithmic data diversity and exploration tendencies; changing the settings that govern diversity can lead to broader exploration. However, the differences are minor, and the default parameters should be sufficient.

\subsection{Validation of the workflow}

The validity of the proposed workflow for draft tube design optimisation is henceforth assessed by comparing the results of CFD simulations with the predictions provided by MO-DNNs for optimal design solutions. The results in Table \ref{tab:comparison_topsis} show that the discrepancy between the DNN-predicted results and those calculated using CFD is consistently less than 3\% for the drag coefficient and 0.5\% for the pressure recovery factor. This indicates that the proposed workflow can provide designs with results that closely match the results of the CFD simulations. The best design variant, Scenario II.a, shows a reduction of $\approx 17\%$ for the drag coefficient and an increase of $\approx 1.5\%$ for the pressure recovery factor compared to the reference design. Even with the errors taken into account, these results are still a significant improvement. CFD results for the best solutions proposed by TOPSIS are shown in Figure \ref{fig:cfd_comparisson}.

\begin{table}[H]
\centering
\footnotesize
\caption{Comparison between $C_p$ and $C_d$ values for TOPSIS-proposed best solutions calculated using CFD and predicted by respective DNN models. The difference between the CFD results and the DNN predictions is less than 3\% for $C_d$ and 0.5\% for $C_p$.}
\label{tab:comparison_topsis}
\begin{tabular}{@{}lcccc@{}}
\toprule
& Scenario I.a & Scenario II.a & Scenario I.b & Scenario II.b \\ \midrule
$C_{p, CFD}$ & 0.813 & 0.831 & 0.824 & 0.826 \\
$C_{p, DNN}$ & 0.815 & 0.829 & 0.820 & 0.825 \\
$\delta_{C_p}~[\%]$ & 0.172 & -0.334 & -0.472 & -0.222 \\ \midrule
$C_{d, CFD}$ & 0.117 & 0.109 & 0.116 & 0.111 \\
$C_{d, DNN}$ & 0.114 & 0.111 & 0.120 & 0.113 \\
$\delta_{C_d}~[\%]$ & -2.693 & 2.055 & 2.793 & 1.470 \\ \bottomrule
\end{tabular}
\end{table}

\subsection{Comparative assessment}

The results obtained using MOEA/D and L-SHADE algorithms for all scenarios are \addedRIII{henceforth} compared to determine the influence of the objectives in single-objective optimisation and the combined effects in multi-objective optimisation. Presented MO designs are TOPSIS optimal.

\begin{figure}[H]
    \centering
    \includegraphics[width=0.9\textwidth]{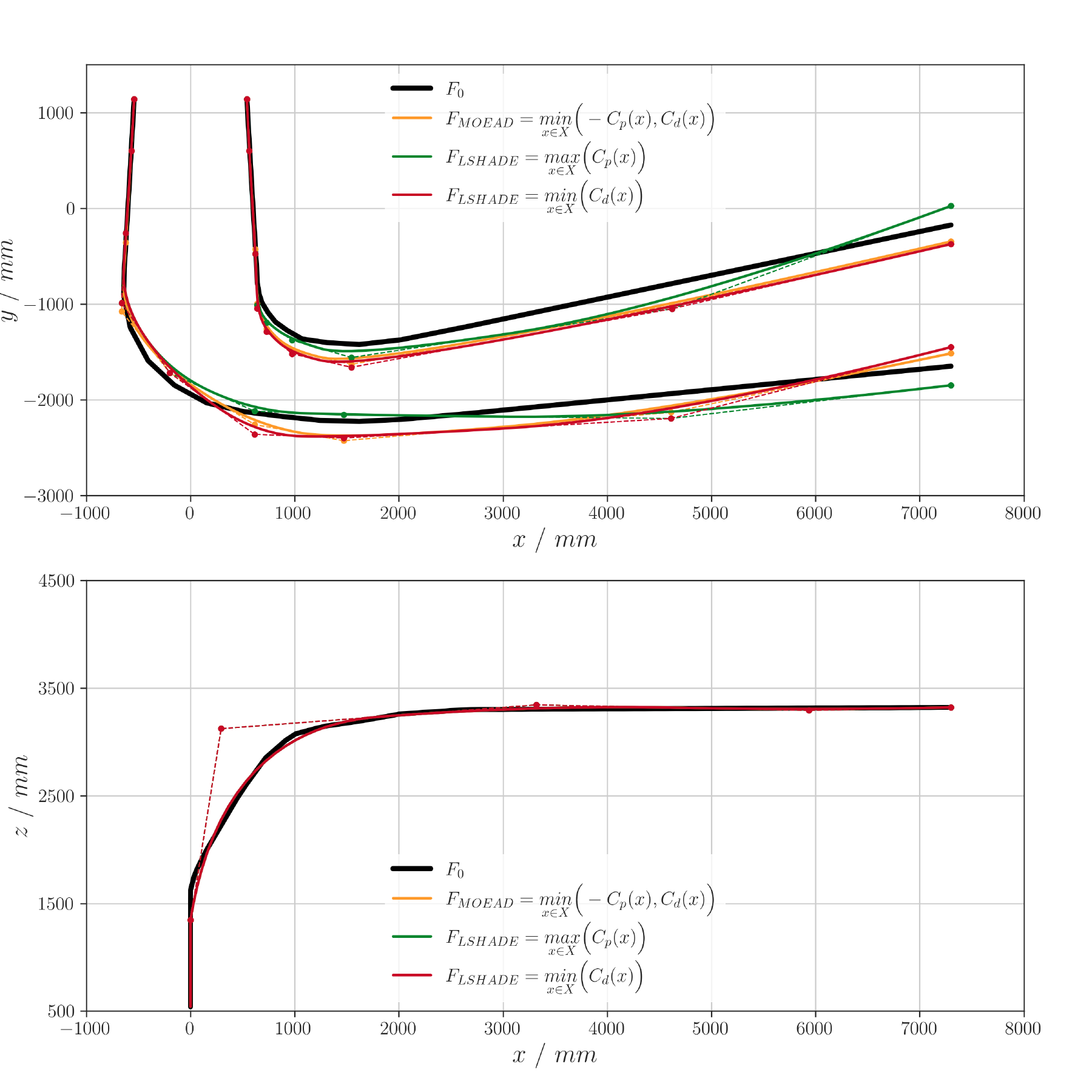}\vspace{-0.5cm}
    \caption{\protect\addedRI{Optimal draft tube designs for Scenario I.a. The reference design is portrayed with black lines. The dots are control points. Dashed lines represent control polygons. The limits for roof and floor control points are set to $-0.25 \leq X_i \leq 0.25$, hence the difference in the resulting curves for different objectives. The width control points are fixed, i.e. the width curve is constant. The MO result is comparable in shape to the SO result for $C_d$, indicating that $C_d$ has greater importance in this scenario. The increase in the cross-sectional area at the outlet for $C_p$ is a direct consequence of the pressure recovery improvement objective (higher static pressure at the outlet).}}
    \label{fig:designs_scenarioIa}
\end{figure}

The Pareto front and SO optimums for $C_p$ and $C_d$ in Scenario I.a are congruent (Figure \ref{fig:compare_scenarioIa}). It can be argued that SO optimisation can be \addedRIII{considered as} an alternative to MO, depending on the \addedRIII{design} requirements. In this design test, only the B-spline control points for the roof and floor are allowed to move, as shown in Figure \ref{fig:designs_scenarioIa}. The elbow of the draft tube is offset downward, and the diffuser points upward. The MO results and the results when $C_d$ is set as an objective are similar, suggesting that $C_d$ has greater importance than $C_p$ in this \addedRIII{MO} optimisation process. \addedRIII{The sudden change in outlet cross-sectional area when $C_p$ is the objective is indicative of the physics of the problem; as the cross-sectional area increases, so does the static pressure, resulting in a greater pressure recovery factor. This increase is at the expense of the drag coefficient.}

\begin{figure}[H]
    \centering
    \includegraphics[width=0.9\textwidth]{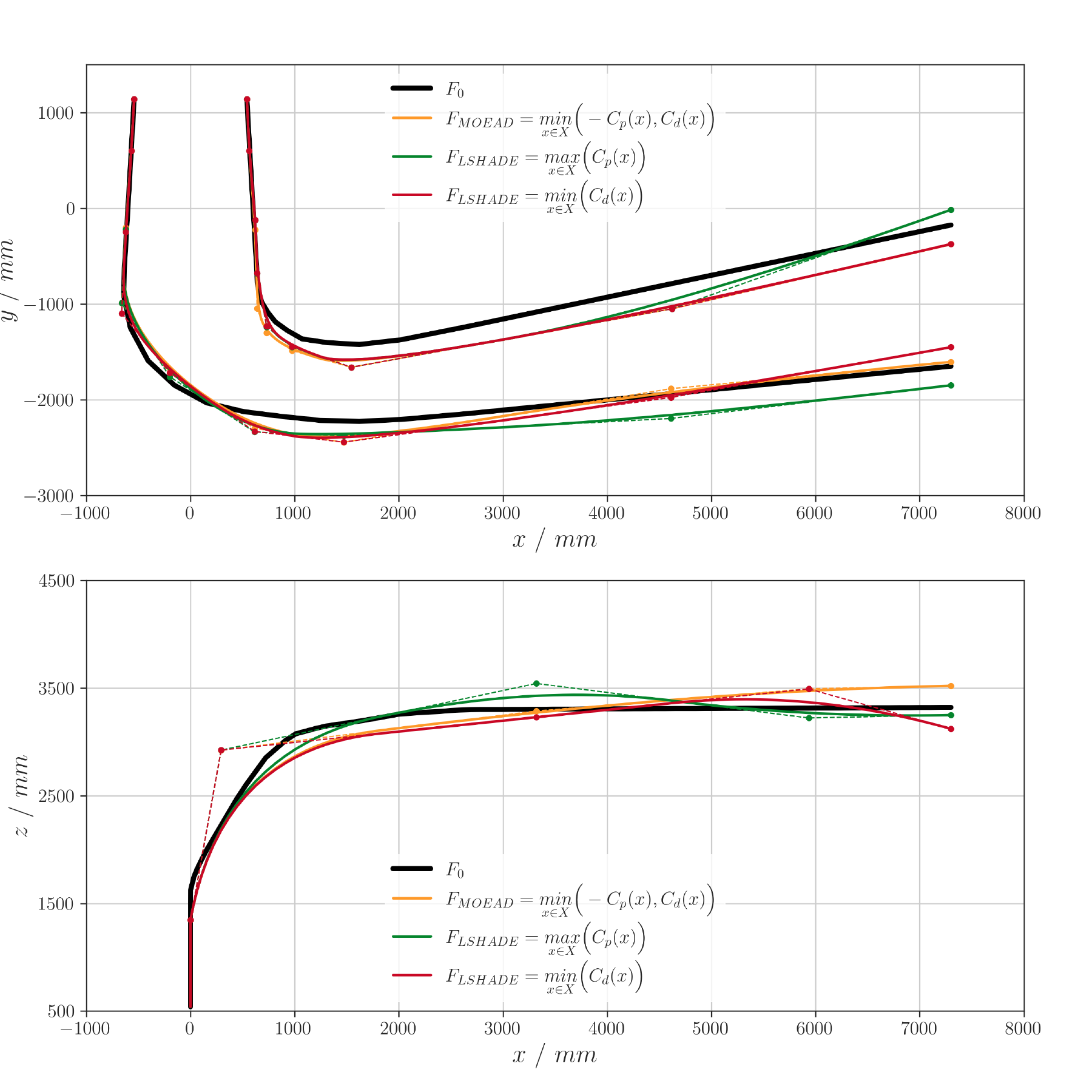}\vspace{-0.5cm}
    \caption{\protect\addedRI{Optimal draft tube designs for Scenario II.a. The reference design is portrayed with black lines. The dots are control points. Dashed lines represent control polygons. The limits for roof, floor and width control points are set to $-0.25 \leq X_i \leq 0.25$. Since there are no imposed constraints, the respective curves vary considerably depending on the set objectives. The cross-sections near the elbow remain consistent, while the cross-sections near the outlet are adjusted to compensate for the pressure drop. When the objective is to maximise the pressure recovery, $C_p$, the cross-sectional area is the largest. The MO design achieves a balance between improved pressure recovery and reduced drag.}}
    \label{fig:designs_scenarioIIa}
\end{figure}

For Scenario II.a, the Pareto front and the SO optima do not coincide (Figure \ref{fig:compare_scenarioIIa}). However, if we extrapolate the front, the SO results would be its anchors. The results for $C_p$ and $C_d$ have improved. In terms of design, the elbows now overlap. The outlet position has \replacedRIII{changed slightly}{been slightly changed}. The diffuser \addedRIII{section} (except for $C_p$) now has a strong upward trend. The width curve for the MO optimisation shows a steady increase towards the outlet, in contrast to the results for \addedRIII{the} SO optimisation. \addedRIII{Compared to Scenario I.a, the pressure recovery factor is increased and the drag coefficient is reduced for TOPSIS optimal design. This is primarily due to the changes in the cross-sectional area governed by the width curve.} Figure \ref{fig:designs_scenarioIIa} shows the optimal shapes generated.

\begin{figure}[H]
    \centering
    \includegraphics[width=0.9\textwidth]{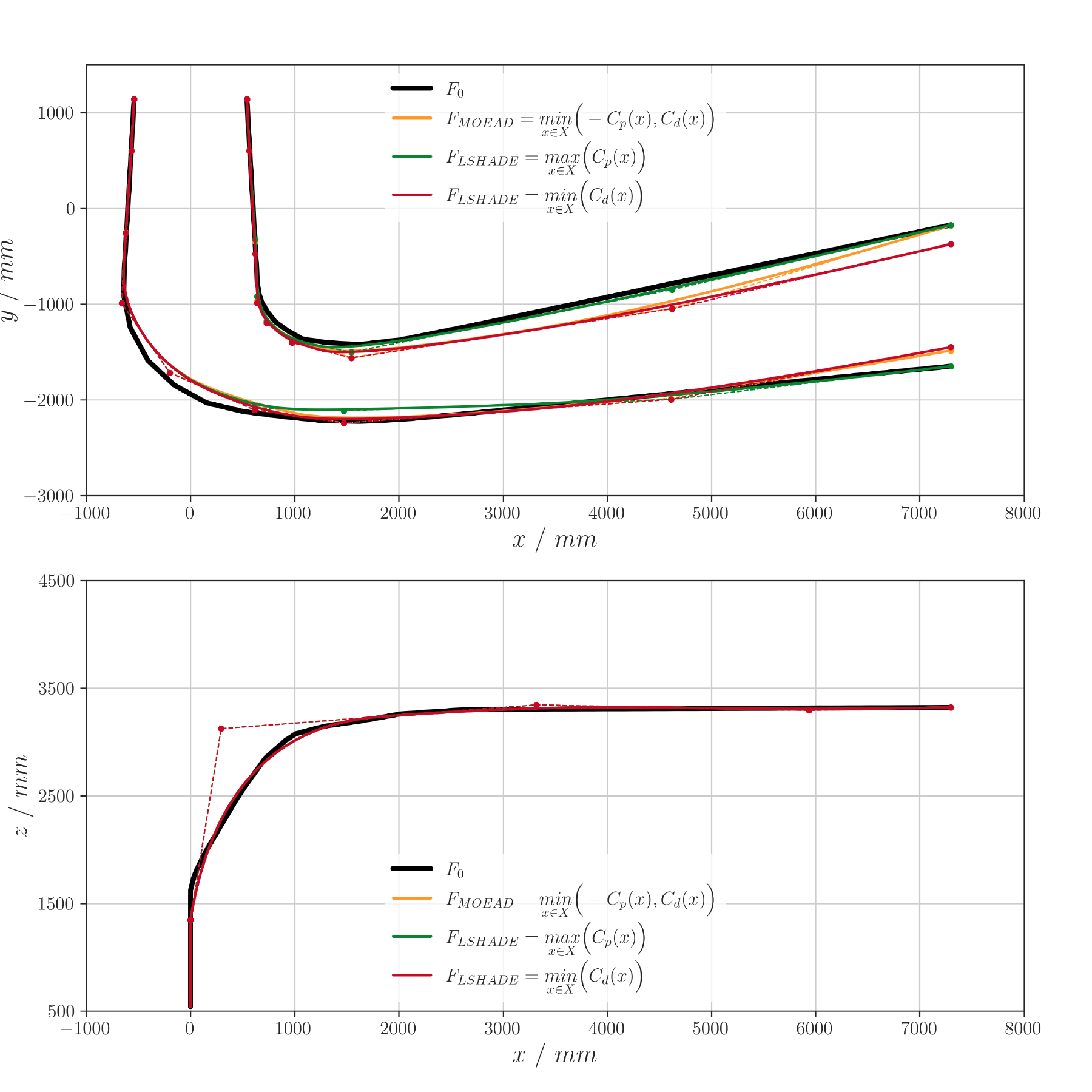}\vspace{-0.5cm}
    \caption{\protect\addedRI{Optimal draft tube designs for Scenario I.b. The reference design is portrayed with black lines. The dots are control points. Dashed lines represent control polygons. All designs are restricted to the interior of the reference design. The width control points are fixed, i.e. the width curve is constant. The design for $C_p$ largely follows the reference design. The design for $C_d$ reduces drag at the expense of pressure recovery, resulting in a smaller cross-sectional area near the outlet. The elbow offset can be linked to the fixed inflow rate.}}
    \label{fig:designs_scenarioIb}
\end{figure}

The \replacedRIII{Pareto front}{results} for Scenario I.b \replacedRIII{is}{are} similar to \deletedRIII{those for} Scenario I.a, but \replacedRIII{narrower}{worse}. This is due to the imposed limits on the optimisation variables\addedRIII{, which restricts the range of valid designs}. Still, the front and SO optima coincide (Figure \ref{fig:compare_scenarioIb}). According to Figure \ref{fig:designs_scenarioIb}, the elbows are slightly \replacedRIII{offset}{moved} to the right. The diffuser section is narrower, except when optimising for $C_p$, where it mostly follows the reference design, suggesting that the reference was designed to improve pressure recovery. \addedRIII{Design changes directly affect the flow and static pressure i.e. maximise the conversion of kinetic energy of the exiting water into potential energy. The change in the drag coefficient is $\approx 5\%$ compared to Scenario I.a.} Due to the imposed model and flow simplifications, these results are reasonable; narrowing is not expected with realistic inlet velocity profiles.

\begin{figure}[H]
    \centering
    \includegraphics[width=0.9\textwidth]{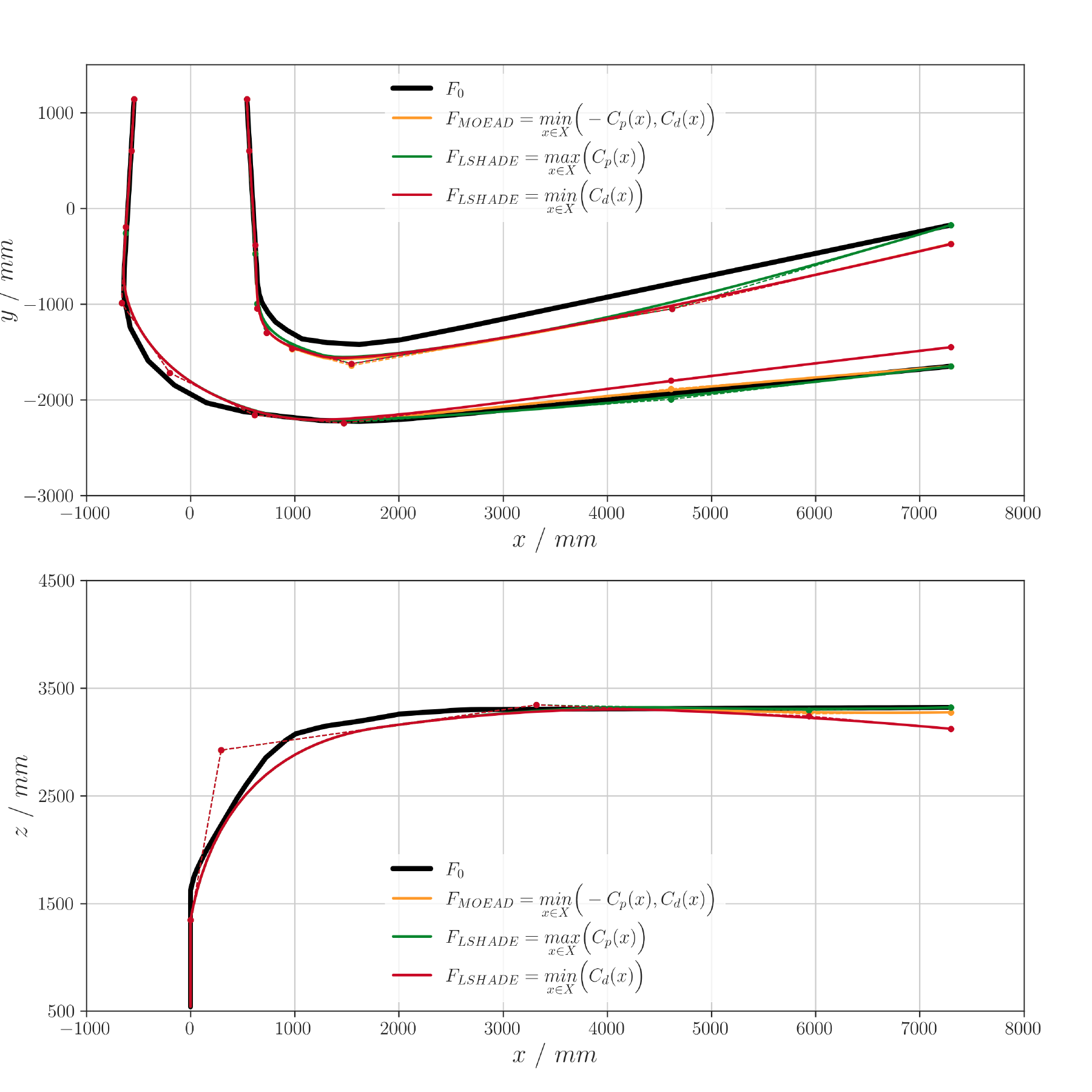}\vspace{-0.5cm}
    \caption{\protect\addedRI{Optimal draft tube designs for Scenario II.b. The reference design is portrayed with black lines. The dots are control points. Dashed lines represent control polygons. All designs are constrained to the interior of the reference design and the width curve may change. The results for $C_p$ and the MO results are largely similar. The shapes of the floor and width curves suggest that the design limits are too restrictive for set objectives.}}
    \label{fig:designs_scenarioIIb}
\end{figure}

The Pareto front for Scenario II.b \replacedRIII{outlines worse results when}{follows a similar pattern to Scenario I.b - worse results are obtained due to imposed limits} compared to a non-limited \addedRIII{(Scenario II.a)} test case (Figure \ref{fig:compare_scenarioIIb}). Compared to I.b \addedRIII{designs}, the MO results and the result when $C_p$ is used as an objective are similar, as demonstrated in Figure \ref{fig:designs_scenarioIIb}. \addedRIII{The results for $C_p$ and $C_d$ are improved, mainly due to the change in the diffuser section and overall width. Due to the limits set for the optimisation variables, the designs perform worse compared to Scenario II.a. Still, they are an improvement over the reference design, hence it is safe to conclude that even under strict limits it is possible to improve overall draft tube performance.}

\subsection{Discussion}

The proposed systematic workflow simplifies the draft tube design optimisation process by combining machine learning and CFD data. \addedRI{Since a large number of design candidates must first be evaluated to generate a comprehensive dataset, the computational effort is required.} \addedRI{However,} this data-driven approach offers more flexibility as it can also be used with experimental or model data. \deletedRI{However, a large number of design candidates must first be evaluated to generate a comprehensive dataset. This means that a computational effort is required.}

The use of DNN surrogates enables rapid evaluation of novel draft tube designs. This significantly accelerates the design exploration process and makes it possible to evaluate many design candidates. The workflow is sufficiently flexible to handle scenarios where only limited geometry adjustments are possible, which can be important for turbine revitalisation as it is usually physically or financially constrained. However, it is important to consider the computational effort required to fine-tune the DNN surrogates and the need for adequate preprocessing of the data.

We have demonstrated that L-SHADE and MOEA/D are suitable algorithms for draft tube design optimisation when standard parameters are used. Nevertheless, it may be necessary to modify the algorithmic parameters to obtain optimal results. Depending on the objective, single-objective algorithms are a viable and efficient alternative to multi-objective algorithms. Since both the SO and MO algorithms use the same surrogate, switching between algorithms is straightforward. In general, however, multi-objective optimisation offers a more comprehensive set of solutions and can be used to identify designs that meet specific criteria. The inclusion of constraints in this context can lead to suboptimal solutions. A natural extension of the proposed multi-objective approach is multi-point optimisation, where performance at different operating points is considered. The growing demand for hydropower turbine operations in a broader range of conditions makes this an important topic. Until now, due to CFD-associated computational requirements, these optimisations have not been carried out. By combining, for example, $C_p$ and $C_d$ or measuring performance in terms of overall efficiency, the MO optimisation can be adapted to provide $N$-dimensional optimal solutions, where $N$ is the number of operating points. The problem can also be formulated in a way that takes into account other factors, such as construction costs, making it a useful tool for designing new draft tubes or revitalising existing ones.

It is important to reflect on the simplifications that were made. A constant flow rate was set at the inlet instead of a velocity profile, which directly affects the flow in the draft tube and thus \replacedRIII{influences}{determines} the optimal shape. If velocity profiles are used, a more realistic shape could be achieved. However, this is not a limitation of the methodology and does not affect the proposed workflow. The use of velocity profiles could also be advantageous for cavitation tracking and prediction. This idea requires further research and may require the use of convolutional neural networks. Additionally, optimal designs are reported only for a single operating point.

\addedRI{The proposed optimised designs have a lower drag coefficient and a higher pressure recovery factor due to changes in the pressure distribution, which implies improved energy extraction from the fluid. This has a direct impact on the total turbine power output. For operational hydropower turbines, where changes can be made to the downstream sections of the turbine as part of planned maintenance, this can be particularly important as efficiency can be improved with minimal investment in a short period of time. Considering the cost and simplicity of implementation, internally constrained designs were explored as this would eliminate the need for complex excavations and other non-localised modifications. However, since the methodology is problem-independent, it can also be adapted for upstream sections and even the runner, and thus provide tangible improvements (albeit at a higher cost).}

\section{Conclusion}

This comprehensive study introduces a systematic workflow for draft tube optimisation using neural network surrogates trained on computational fluid dynamics data. \deletedRI{Latin hypercube sampling was used to generate design candidates, and the validated surrogate was}\addedRI{Validated surrogates were} coupled with single-objective and multi-objective algorithms to propose optimal draft tube designs based on flow conditions and specified constraints. The pressure recovery factor and drag coefficient were used as objectives in this assessment. The results of this study demonstrate the potential of this methodology to improve the performance of hydropower turbines. An in-depth evaluation of the optimisation algorithms was also conducted. The main findings are:

\begin{itemize}
    \item The use of deep neural network surrogates reduces \addedRIII{the} computational requirements \deletedRI{when optimising the design of draft tubes, thus facilitating}\addedRI{and facilitates} the exploration of \addedRIII{different} design candidates. \deleted{This makes them a valuable tool for both the revitalisation of ageing turbines and the design of efficient new turbines.} \addedRI{The computational requirements associated with data acquisition remain unchanged.}
    
    \item The difference between the predictions for the optimal candidates and the respective numerical results was less than 0.5\% for the pressure recovery factor and 3\% for the drag coefficient. The tuned DNNs show excellent reliability and accuracy.
    
    \item The multi-criteria decision analysis method was used to determine the optimal \addedRIII{MO} design. The proposed optimum provided an \replacedRIII{improvement}{increase} of 1.5\% and 17\% for the pressure recovery factor and \deletedRIII{the} drag coefficient, respectively. \addedRIII{Constrained designs achieved improvements of 1\% and 15\% for $C_p$ and $C_d$. The main cause is the change in the draft tube's elbow. This has a strong influence on the velocity and thus on the total pressure. Although the change in the cross-sectional area leads to an increase in the static pressure drop in the draft tube, the effect is not as significant.}
    
    \item \deletedRIII{We have demonstrated the effectiveness of the L-SHADE and MOEA/D algorithms for draft tube design optimisation.} Depending on the goals, it is valid to consider single-objective optimisation. \addedRIII{MOEA/D is consistently the best-performing MO algorithm, while L-SHADE is the best SO algorithm.}
    
    \deletedRI{The proposed methodology is evaluated at a single operating point. Furthermore, flow simplifications were introduced. Nevertheless, a foundation is laid for future research that would consider velocity profiles and behaviour at multiple points. The workflow is still valid in this context.}
\end{itemize}

The proposed methodology has the potential to improve the effectiveness, environmental sustainability, and operational lifetime of hydropower turbines. As the demand for renewable energy increases, the inclusion of data-driven optimisation techniques such as the one presented here \replacedRIII{will be}{becomes} increasingly important to achieve a sustainable energy future. \addedRI{Key aspects that should be investigated in future studies are:}
\begin{itemize}
    \item \addedRI{The methodology should be extended to include data for different operating conditions. By taking into account different flow regimes, a design that works efficiently under a range of operating conditions can be obtained. This is not feasible using conventional approaches.}
    
    \item \addedRI{The impact of velocity profiles on the design optimisation process should be investigated. To simplify the process (and due to a lack of validation data), this has been omitted from this study. Velocity profiles introduce variations in flow behaviour across the draft tube, directly influencing pressure distribution and consequently pressure recovery and drag.}

    \item \addedRI{A robustness assessment should be performed to determine how sensitive the optimised designs are to changes in input parameters and operating conditions. This will provide information on the reliability and stability of the proposed solutions and the overall methodology.}

    \item \addedRI{Cavitation effects should be taken into account. The collapse of vapour bubbles can create localised high-velocity jets and shock waves. Changes in flow patterns caused by cavitation can lead to fluctuations in the output of a hydropower system. Although this aspect does not affect the methodology, it should be included as a feature in the input data.}
    
    \item \addedRI{The applicability of advanced machine learning methods, e.g. reinforcement learning, should be investigated. By using a standardised approach to parameterisation (and data acquisition), a unified model could be trained that is flexible and applicable to different designs and conditions.}
\end{itemize}

\deletedRI{ Future studies should focus on investigating the influence of velocity profiles for realistic draft tube designs and cavitation prediction. In addition, multi-point optimisation should be considered to improve efficiency over a wider operating range. Finally, the approach should be extended to other turbine components to achieve comprehensive efficiency improvements.}

\addedRI{In summary, this study introduces a machine learning-based approach to enhance the efficiency and performance of draft tubes that can be extended to other turbine components. The implications go beyond theoretical advances and offer a practical approach to the revitalisation of existing systems.}

\appendix
\setcounter{figure}{0}
\section{}
\label{sec:appendix:a}

\begin{figure}[H]
    \centering
    \includegraphics[trim={7cm 1cm 7cm 1cm}, clip, width=\textwidth]{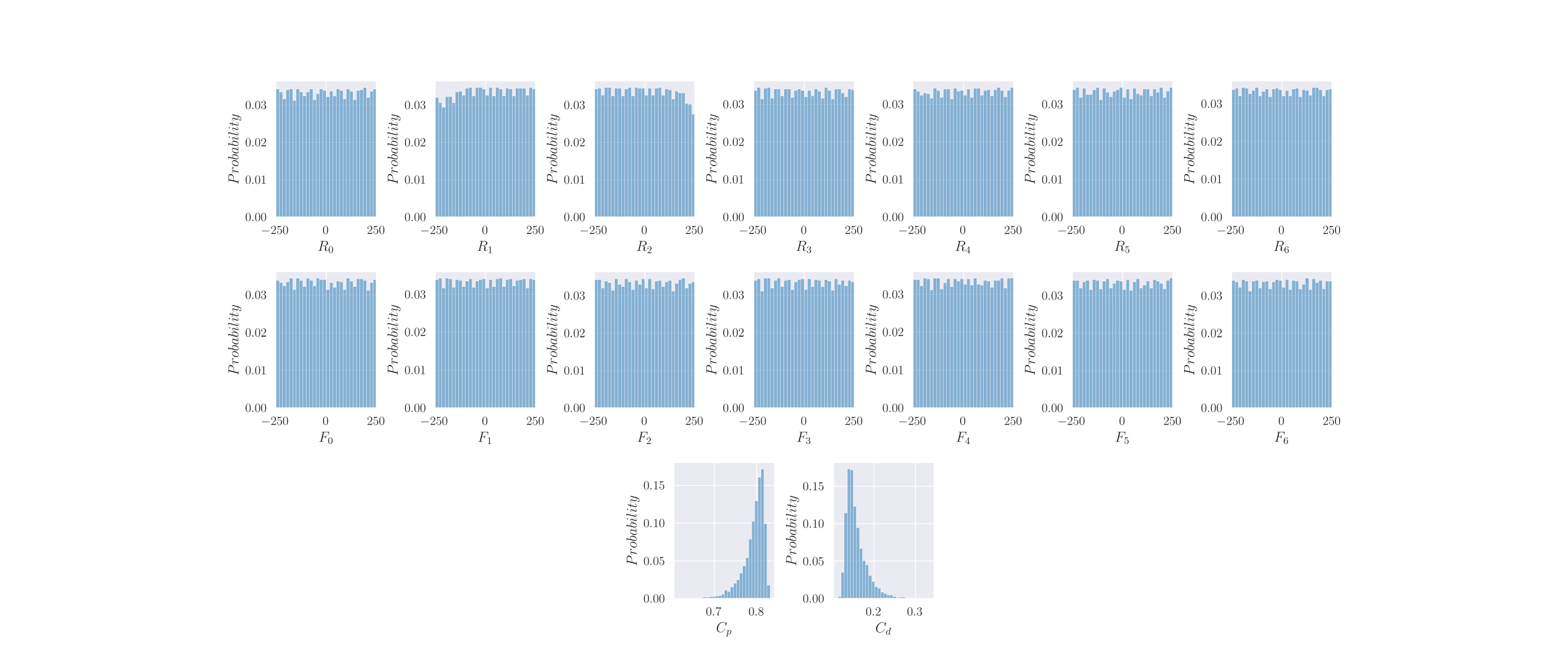}\vspace{-0.5cm}
    \caption{\protect\addedRI{Distribution of values for each feature and targets in the dataset used in Scenario I.}}
    \label{fig:distribution_fixed}
\end{figure}

\begin{figure}[H]
    \centering
    \includegraphics[trim={7cm 1cm 7cm 1cm}, clip, width=\textwidth]{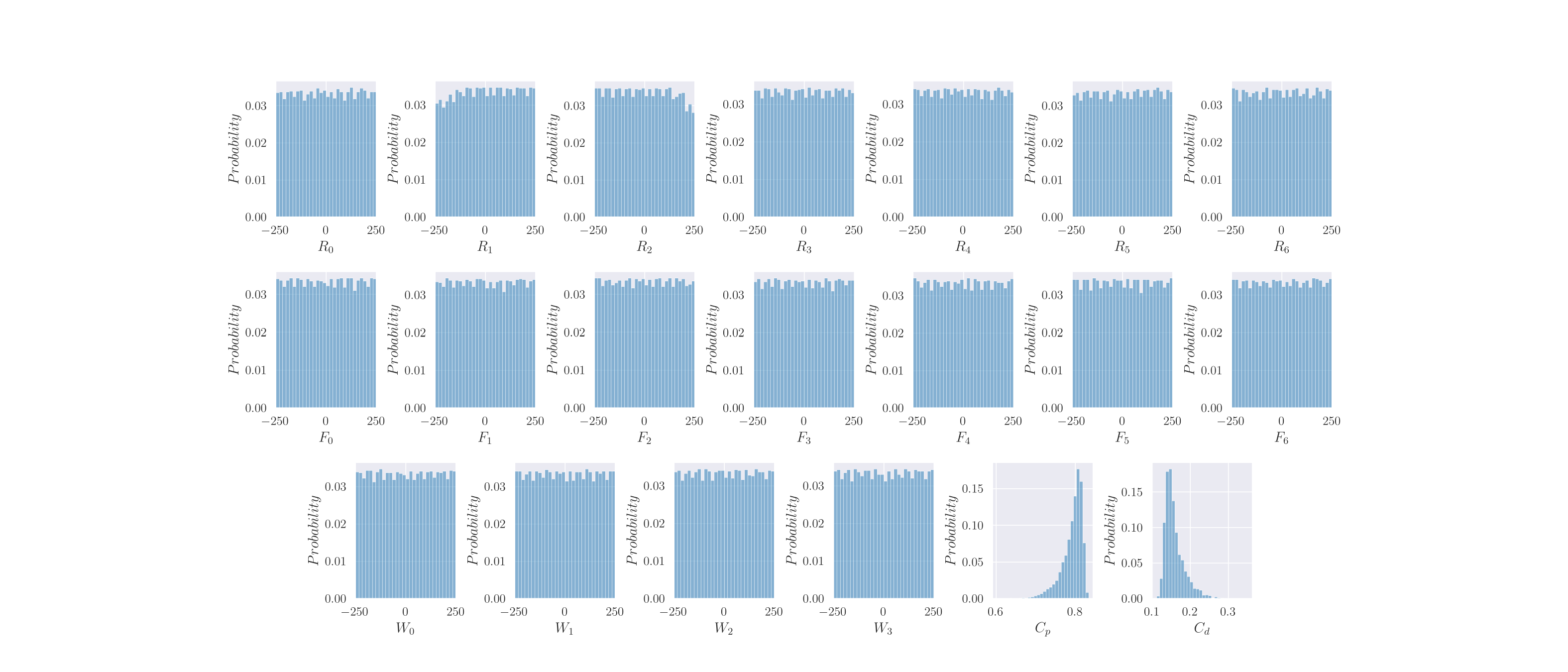}\vspace{-0.5cm}
    \caption{\protect\addedRI{Distribution of values for each feature and targets in the dataset used in Scenario II.}}
    \label{fig:distribution_free}
\end{figure}

\begin{figure}[H]
\centering
\begin{subfigure}[b]{0.495\textwidth}
    \centering
    \includegraphics[width=0.5\textwidth]{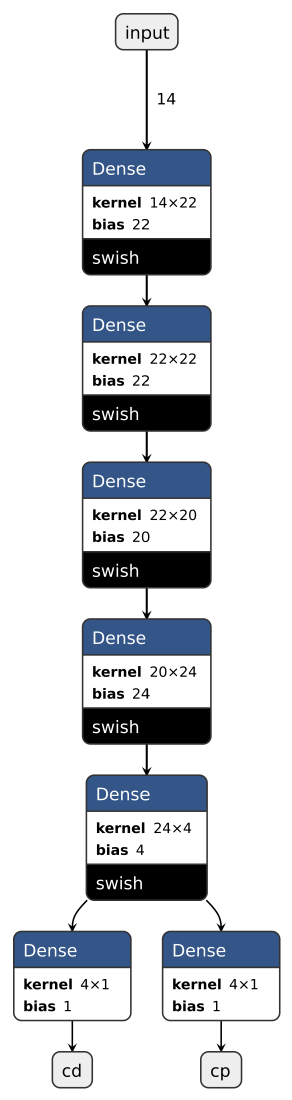}
    \caption{}
    \label{fig:DNN:fixed}
\end{subfigure}
\begin{subfigure}[b]{0.495\textwidth}
    \centering
    \includegraphics[width=0.5\textwidth]{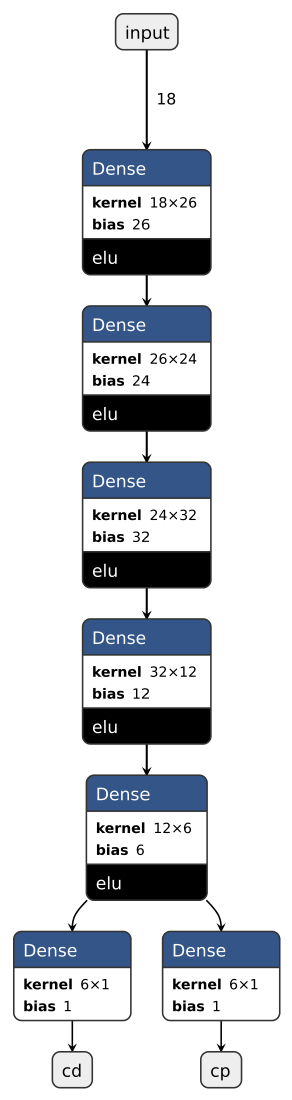}
    \caption{}
    \label{fig:DNN:free}
\end{subfigure}
\caption{Topology of the tuned DNNs for Scenario I (a) and Scenario II (b). Both networks have an input layer, five hidden layers and an output layer. Networks are multi-output regressors and predict $C_p$ and $C_d$.}
\label{fig:DNN}
\end{figure}

\begin{figure}[H]
\centering
\begin{subfigure}[b]{0.325\textwidth}
    \centering
    \includegraphics[trim={6cm 6cm 6cm 0cm}, clip, width=\textwidth]{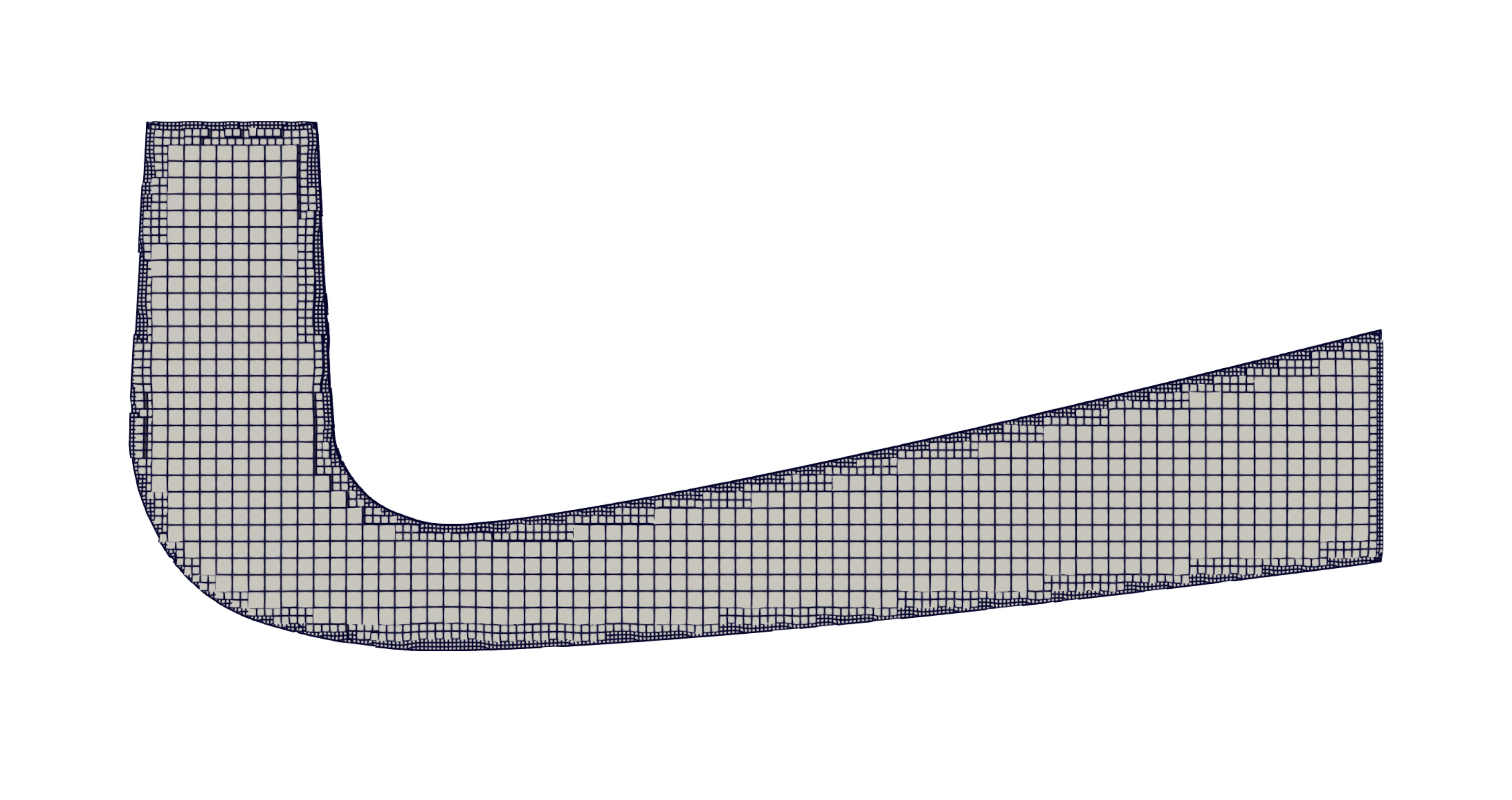}
    \label{fig:geom_mesh:o}
\end{subfigure}
\begin{subfigure}[b]{0.325\textwidth}
    \centering
    \includegraphics[trim={6cm 6cm 6cm 0cm}, clip, width=\textwidth]{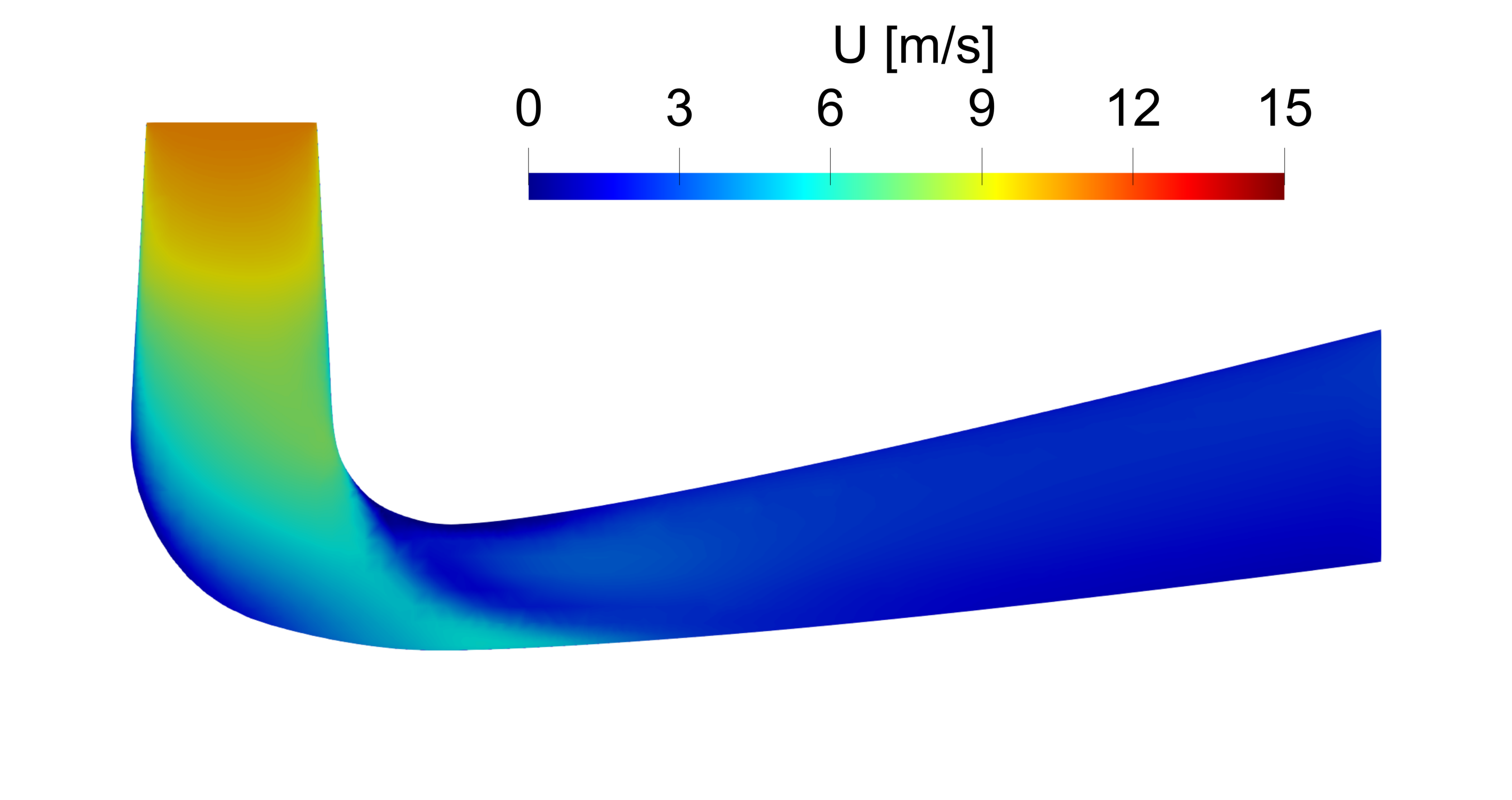}
    \label{fig:velocity:o}
\end{subfigure}
\begin{subfigure}[b]{0.325\textwidth}
    \centering
    \includegraphics[trim={6cm 6cm 6cm 0cm}, clip, width=\textwidth]{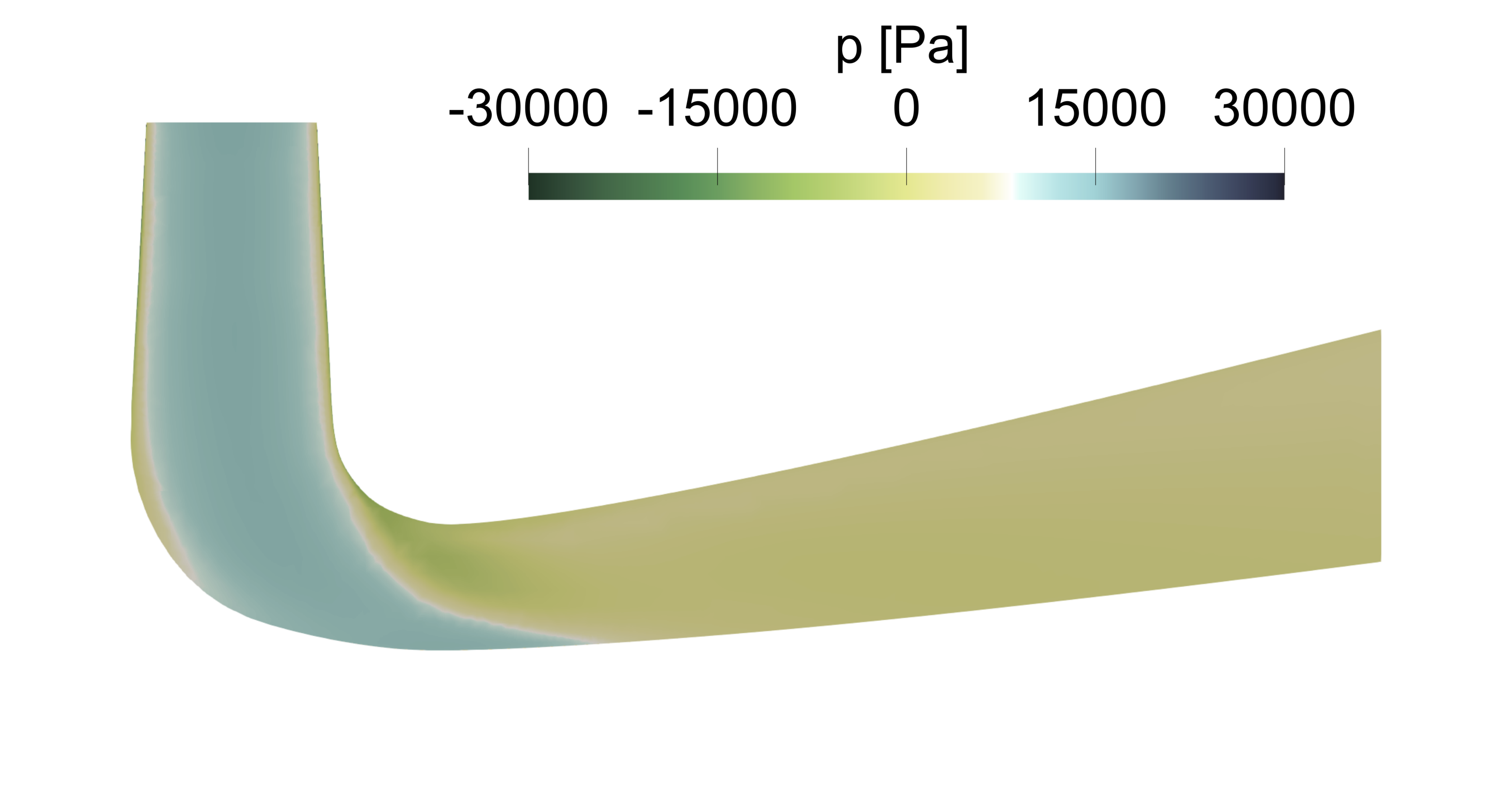}
    \label{fig:pressure:o}
\end{subfigure} \\\vspace{-1cm}

\begin{subfigure}[b]{0.325\textwidth}
    \centering
    \includegraphics[trim={6cm 6cm 6cm 0cm}, clip, width=\textwidth]{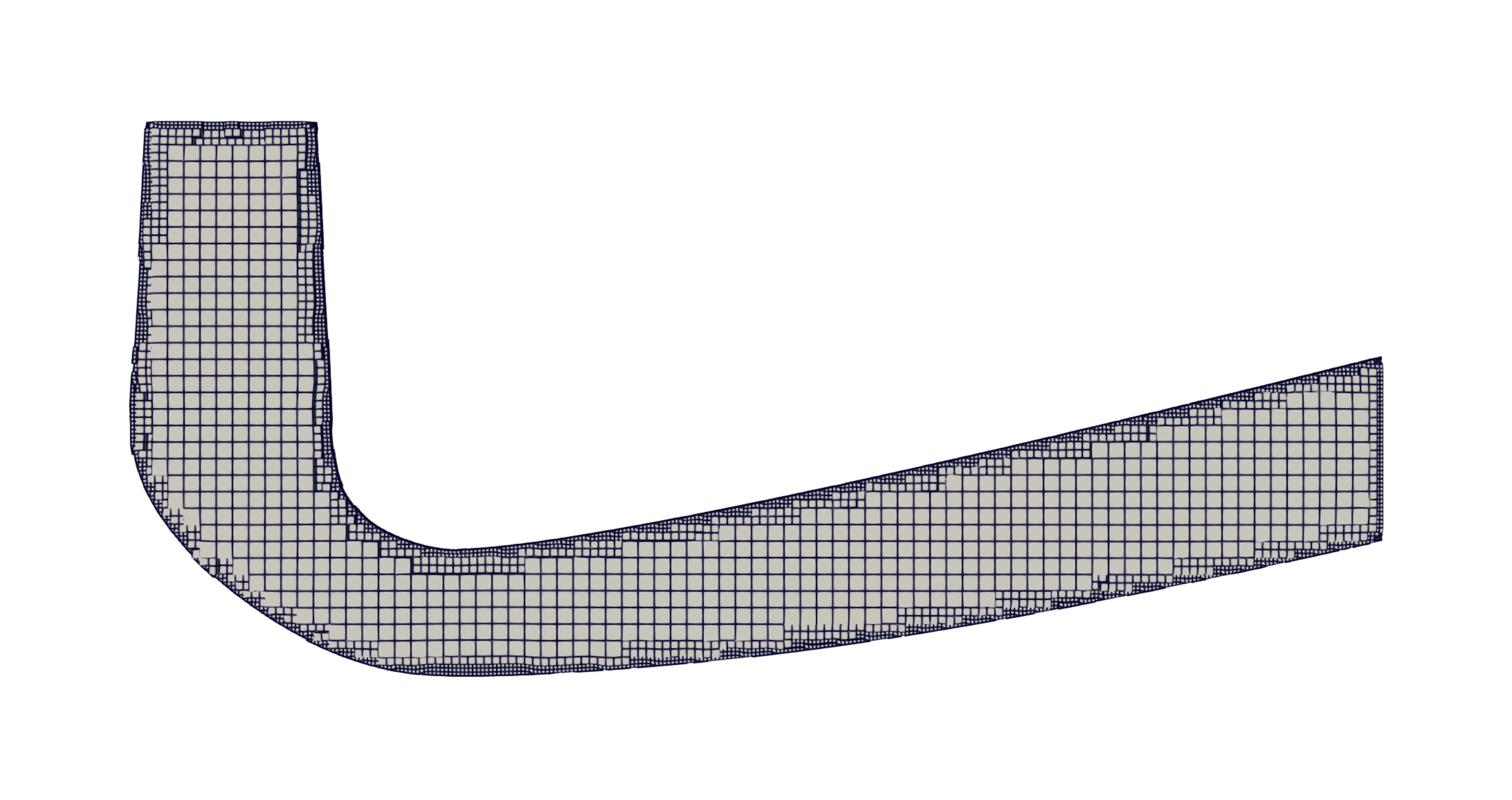}
    \label{fig:geom_mesh:g_fixed}
\end{subfigure}
\begin{subfigure}[b]{0.325\textwidth}
    \centering
    \includegraphics[trim={6cm 6cm 6cm 0cm}, clip, width=\textwidth]{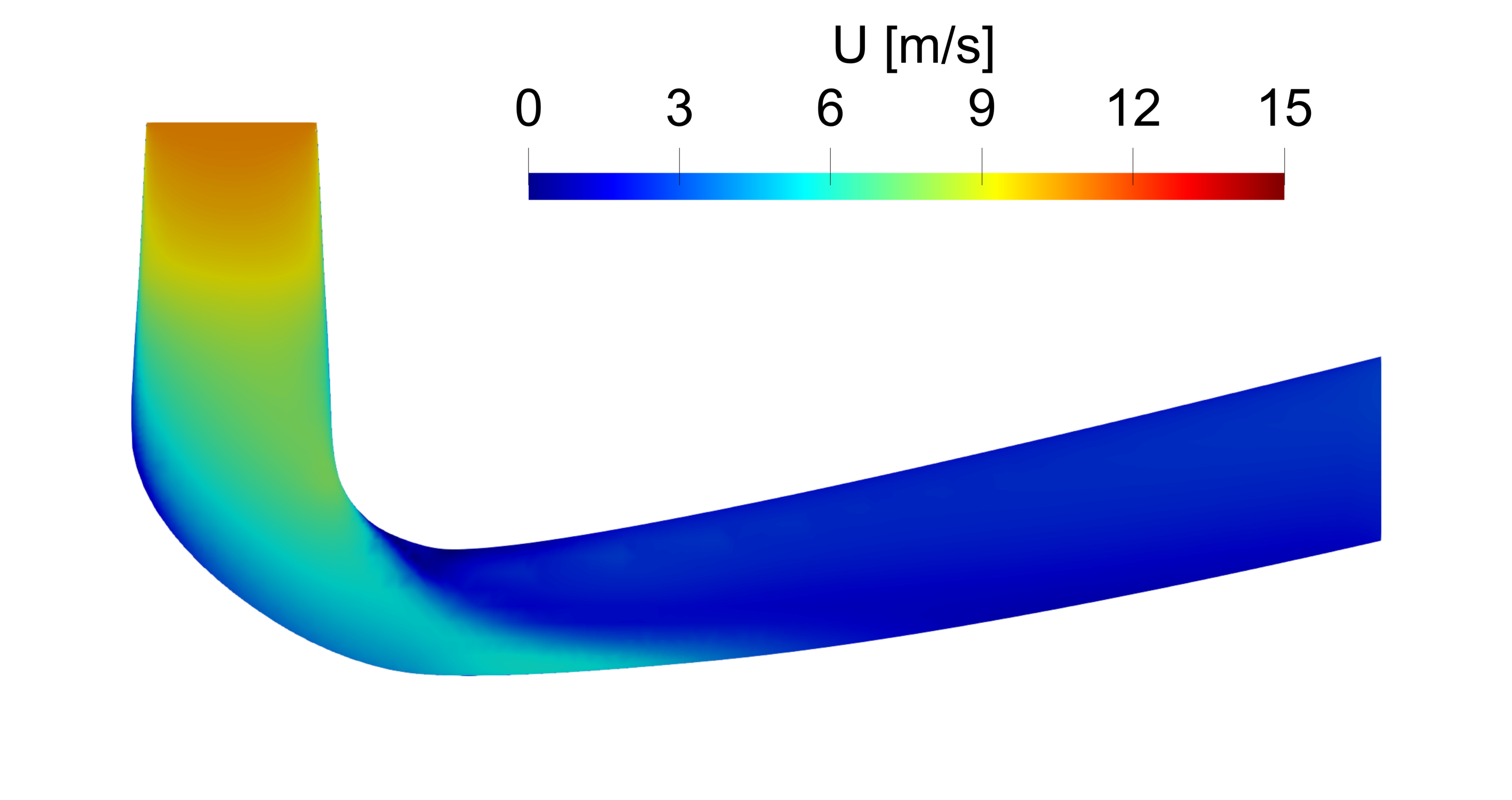}
    \label{fig:velocity:g_fixed}
\end{subfigure}
\begin{subfigure}[b]{0.325\textwidth}
    \centering
    \includegraphics[trim={6cm 6cm 6cm 0cm}, clip, width=\textwidth]{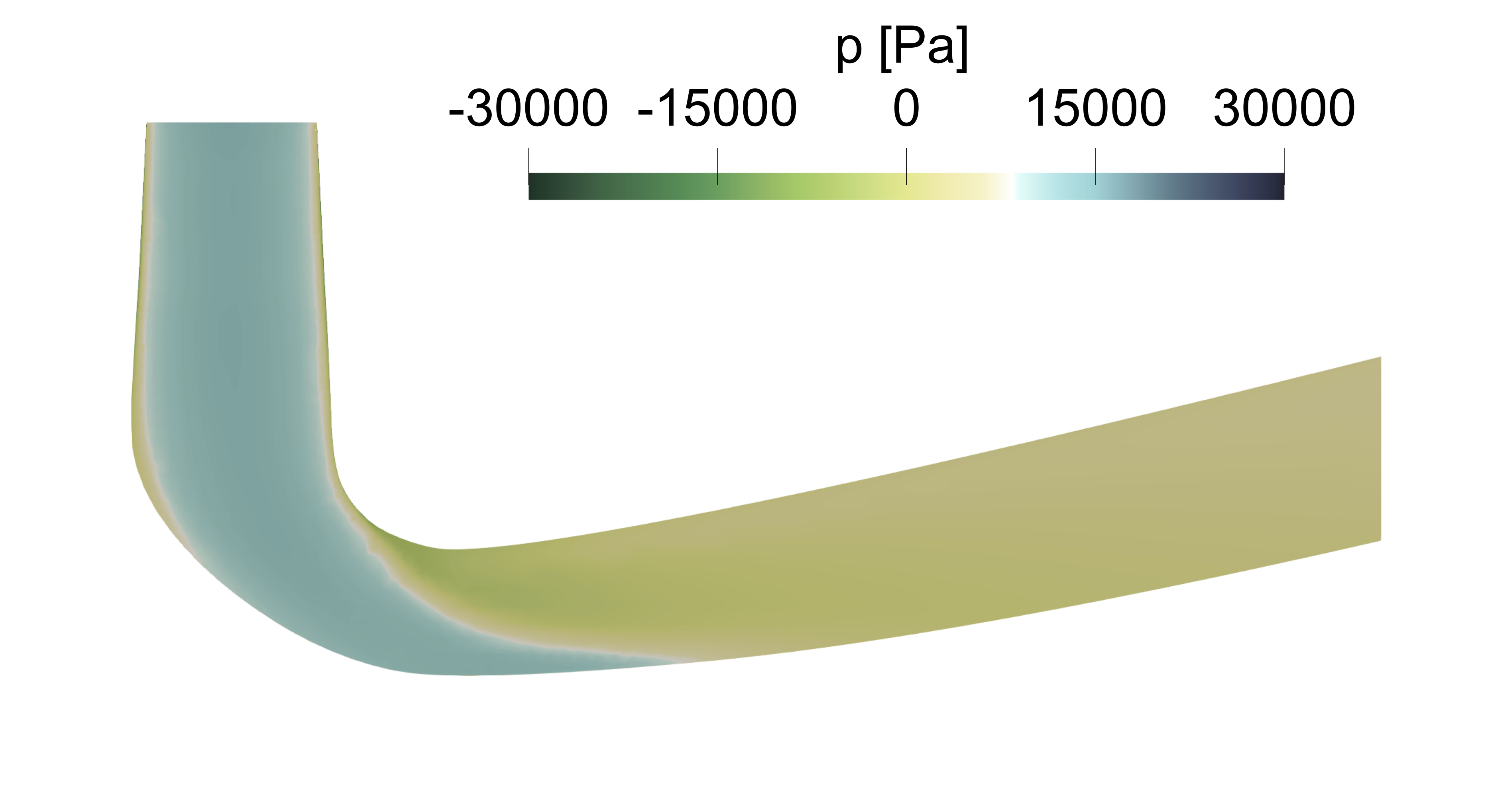}
    \label{fig:pressure:g_fixed}
\end{subfigure} \\\vspace{-1cm}

\begin{subfigure}[b]{0.325\textwidth}
    \centering
    \includegraphics[trim={6cm 6cm 6cm 0cm}, clip, width=\textwidth]{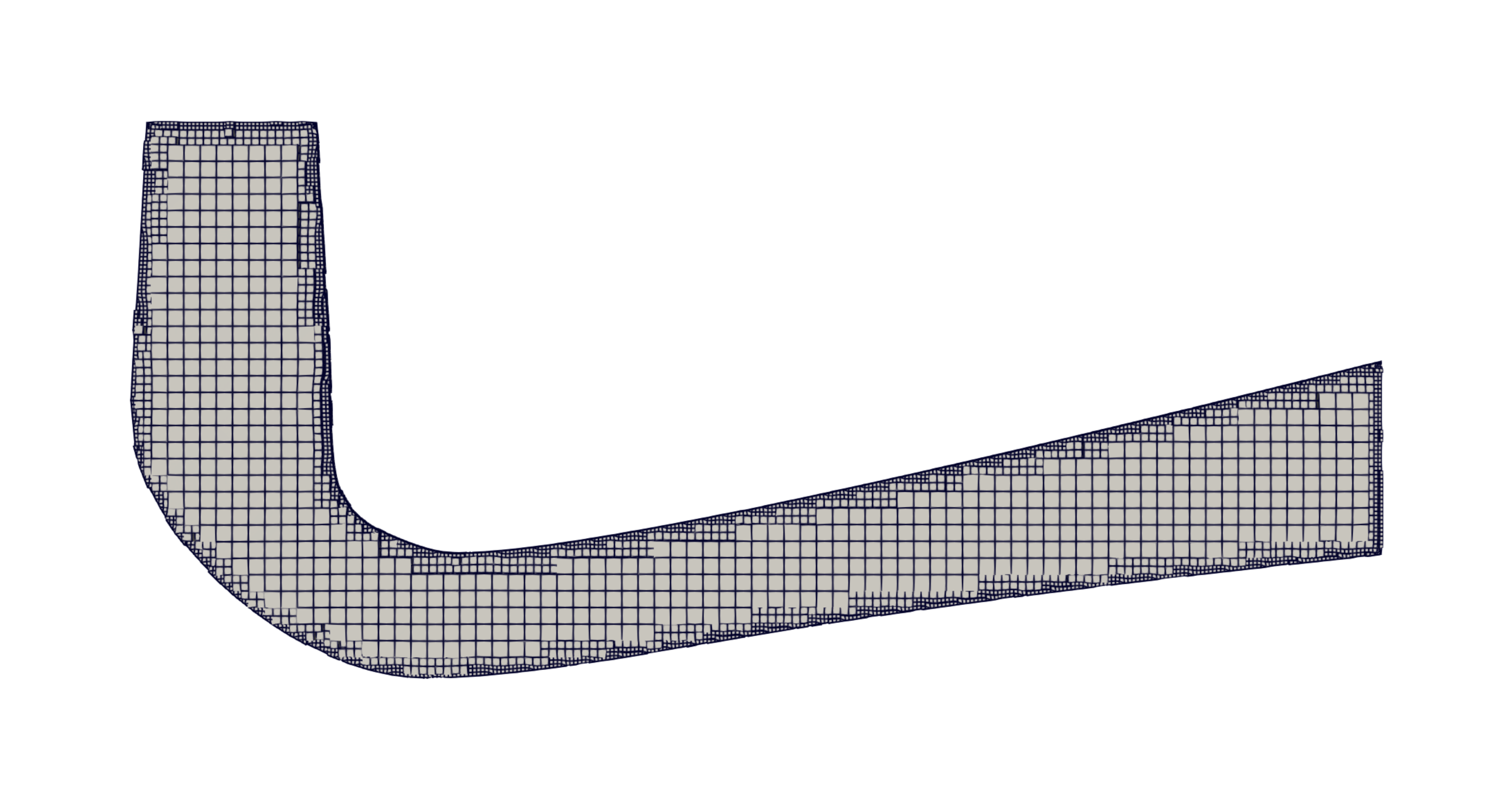}
    \label{fig:geom_mesh:g_free}
\end{subfigure}
\begin{subfigure}[b]{0.325\textwidth}
    \centering
    \includegraphics[trim={6cm 6cm 6cm 0cm}, clip, width=\textwidth]{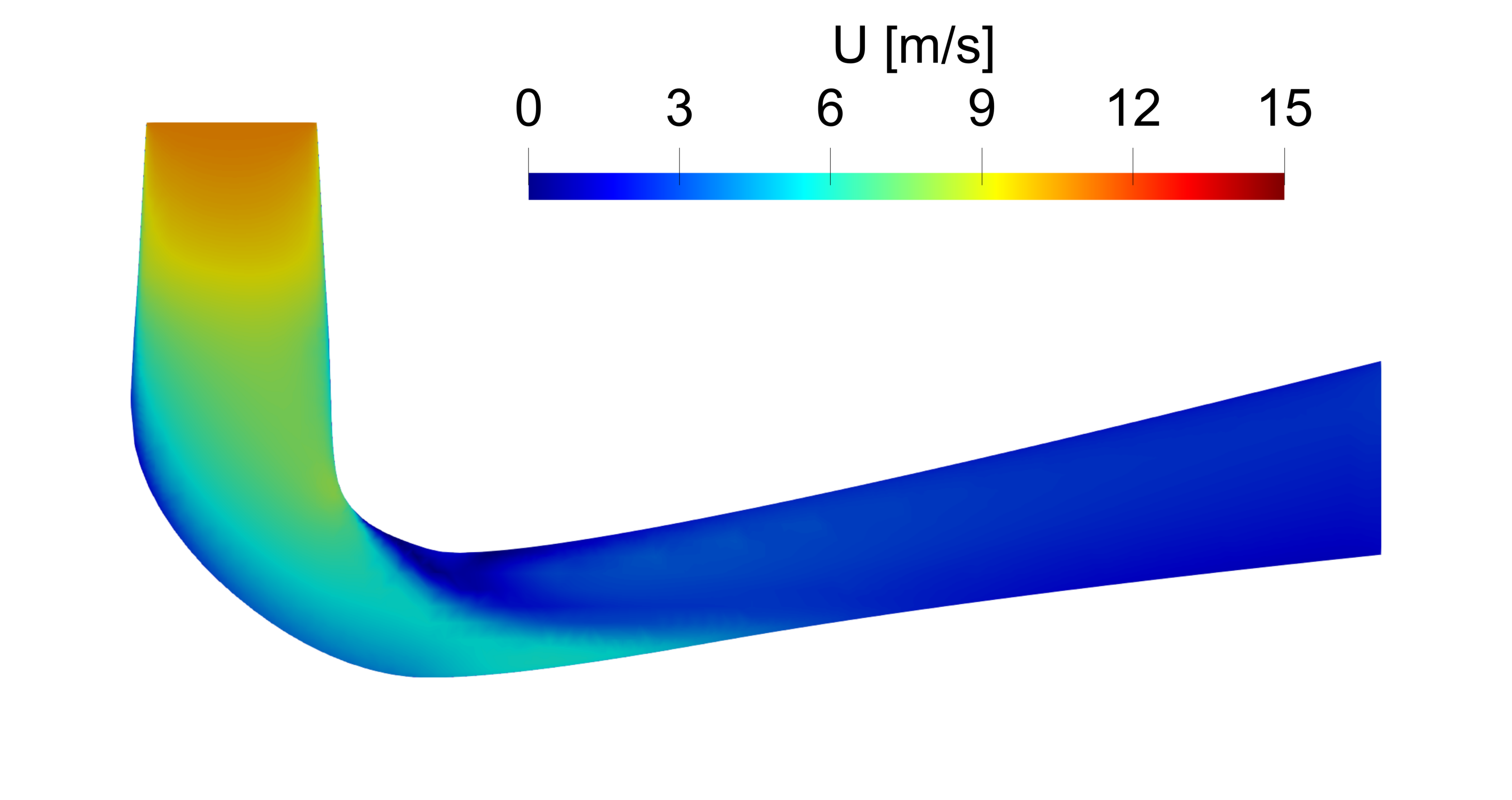}
    \label{fig:velocity:g_free}
\end{subfigure}
\begin{subfigure}[b]{0.325\textwidth}
    \centering
    \includegraphics[trim={6cm 6cm 6cm 0cm}, clip, width=\textwidth]{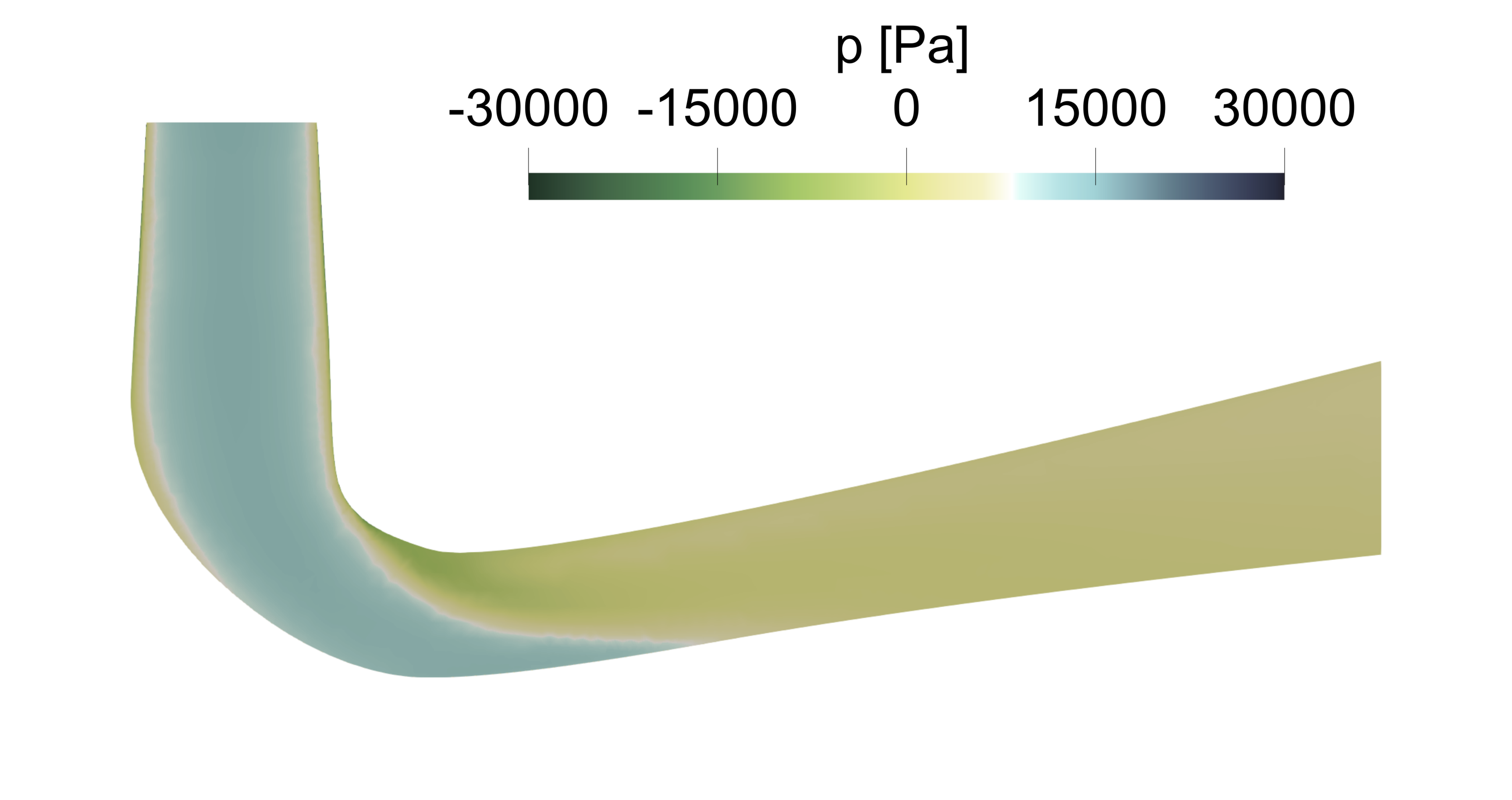}
    \label{fig:pressure:g_free}
\end{subfigure} \\\vspace{-1cm}

\begin{subfigure}[b]{0.325\textwidth}
    \centering
    \includegraphics[trim={6cm 6cm 6cm 0cm}, clip, width=\textwidth]{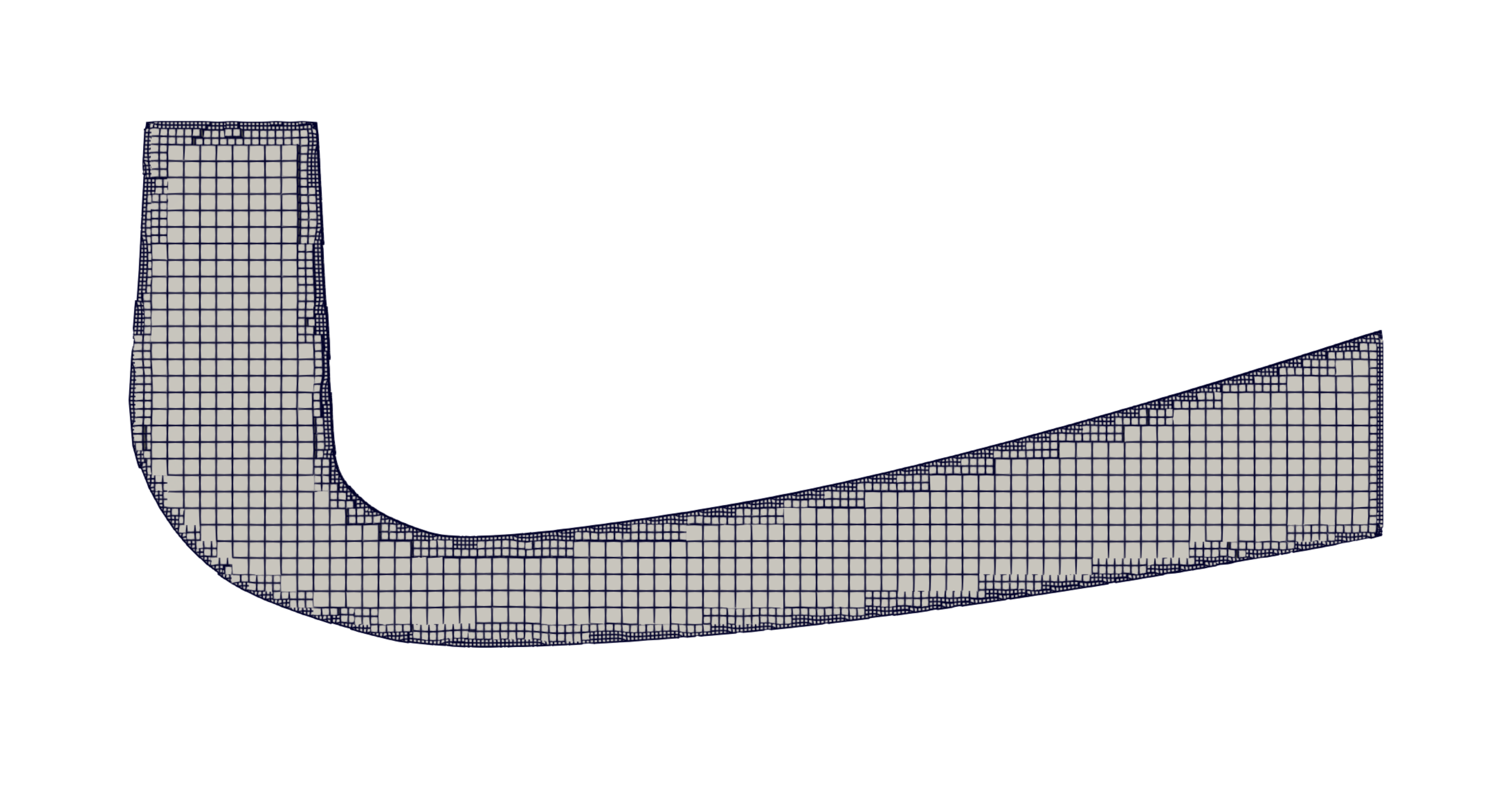}
    \label{fig:geom_mesh:i_fixed}
\end{subfigure}
\begin{subfigure}[b]{0.325\textwidth}
    \centering
    \includegraphics[trim={6cm 6cm 6cm 0cm}, clip, width=\textwidth]{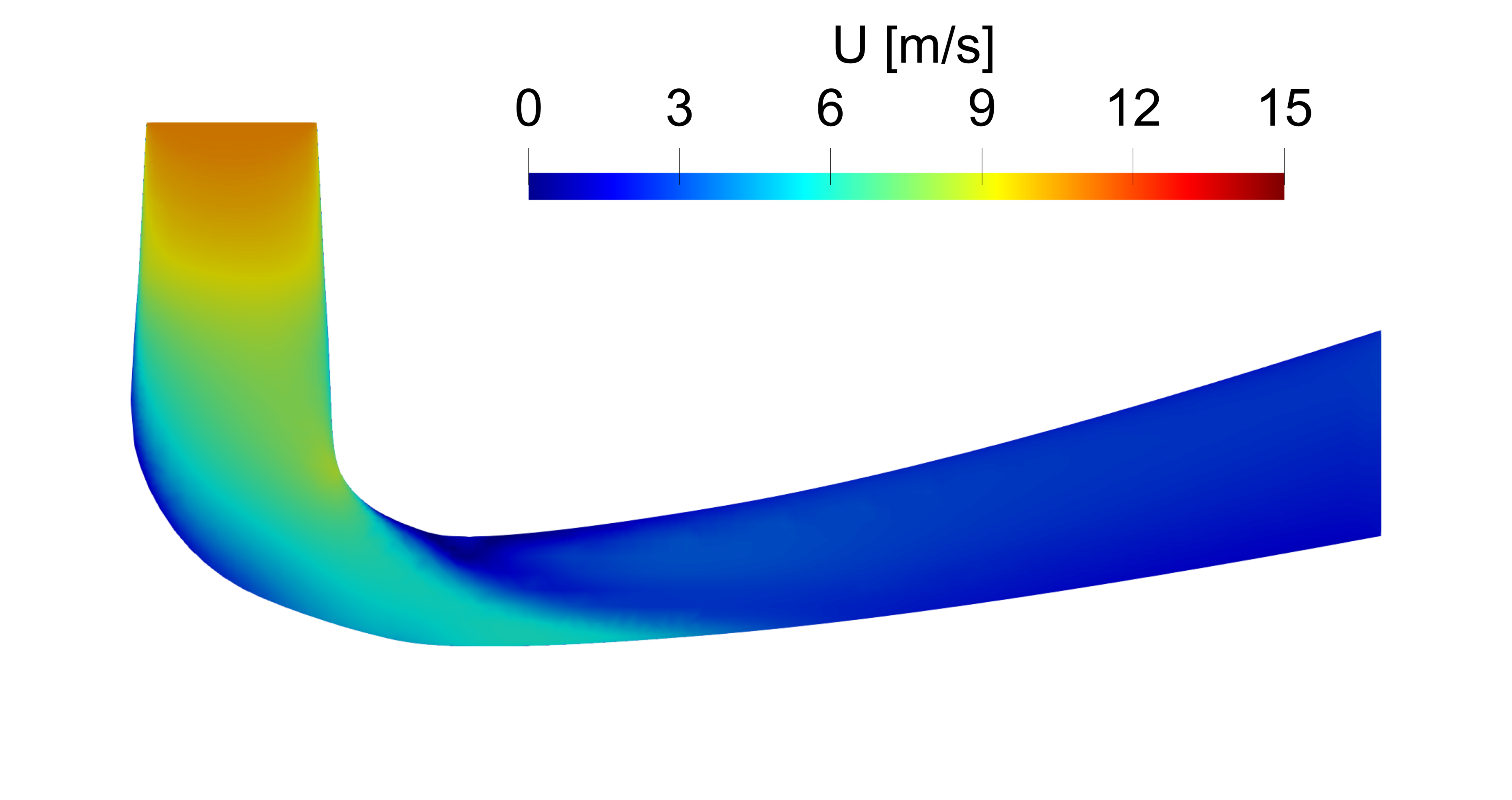}
    \label{fig:velocity:i_fixed}
\end{subfigure}
\begin{subfigure}[b]{0.325\textwidth}
    \centering
    \includegraphics[trim={6cm 6cm 6cm 0cm}, clip, width=\textwidth]{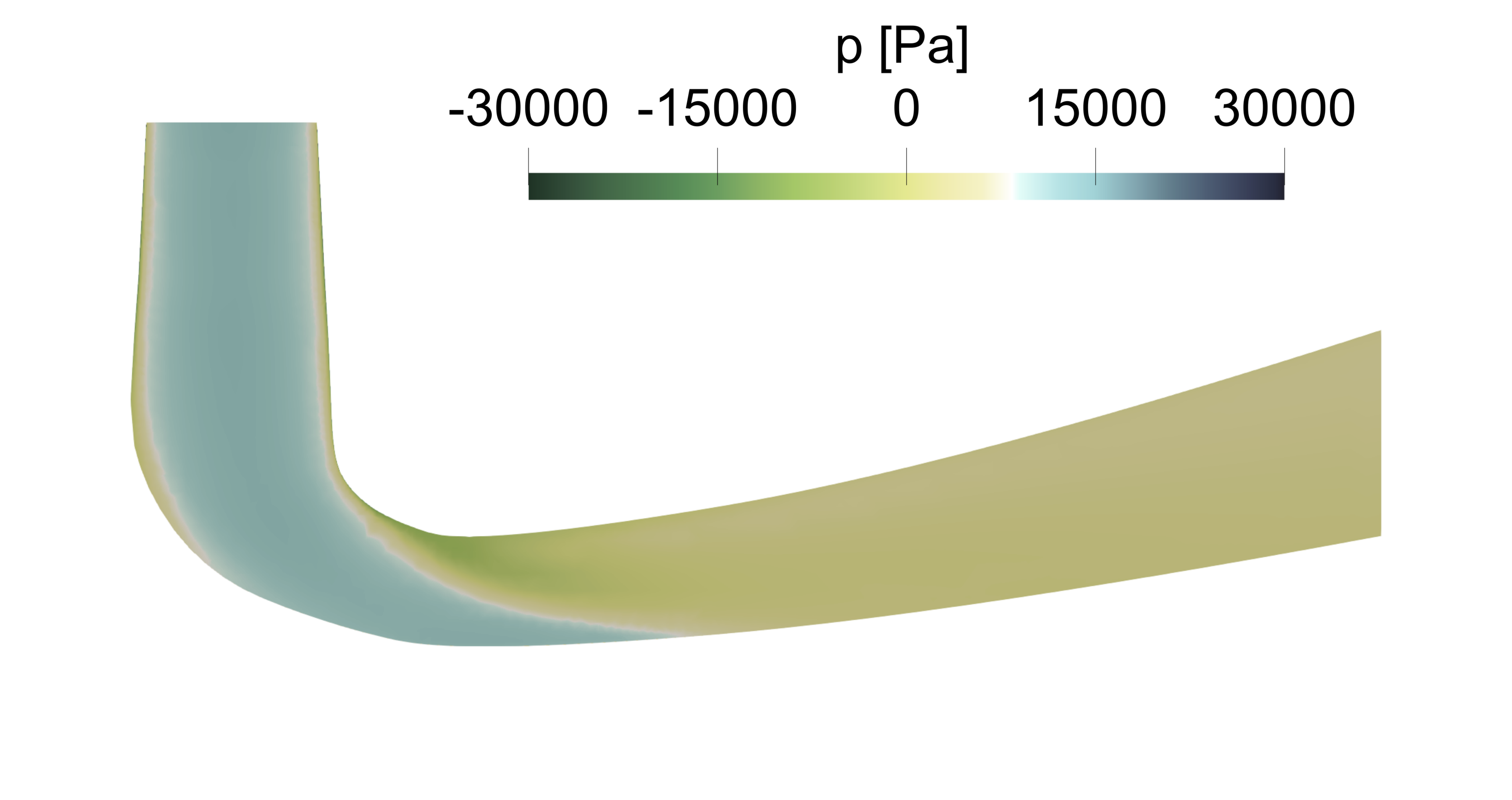}
    \label{fig:pressure:i_fixed}
\end{subfigure} \\\vspace{-1cm}

\begin{subfigure}[b]{0.325\textwidth}
    \centering
    \includegraphics[trim={6cm 6cm 6cm 0cm}, clip, width=\textwidth]{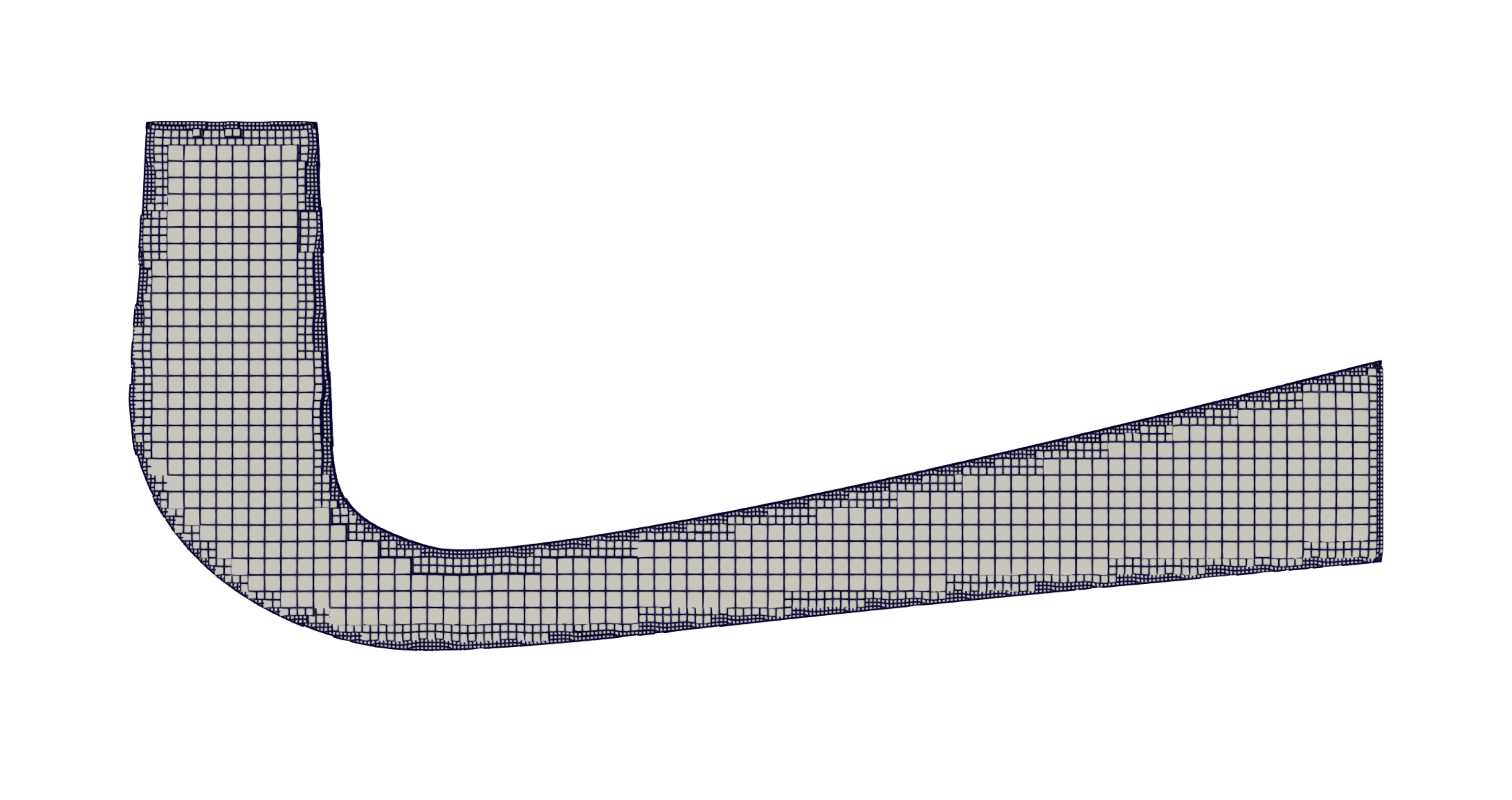}
    \caption{}
    \label{fig:geom_mesh:i_free}
\end{subfigure}
\begin{subfigure}[b]{0.325\textwidth}
    \centering
    \includegraphics[trim={6cm 6cm 6cm 0cm}, clip, width=\textwidth]{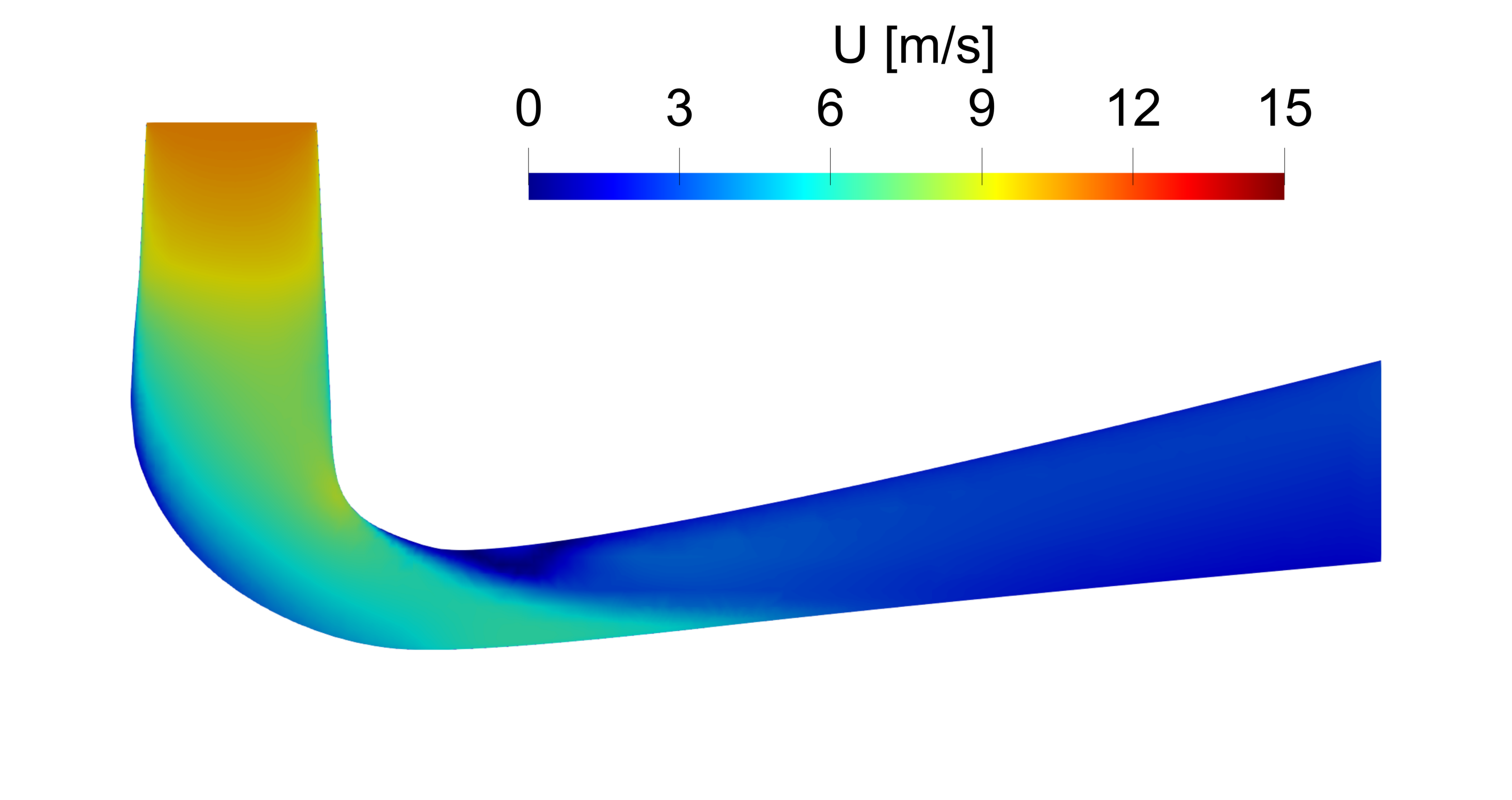}
    \caption{}
    \label{fig:velocity:i_free}
\end{subfigure}
\begin{subfigure}[b]{0.325\textwidth}
    \centering
    \includegraphics[trim={6cm 6cm 6cm 0cm}, clip, width=\textwidth]{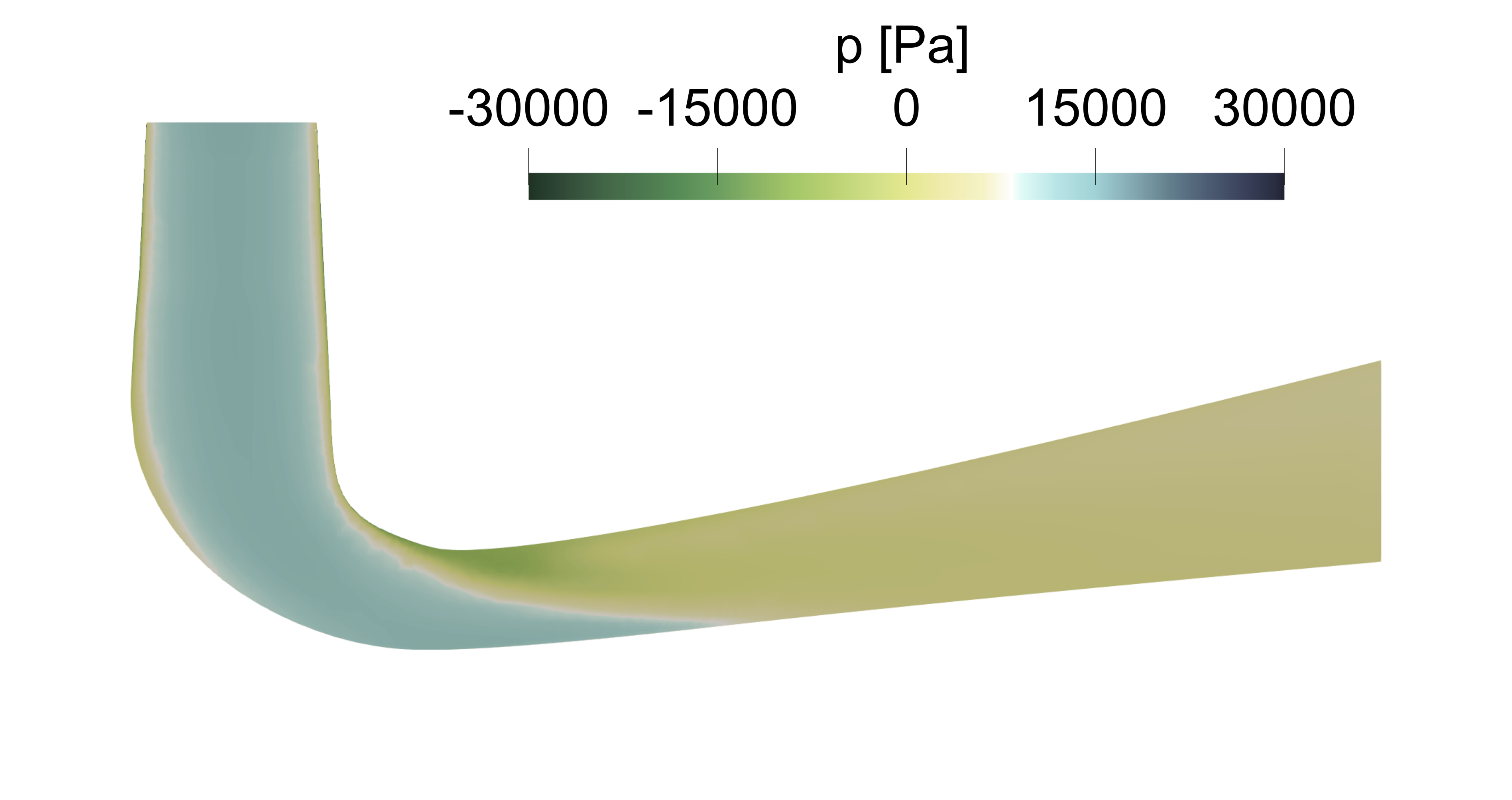}
    \caption{}
    \label{fig:pressure:i_free}
\end{subfigure}
\caption{Cross-sectional cuts ($xy$ plane) of the TOPSIS-proposed optimal designs for Scenario I.a, II.a, I.b and II.b (from top to bottom). The numerical grid (a), the flow field (b) and the pressure distribution (c) are shown.}
\label{fig:cfd_comparisson}
\end{figure}

\setcounter{figure}{0}
\section{}
\label{sec:appendix:b}

\begin{figure}[H]
\centering
\begin{subfigure}[b]{\textwidth}
    \centering
    \includegraphics[trim={0cm 0.5cm 0cm 1cm}, clip, width=0.9\textwidth]{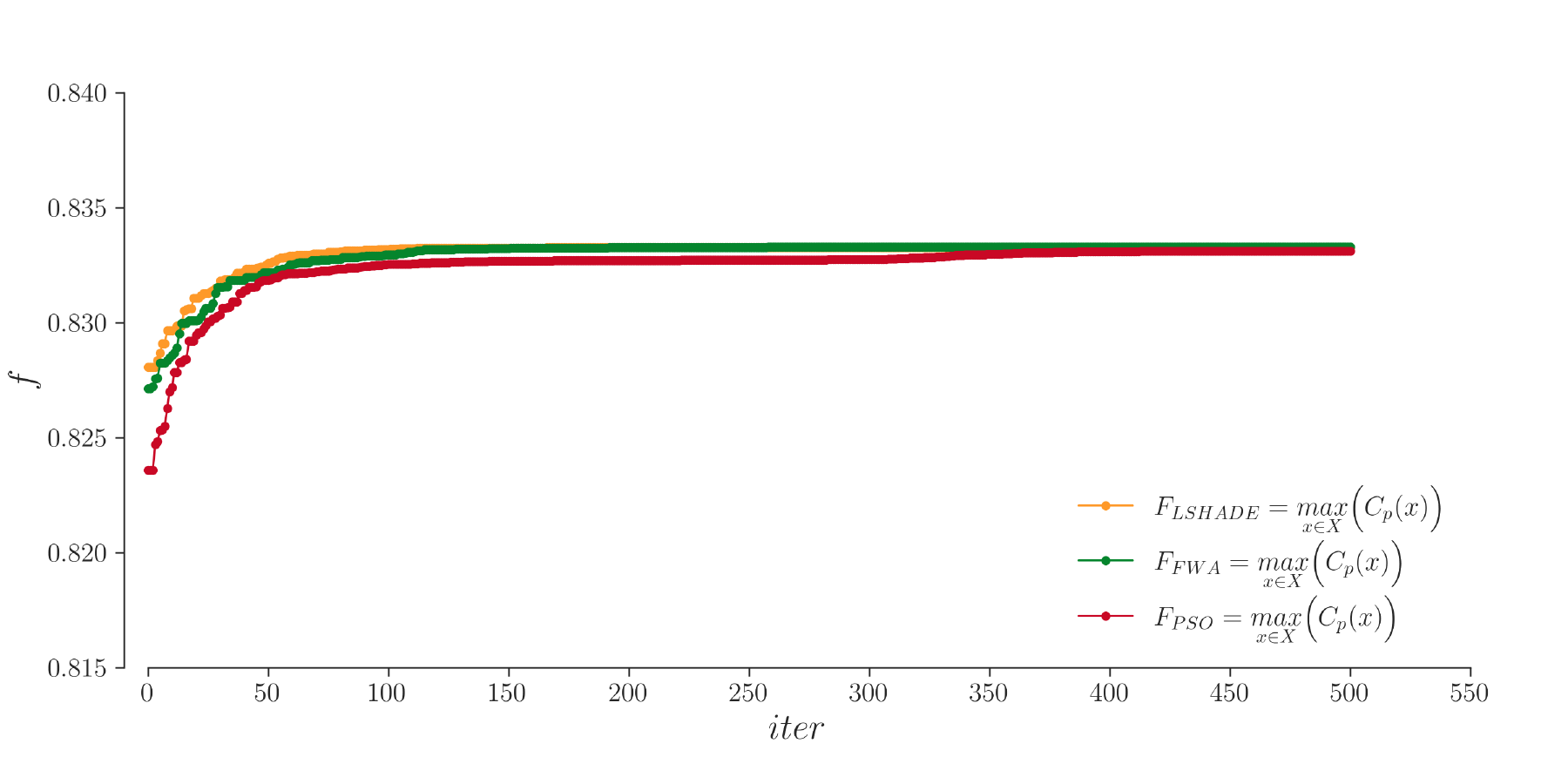}
    \caption{}
    \label{fig:cvg_scenarioIa:CP}
\end{subfigure} \\
\begin{subfigure}[b]{\textwidth}
    \centering
    \includegraphics[trim={0cm 0.5cm 0cm 1cm}, clip, width=0.9\textwidth]{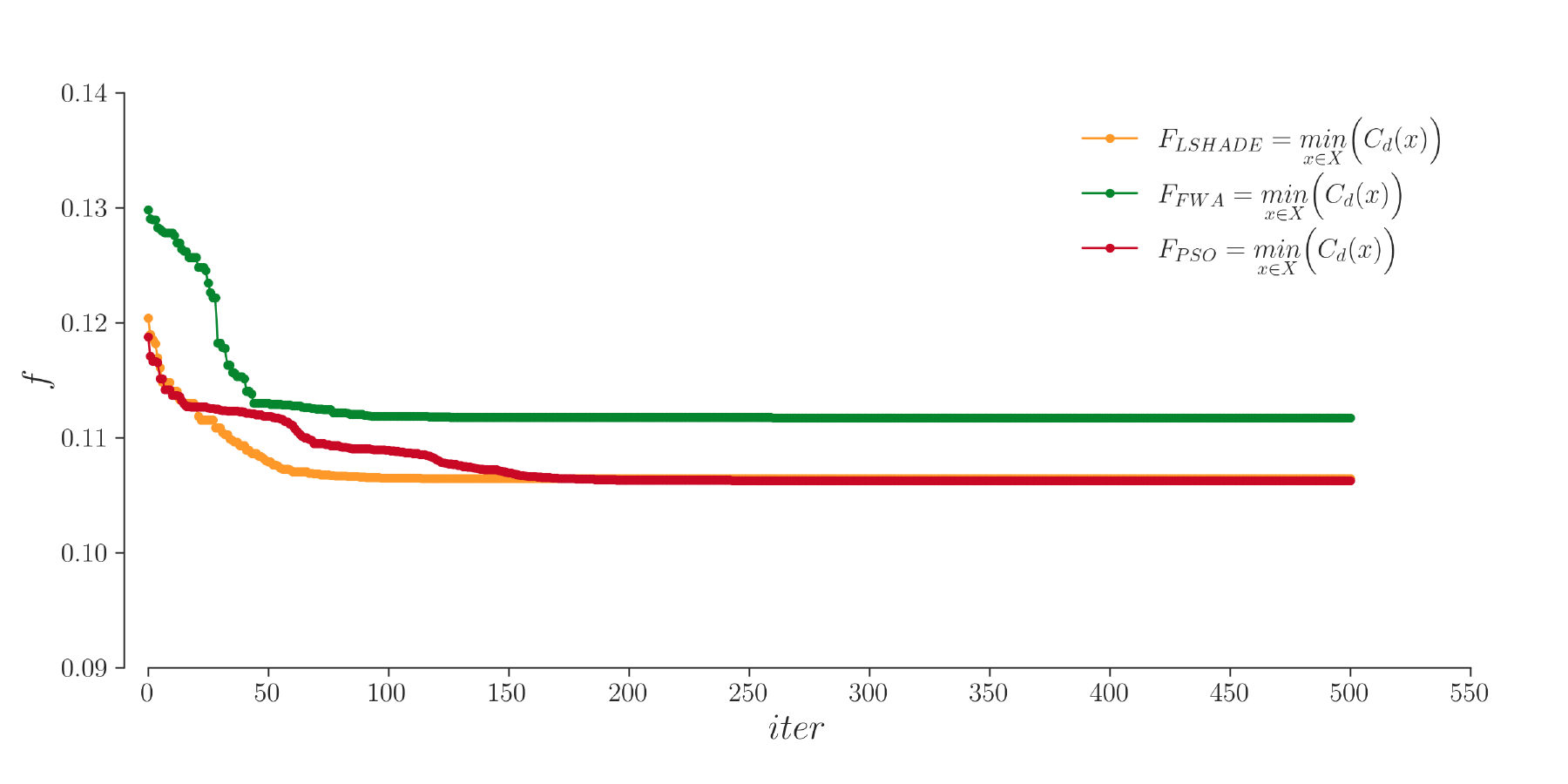}
    \caption{}
    \label{fig:cvg_scenarioIa:CD}
\end{subfigure}
\caption{Convergence graphs for FWA, PSO and L-SHADE when optimising for $C_p$ (a) and $C_d$ (b) in test Scenario I.a. L-SHADE performs better than the competition. The differences for $C_p$ are negligible, while for $C_d$, a difference of $\approx 5\%$ can be observed.}
\label{fig:cvg_scenarioIa}
\end{figure}

\begin{figure}[H]
\centering
\begin{subfigure}[b]{\textwidth}
    \centering
    \includegraphics[trim={0cm 0.5cm 0cm 1cm}, clip, width=0.9\textwidth]{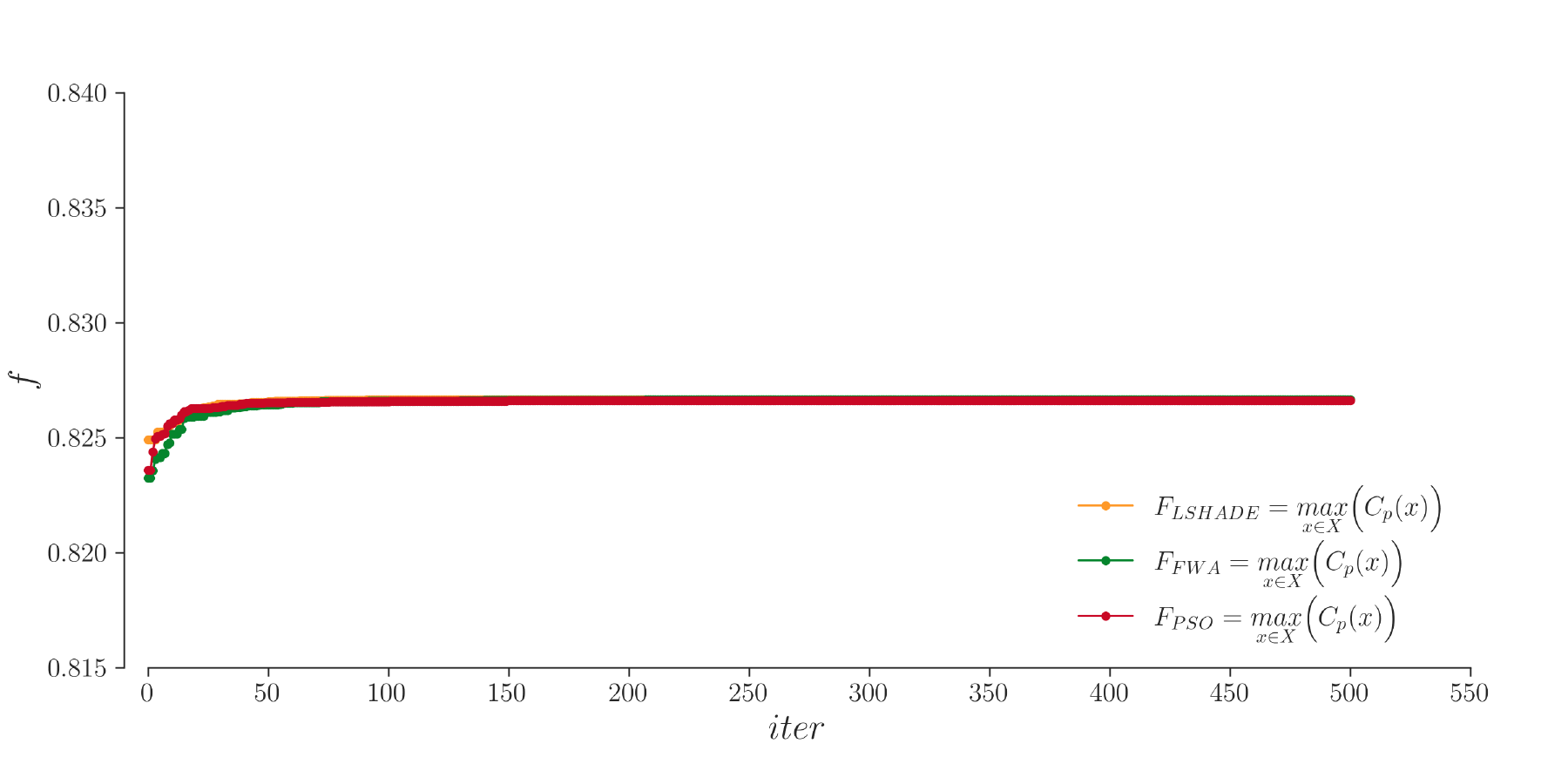}
    \caption{}
    \label{fig:cvg_scenarioIb:CP}
\end{subfigure} \\
\begin{subfigure}[b]{\textwidth}
    \centering
    \includegraphics[trim={0cm 0.5cm 0cm 1cm}, clip, width=0.9\textwidth]{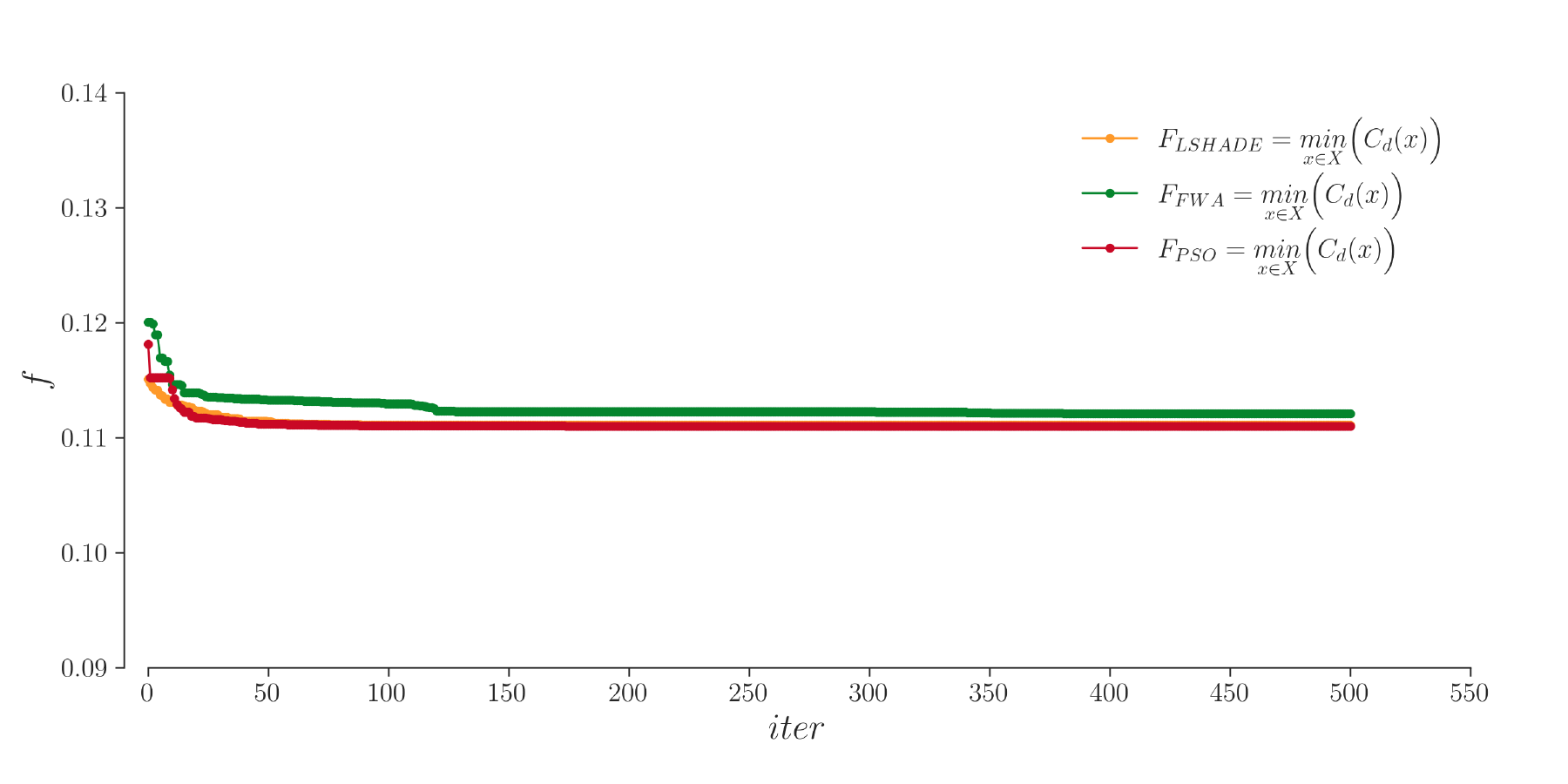}
    \caption{}
    \label{fig:cvg_scenarioIb:CD}
\end{subfigure}
\caption{Convergence graphs for FWA, PSO and L-SHADE when optimising for $C_p$ (a) and $C_d$ (b) in test Scenario I.b. L-SHADE performs better than the competition. FWA performs worst overall.}
\label{fig:cvg_scenarioIb}
\end{figure}

\begin{figure}[H]
\centering
\begin{subfigure}[b]{\textwidth}
    \centering
    \includegraphics[trim={0cm 0.5cm 0cm 1cm}, clip, width=0.9\textwidth]{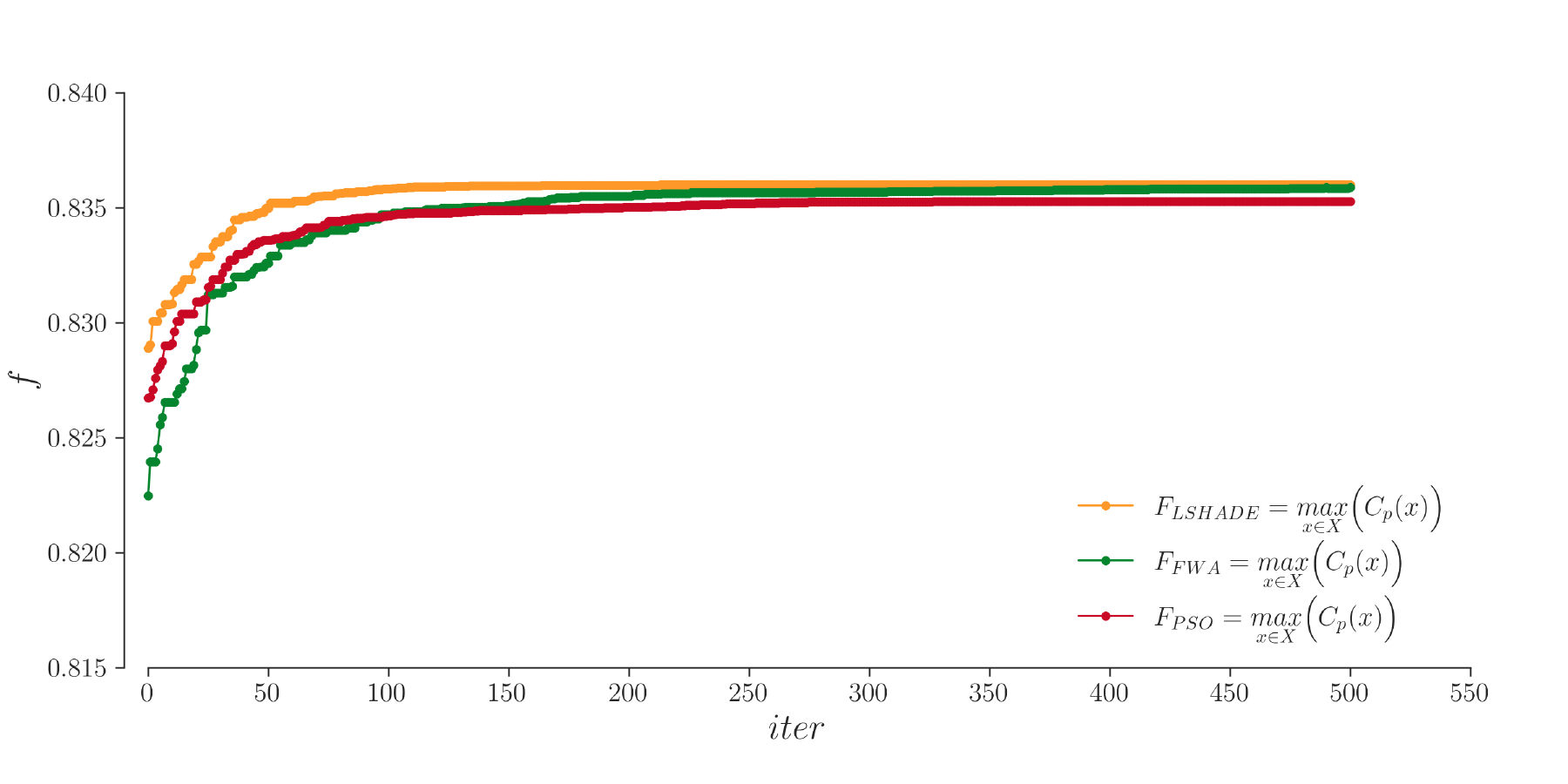}
    \caption{}
    \label{fig:cvg_scenarioIIa:CP}
\end{subfigure} \\
\begin{subfigure}[b]{\textwidth}
    \centering
    \includegraphics[trim={0cm 0.5cm 0cm 1cm}, clip, width=0.9\textwidth]{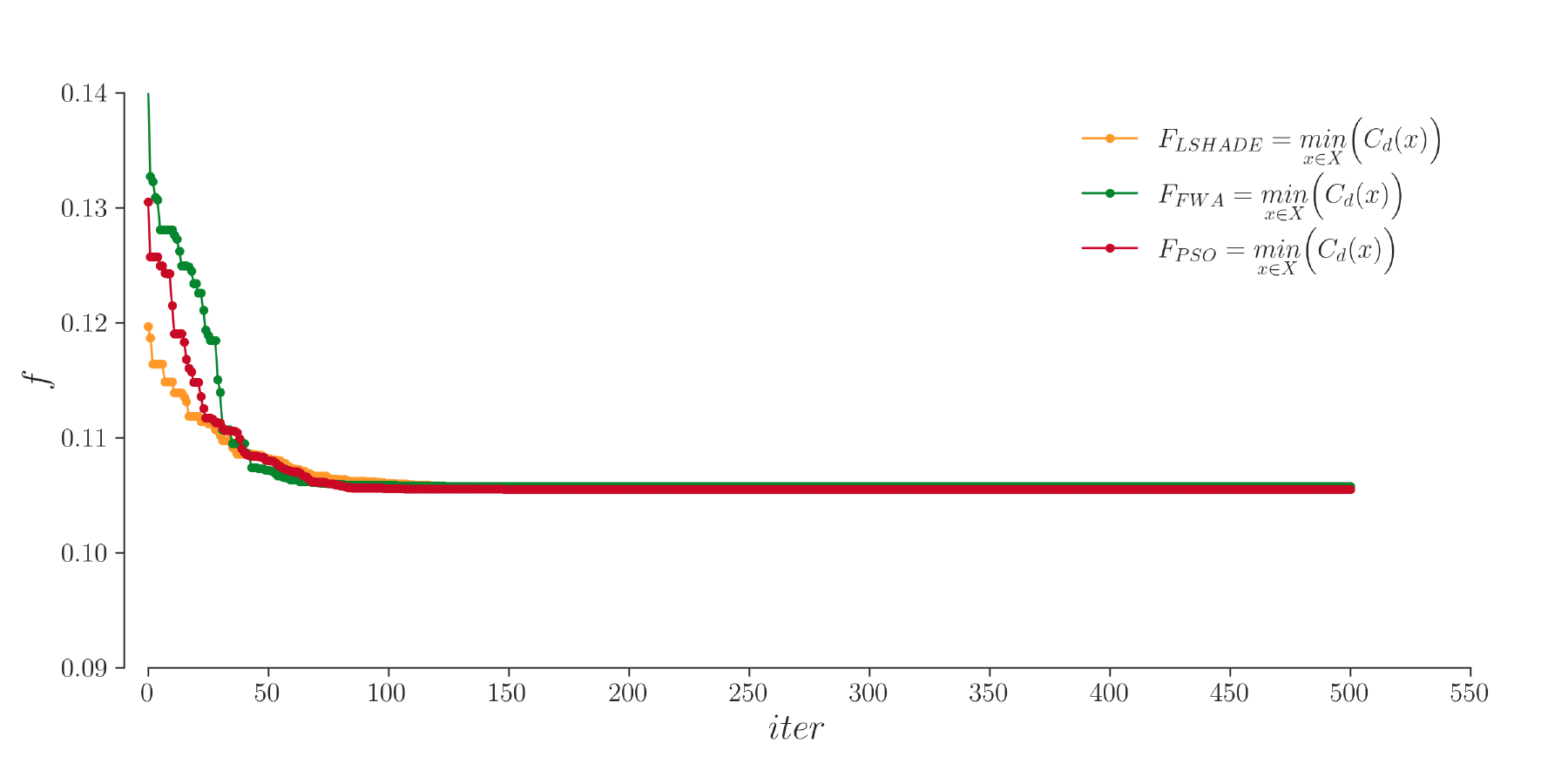}
    \caption{}
    \label{fig:cvg_scenarioIIa:CD}
\end{subfigure}
\caption{Convergence graphs for FWA, PSO and L-SHADE when optimising for $C_p$ (a) and $C_d$ (b) in test Scenario II.a. L-SHADE converges rapidly when optimising for $C_p$. All algorithms show acceptable convergence when optimising for $C_d$.}
\label{fig:cvg_scenarioIIa}
\end{figure}

\begin{figure}[H]
\centering
\begin{subfigure}[b]{\textwidth}
    \centering
    \includegraphics[trim={0cm 0.5cm 0cm 1cm}, clip, width=0.9\textwidth]{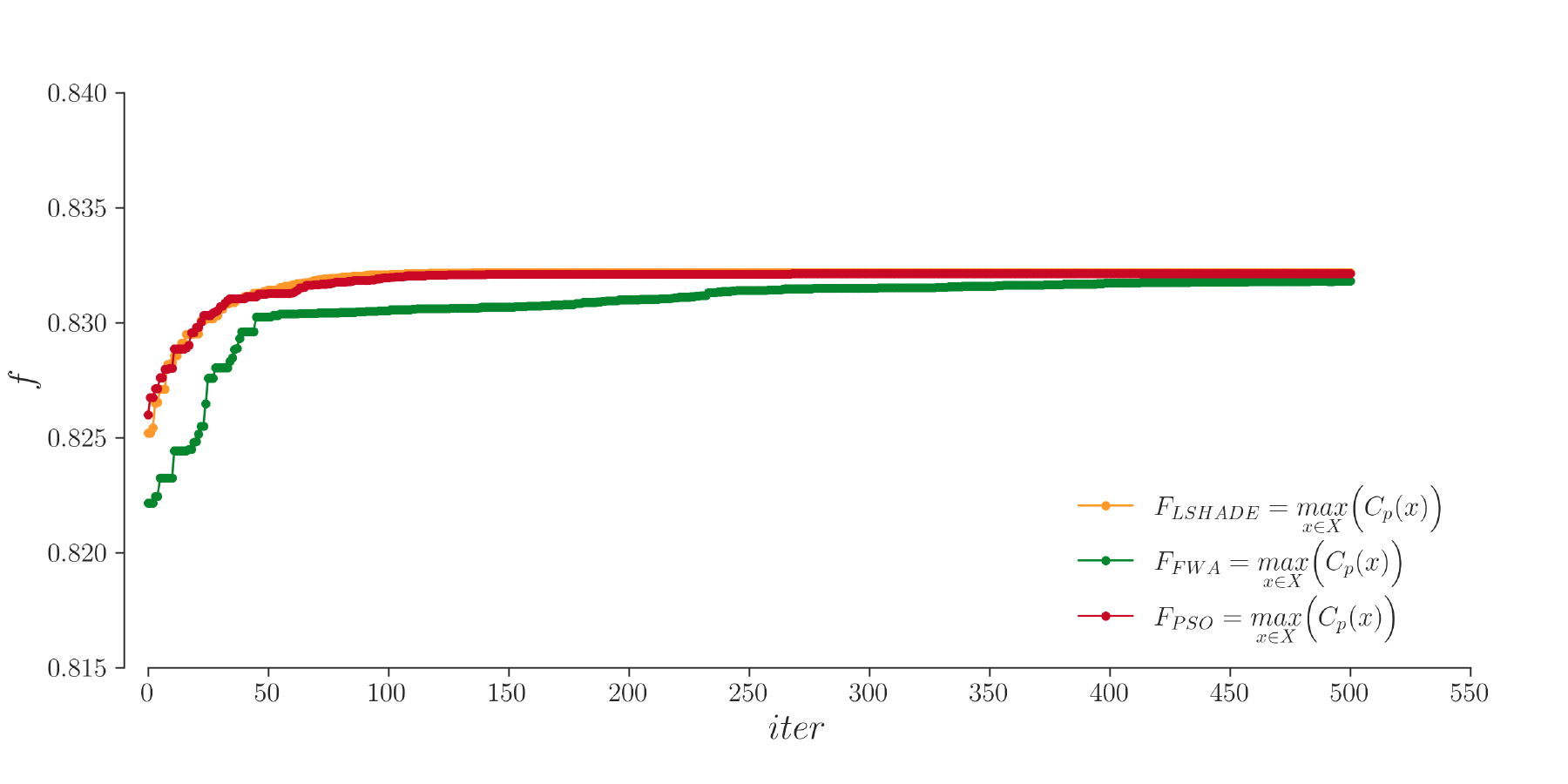}
    \caption{}
    \label{fig:cvg_scenarioIIb:CP}
\end{subfigure} \\
\begin{subfigure}[b]{\textwidth}
    \centering
    \includegraphics[trim={0cm 0.5cm 0cm 1cm}, clip, width=0.9\textwidth]{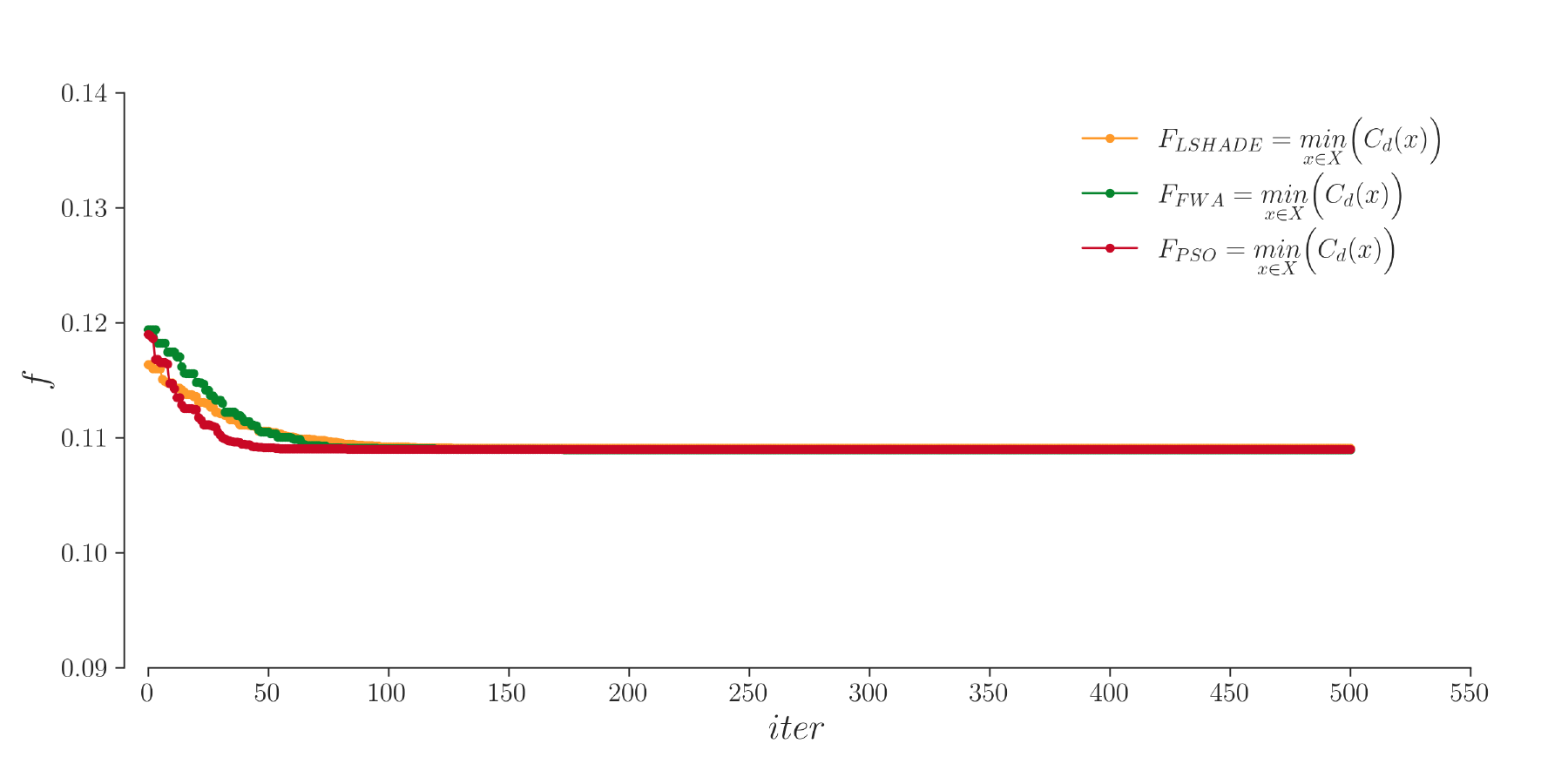}
    \caption{}
    \label{fig:cvg_scenarioIIb:CD}
\end{subfigure}
\caption{Convergence graphs for FWA, PSO and L-SHADE when optimising for $C_p$ (a) and $C_d$ (b) in test Scenario II.b. L-SHADE and PSO perform similarly, with PSO outperforming L-SHADE when optimising for $C_d$. FWA struggles with convergence for $C_p$.}
\label{fig:cvg_scenarioIIb}
\end{figure}

\begin{figure}[H]
    \centering
    \includegraphics[width=0.9\textwidth]{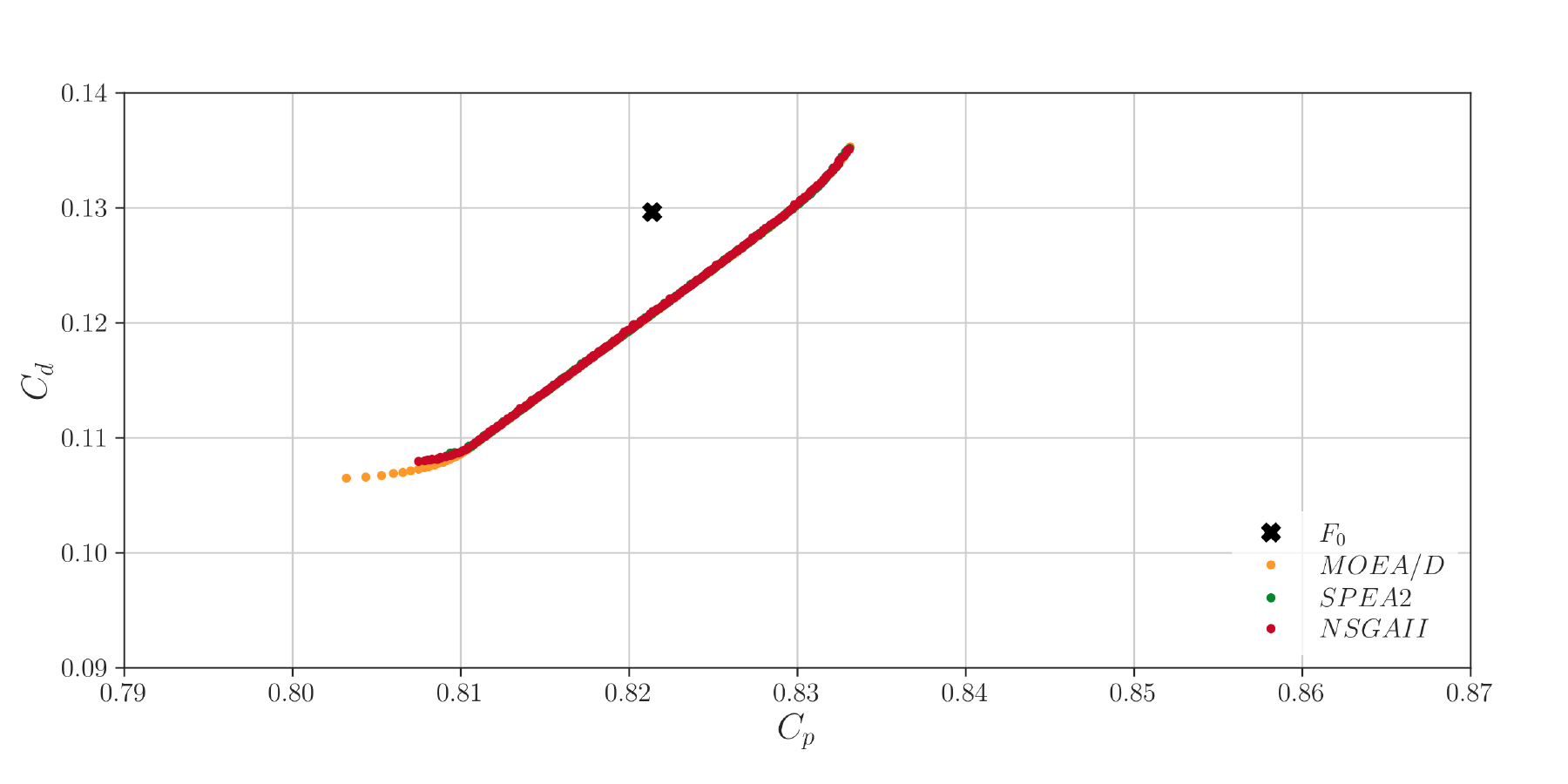}\vspace{-0.5cm}
    \caption{Pareto fronts for MOEA/D, SPEA2 and NSGA-II in test Scenario I.a. The fronts overlap to a large extent. However, as $C_p$ decreases, MOEA/D can provide additional Pareto-optimal solutions.}
    \label{fig:pf_scenarioIa}
\end{figure}

\begin{figure}[H]
    \centering
    \includegraphics[width=0.9\textwidth]{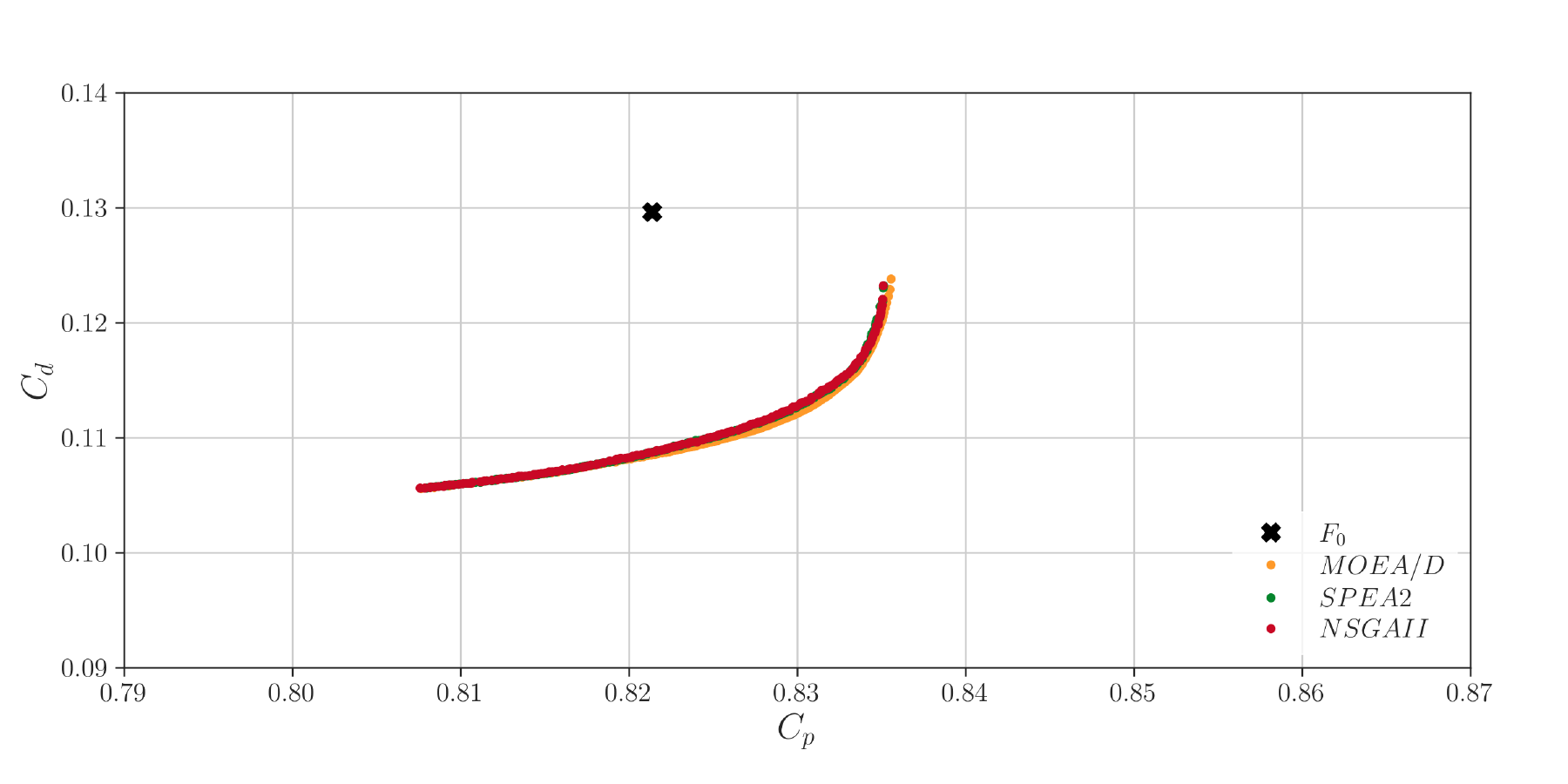}\vspace{-0.5cm}
    \caption{Pareto fronts for MOEA/D, SPEA2 and NSGA-II in test Scenario II.a. The fronts overlap to a large extent. As the width of the draft tube can change \protect\addedRIII{(introduced design flexibility)}, the fronts have moved downwards and to the right. The front for MOEA/D protrudes slightly to the right with the increase in $C_d$.}
    \label{fig:pf_scenarioIIa}
\end{figure}

\begin{figure}[H]
    \centering
    \includegraphics[width=0.9\textwidth]{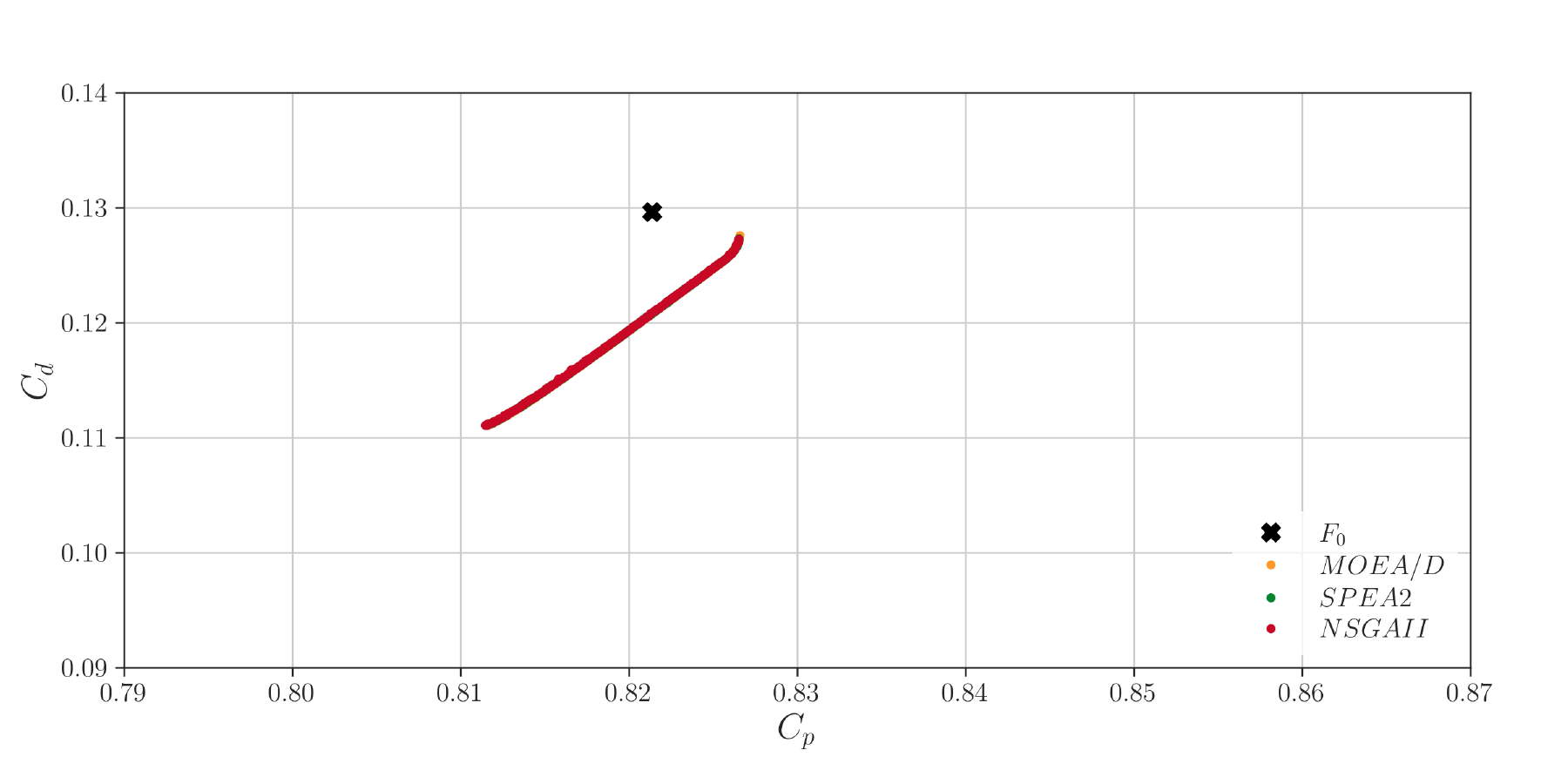}\vspace{-0.5cm}
    \caption{Pareto fronts for MOEA/D, SPEA2 and NSGA-II in test Scenario I.b. The fronts coincide. The fronts have moved upwards and to the left when compared to Scenario I.a. \protect\addedRIII{The movement of the fronts can be directly linked to the constrained design space, i.e. the limits set for the optimisation variables allow for designs with poorer overall performance.}}
    \label{fig:pf_scenarioIb}
\end{figure}

\begin{figure}[H]
    \centering
    \includegraphics[width=0.9\textwidth]{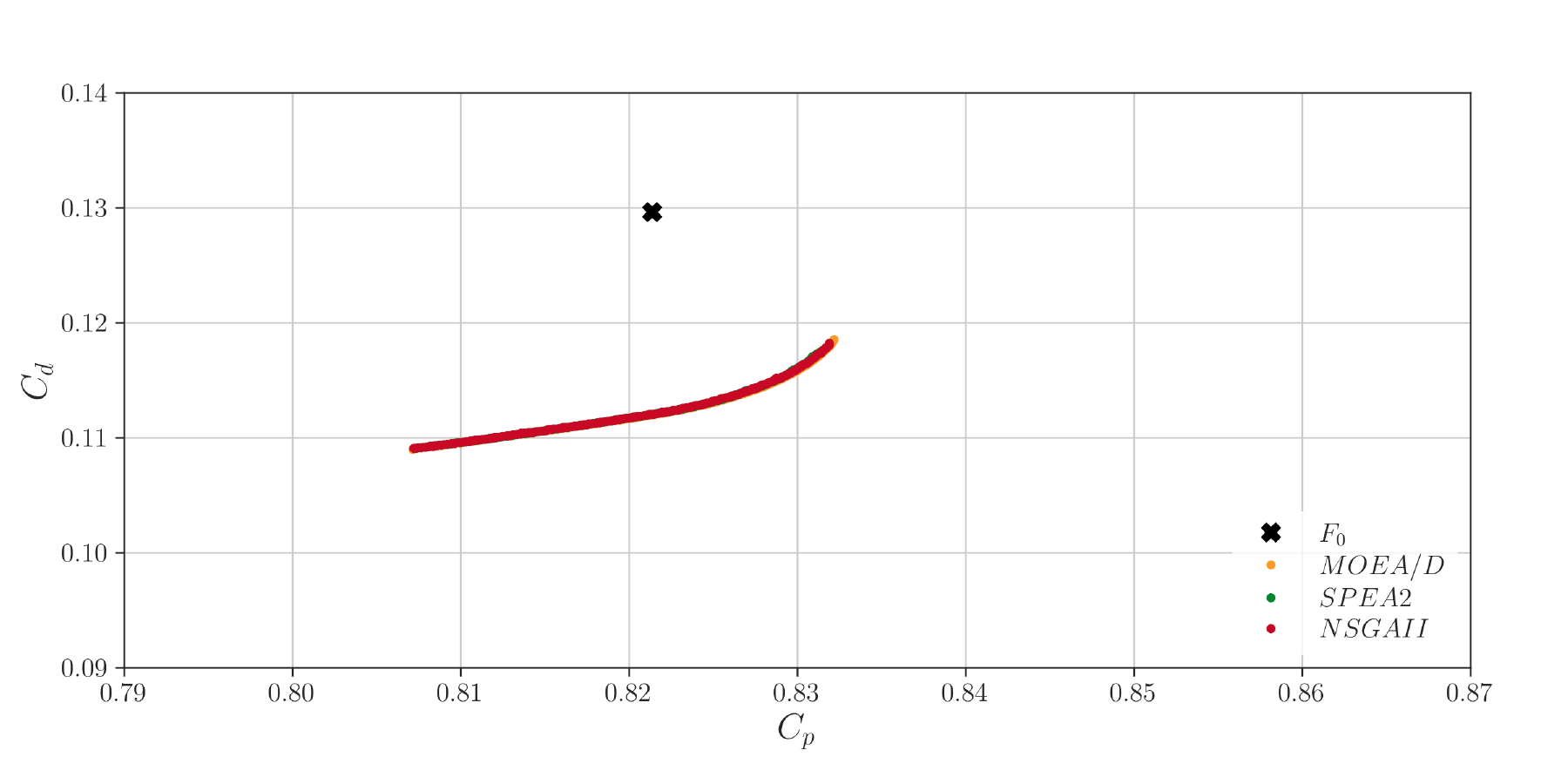}\vspace{-0.5cm}
    \caption{Pareto fronts for MOEA/D, SPEA2 and NSGA-II in test Scenario II.b. The fronts coincide. The fronts have moved upwards and to the left when compared to Scenario I.b. \protect\addedRIII{Analogously to Figure \ref{fig:pf_scenarioIb}, fixed limits restrict the flexibility of the design and thus the performance.}}
    \label{fig:pf_scenarioIIb}
\end{figure}

\begin{figure}[H]
\centering
\begin{subfigure}[b]{0.495\textwidth}
    \centering
    \includegraphics[trim={0cm 0.5cm 0cm 1cm}, clip, width=\textwidth]{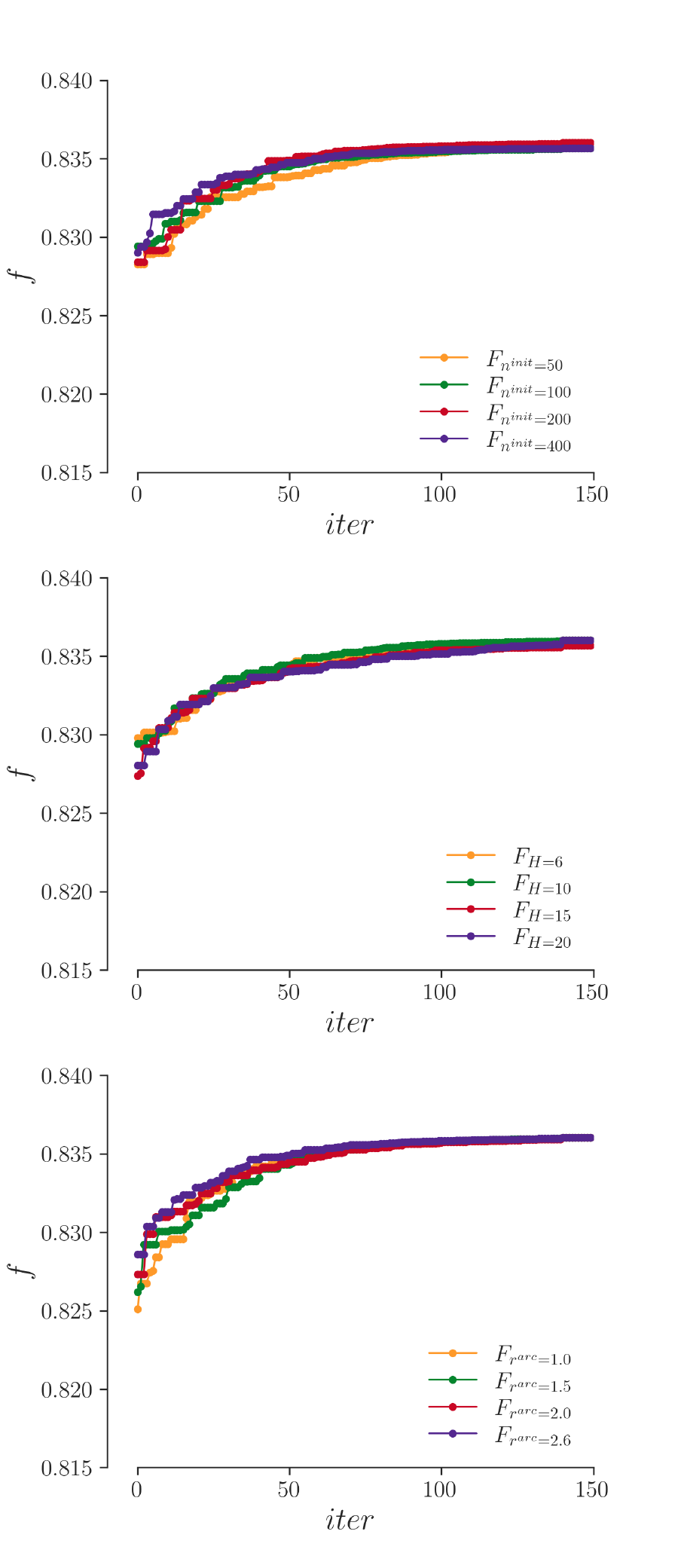}
    \caption{}
    \label{fig:comparisson_LSHADE:CP}
\end{subfigure}
\begin{subfigure}[b]{0.495\textwidth}
    \centering
    \includegraphics[trim={0cm 0.5cm 0cm 1cm}, clip, width=\textwidth]{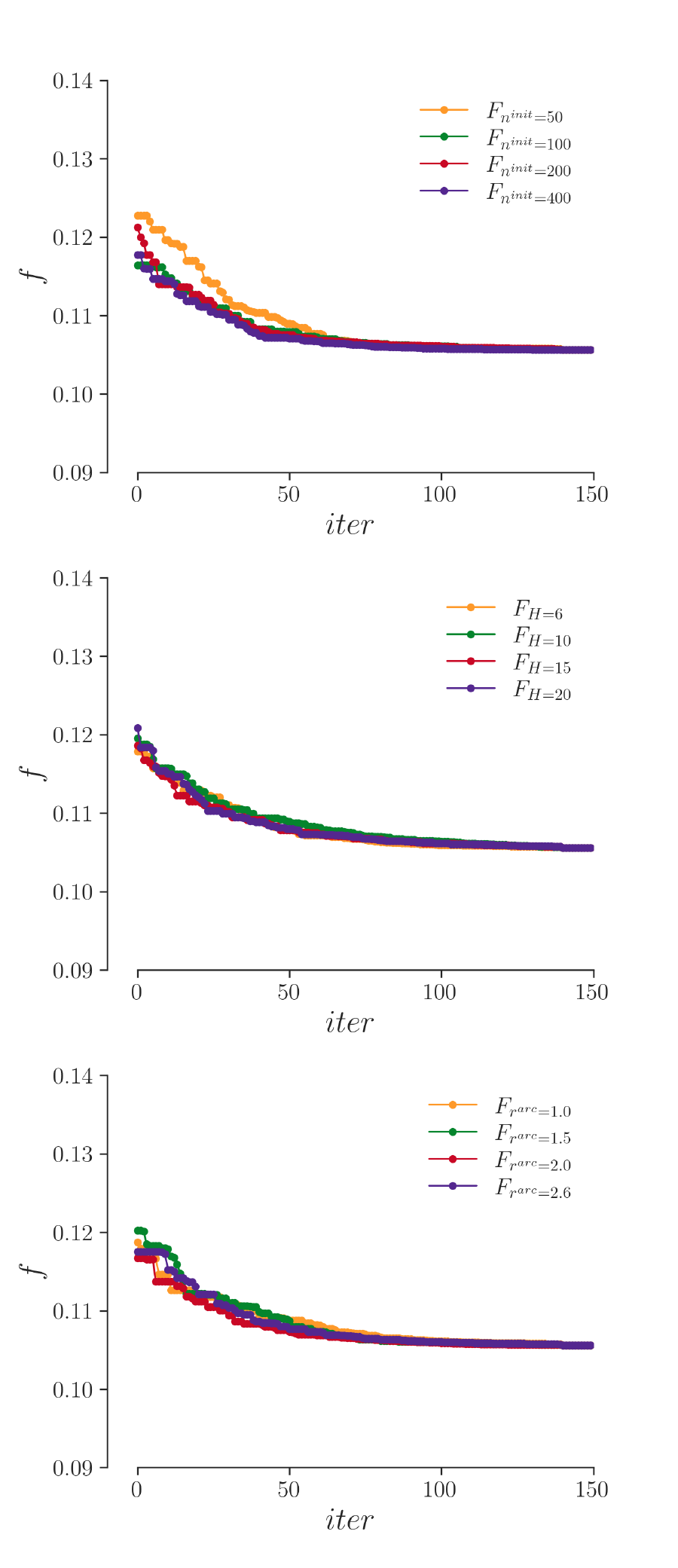}
    \caption{}
    \label{fig:comparisson_LSHADE:CD}
\end{subfigure}
\caption{Convergence graphs for L-SHADE when optimising for $C_p$ (a) and $C_d$ (b) in test Scenario II.a. Various parameters were evaluated, and their influence on convergence was assessed. The results show that $n^{init}$ has the greatest influence on convergence.}
\label{fig:comparisson_LSHADE}
\end{figure}

\begin{figure}[H]
    \centering
    \includegraphics[width=\textwidth]{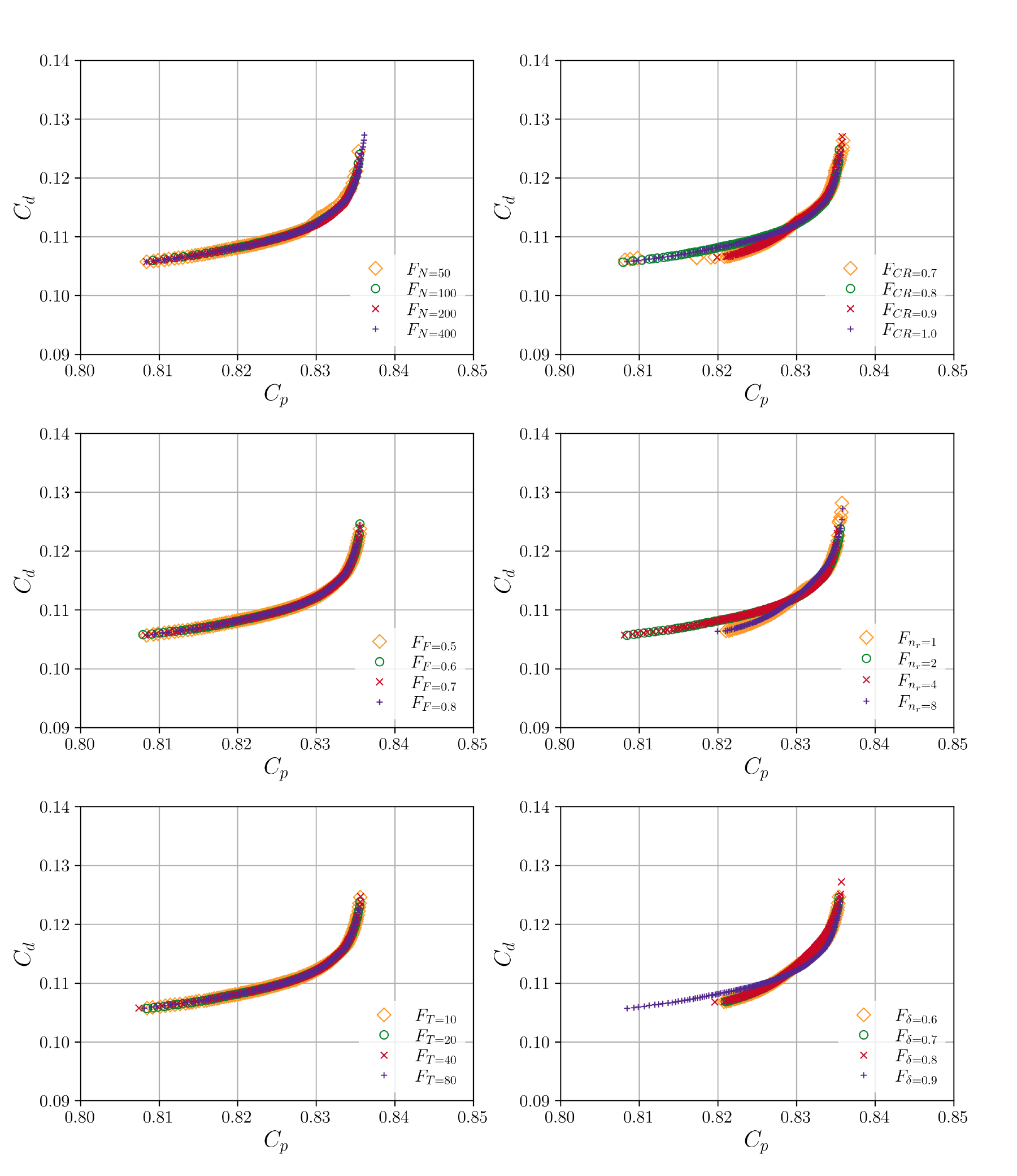}\vspace{-0.5cm}
    \caption{Pareto fronts for MOEA/D in test Scenario II.a. The influence of various \protect\addedRIII{algorithmic} parameters is examined.}
    \label{fig:comparisson_MOEAD}
\end{figure}

\begin{figure}[H]
    \centering
    \includegraphics[width=0.9\textwidth]{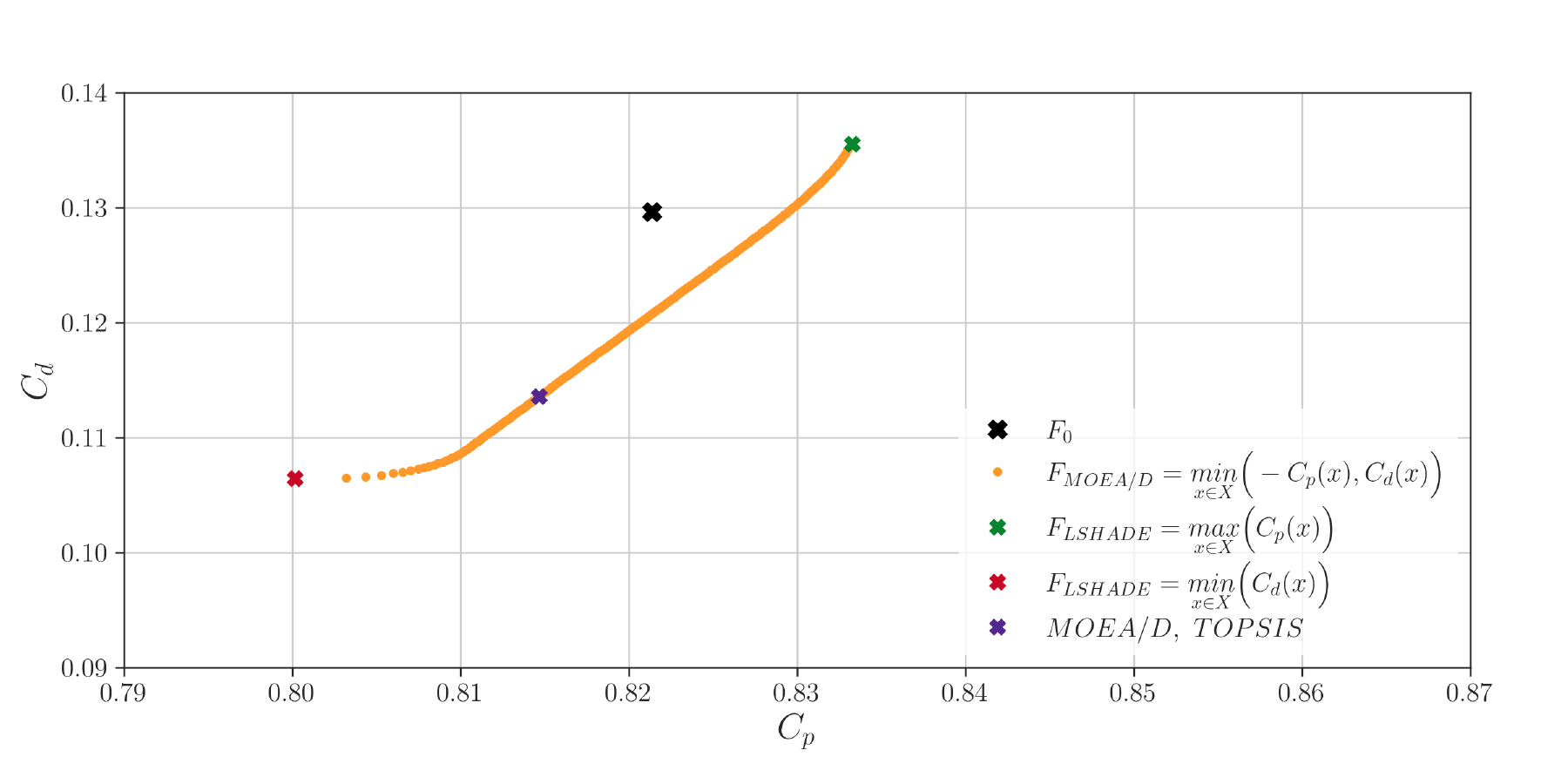}\vspace{-0.5cm}
    \caption{A comparison between SO and MO results for Scenario I.a. The Pareto front obtained using MOEA/D is compared to the SO results for $C_p$ and $C_d$ obtained using L-SHADE. The optimal MO result is determined using the TOPSIS approach. The front and SO optima are congruent.}
    \label{fig:compare_scenarioIa}
\end{figure}

\begin{figure}[H]
    \centering
    \includegraphics[width=0.9\textwidth]{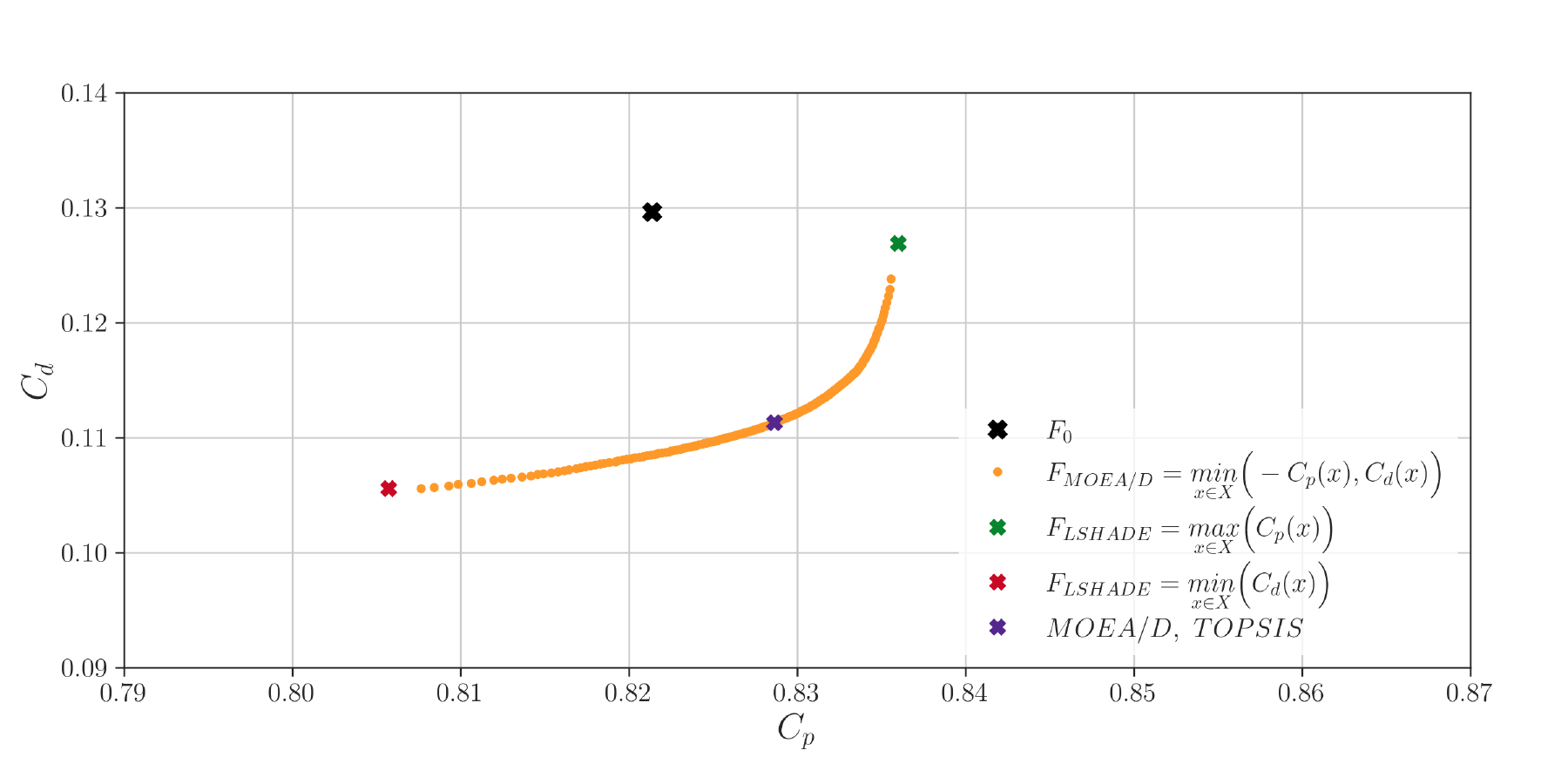}\vspace{-0.5cm}
    \caption{A comparison between SO and MO results for Scenario II.a. The Pareto front obtained using MOEA/D is compared to the SO results for $C_p$ and $C_d$ obtained using L-SHADE. The optimal MO result is determined using the TOPSIS approach.}
    \label{fig:compare_scenarioIIa}
\end{figure}

\begin{figure}[H]
    \centering
    \includegraphics[width=0.9\textwidth]{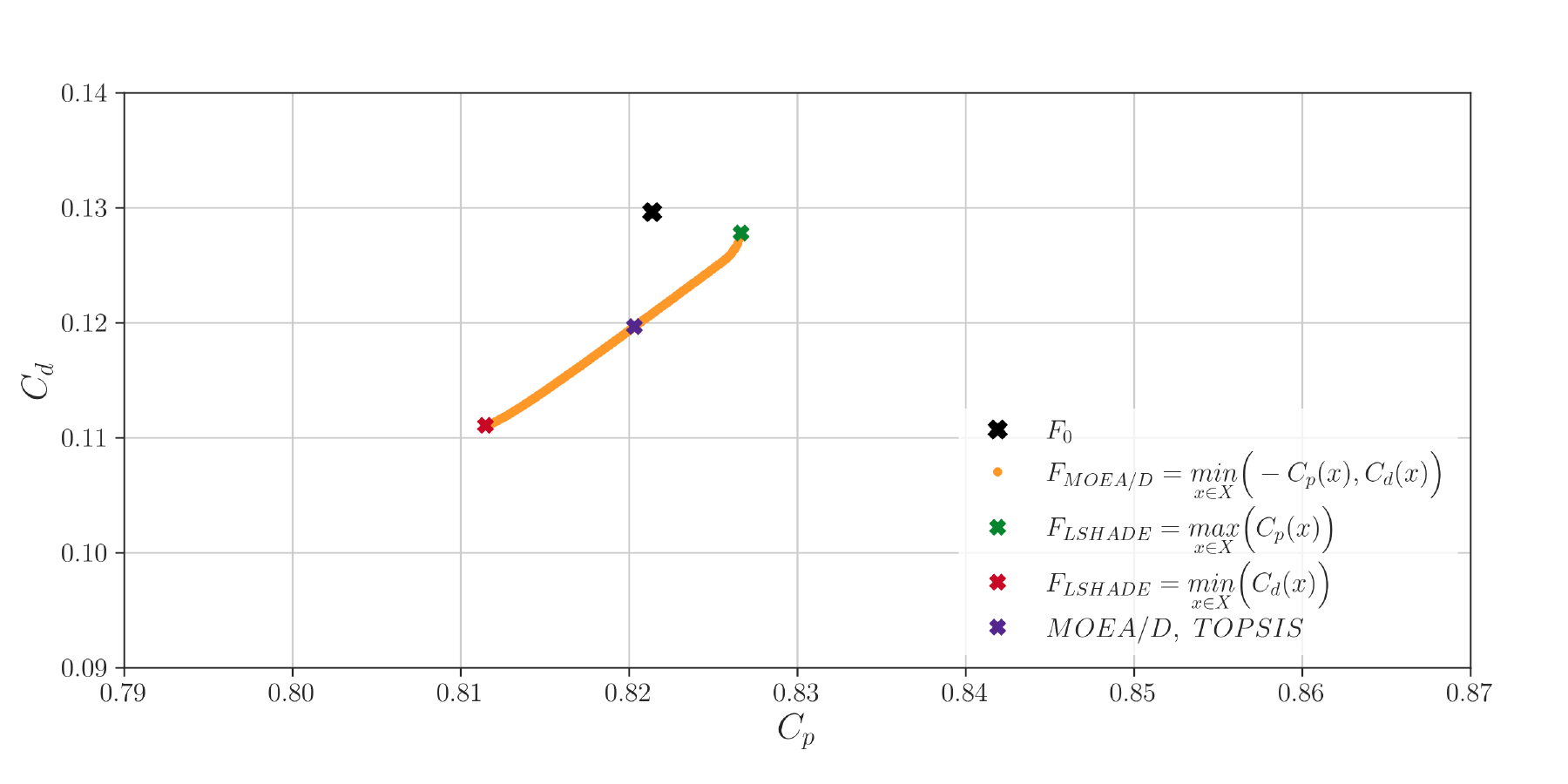}\vspace{-0.5cm}
    \caption{A comparison between SO and MO results for Scenario I.b. The Pareto front obtained using MOEA/D is compared to the SO results for $C_p$ and $C_d$ obtained using L-SHADE. The optimal MO result is determined using the TOPSIS approach.}
    \label{fig:compare_scenarioIb}
\end{figure}

\begin{figure}[H]
    \centering
    \includegraphics[width=0.9\textwidth]{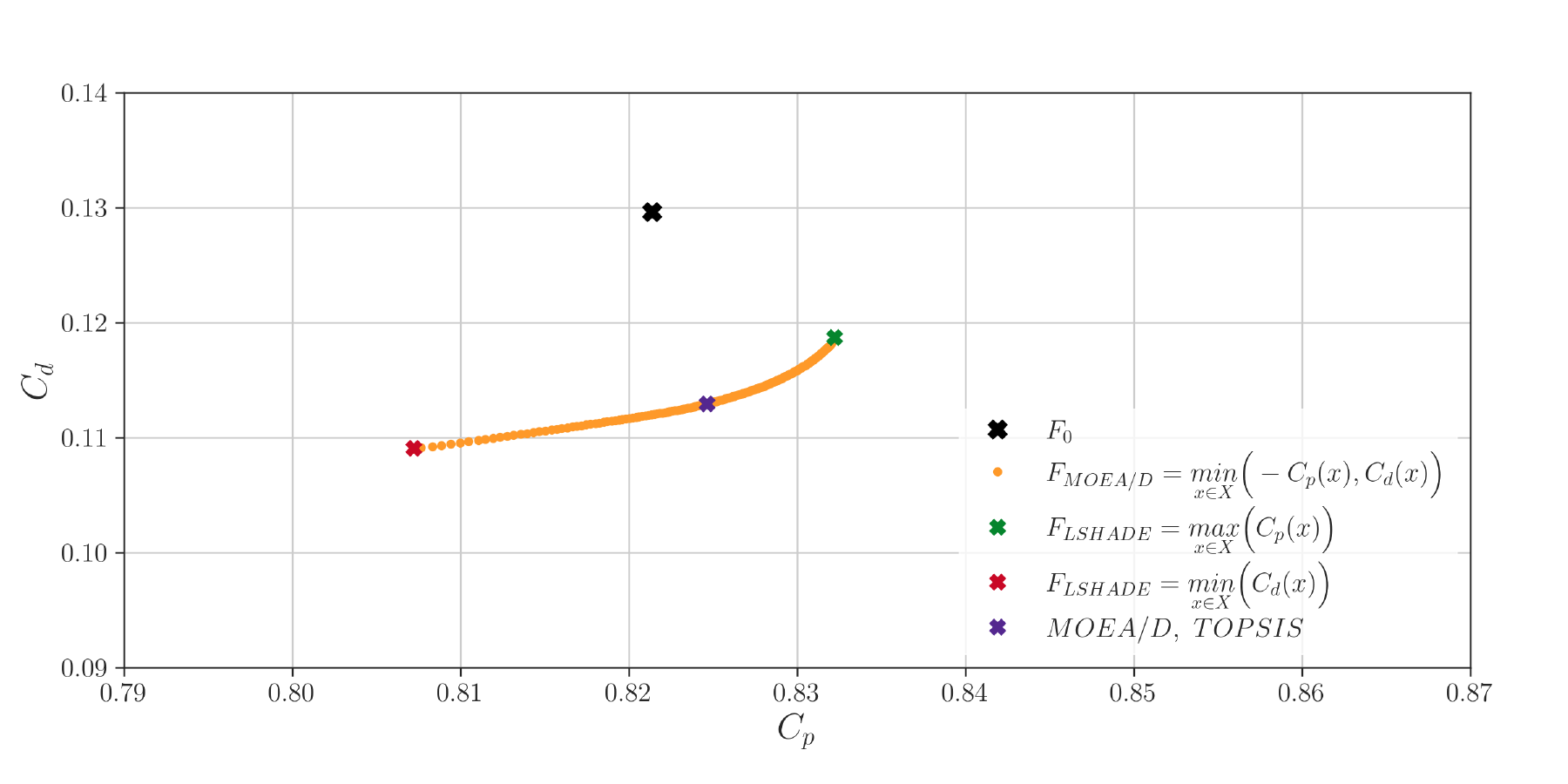}\vspace{-0.5cm}
    \caption{A comparison between SO and MO results for Scenario II.b. The Pareto front obtained using MOEA/D is compared to the SO results for $C_p$ and $C_d$ obtained using L-SHADE. The optimal MO result is determined using the TOPSIS approach.}
    \label{fig:compare_scenarioIIb}
\end{figure}

\bibliographystyle{elsarticle-num-names} 
\bibliography{bibliography}

\end{document}